\documentclass[11pt]{article}
\usepackage{amsthm}
\usepackage{amssymb}
\usepackage{amsmath}
\usepackage{amsfonts}
\usepackage{amscd}
\usepackage{graphicx}
\usepackage{microtype}
\usepackage{enumerate}
\usepackage[all,cmtip]{xy}
\usepackage[left=1in,right=1in,top=1.1in,bottom=1in]{geometry}

\newdir^{ (}{{}*!/-5pt/@^{(}}
\entrymodifiers={+!!<0pt,\fontdimen22\textfont2>}

\theoremstyle{plain}

\newtheorem{theorem}{Theorem}[section]
\newtheorem{question}[theorem]{Question}
\newtheorem{proposition}[theorem]{Proposition}
\newtheorem{lemma}[theorem]{Lemma}
\newtheorem{corollary}[theorem]{Corollary}
\newtheorem*{proposition:Tbracket}{Proposition~\ref{prop:Tbracket}}
\newtheorem*{theorem:nonsep}{Theorem~\ref{thm:intrononsep}}
\newtheorem*{theorem:sep}{Theorem~\ref{thm:introseparating}}

\theoremstyle{definition}
\newtheorem{definition}[theorem]{Definition}
\newtheorem{remark}[theorem]{Remark}

\newcommand{\iso}{\cong}
\renewcommand{\phi}{\varphi}
\renewcommand{\epsilon}{\varepsilon}

\DeclareMathOperator{\im}{im}

\DeclareMathOperator{\Mod}{Mod}

\DeclareMathOperator{\Sp}{Sp}
\DeclareMathOperator{\Hom}{Hom}

\DeclareMathOperator{\Aut}{Aut}

\DeclareMathOperator{\Stab}{Stab}
\DeclareMathOperator{\Grp}{Groups}

\DeclareMathOperator{\id}{id}

\DeclareMathOperator{\coker}{coker}

\newcommand{\Z}{\mathbb{Z}}

\newcommand{\Q}{\mathbb{Q}}
\newcommand{\I}{\mathcal{I}}
\newcommand{\K}{\mathcal{K}}
\newcommand{\cS}{\mathcal{S}}
\renewcommand{\L}{\mathcal{L}}

\renewcommand{\S}{\Sigma}
\renewcommand{\P}{\mathcal{P}}
\newcommand{\Shat}{\widehat{\S}}
\newcommand{\Surhat}{\widehat{S}}
\newcommand{\Phat}{\widehat{\P}}
\newcommand{\ihat}{\widehat{i}}
\newcommand{\Sbar}{\overline{\S}}
\newcommand{\tS}{\S'}
\newcommand{\tP}{P'}
\newcommand{\tPP}{\P'}

\newcommand{\IS}{\mathcal{I}(\S)}
\newcommand{\KS}{\mathcal{K}(\S)}
\newcommand{\HS}{H(\S)}
\newcommand{\NS}{N(\S)}
\newcommand{\ES}{E(\S)}
\newcommand{\DS}{D(\S)}

\newcommand{\MS}[1]{\Mod_{(#1)}(\S)}

\newcommand{\Surf}{\mathcal{T}\text{Surf}}
\newcommand{\pihat}{\widehat{\pi}}

\newcommand{\wedgetwo}{{\textstyle \bigwedge^2}}

\newcommand{\inj}{\hookrightarrow}
\newcommand{\phibar}{\overline{\phi}}
\newcommand{\ISs}{\I(S)}
\newcommand{\KSs}{\K(S)}

\newcommand{\coloneq}{\mathrel{\mathop:}\mkern-1.2mu=}

\newcommand{\abs}[1]{\lvert#1\rvert}
\newcommand{\bwedge}{\textstyle{\bigwedge}}

\newcommand{\para}[1]{\medskip\noindent\textbf{#1.}\ \ }
\newcommand{\XXXX}[1]{}

\newcommand{\reffigurearcs}{1}
\newcommand{\reffigurearcsmoved}{2}

\author{Thomas Church  \thanks{The author
    gratefully acknowledges support from the National Science
    Foundation.}}
\title{Orbits of curves under the Johnson kernel}

\begin{document}
\maketitle
\begin{abstract}
This paper has two main goals. First, we give a complete, explicit, and computable solution to the problem of when two simple closed curves on a surface are equivalent under the Johnson kernel. Second, we show that the Johnson filtration and the Johnson homomorphism can be defined intrinsically on subsurfaces and prove that both are functorial under inclusions of subsurfaces. The key point is that the latter reduces the former to a finite computation, which can be carried out by hand. In particular this solves the conjugacy problem in the Johnson kernel for separating twists. Using a theorem of Putman, we compute the first Betti number of the Torelli group of a subsurface. 
\end{abstract}

\maketitle
\tableofcontents

\section{Introduction}
Let $S=S_{g,1}$ be a surface of genus $g$ with one boundary component, with basepoint $\ast\in \partial S$. The \emph{mapping class group} $\Mod(S)$ is the group of self-homeomorphisms of $S$ fixing $\partial S$, up to isotopy fixing $\partial S$. The mapping class group is filtered by the \emph{Johnson filtration} $\Mod_{(k)}(S)$, consisting of those mapping classes that act trivially on the universal $(k-1)$-step nilpotent quotient of $\pi_1(S,\ast)$. Of particular interest are the \emph{Torelli group} $\ISs=\Mod_{(2)}(S)$ and the \emph{Johnson kernel} $\Mod_{(3)}(S)$. Johnson~\cite{JohnsonII} proved that $\Mod_{(3)}(S)$ is equal to the subgroup $\KSs$  of $\Mod(S)$ generated by Dehn twists about separating curves.

The mapping class group $\Mod(S)$ acts on the set of all simple closed curves on $S$ (more precisely, their isotopy classes), and we say that two curves are \emph{equivalent} under a subgroup $\Gamma<\Mod(S)$ if they lie in the same $\Gamma$--orbit under this action. Two curves $C$ and $D$ are equivalent under $\Mod(S)$ if and only if the complements $S-C$ and $S-D$ are homeomorphic. One of the main goals of this paper is to describe precisely two simple closed curves are equivalent under $\KSs$; in other words, we determine when one simple closed curve can be taken to another by applying a sequence of separating Dehn twists. 

\para{Orbits of nonseparating curves} Our first theorem describes when two \emph{nonseparating} curves are equivalent under $\KSs$. Any two nonseparating curves are equivalent under $\Mod(S)$, even when considered as oriented curves. Johnson~\cite{JohnsonConjugacy} proved the deeper result that two oriented nonseparating curves are equivalent under $\ISs$ if and only if they are homologous.

It is easy to show that if $C$ and $D$ lie in the same $\KSs$-orbit, then the mapping class $T_CT_D^{-1}$ lies in $\KSs$. 
Our first main result shows that this condition is also sufficient.
We also obtain an alternate condition in terms of based loops $\gamma,\delta\in \pi_1(S,\ast)$ representing the curves $C$ and $D$.
In the following theorem, $\Gamma_k(S)$ denotes the $k$-th term of the lower central series of $\pi_1(S,\ast)$, indexed so that $\pi_1(S,\ast)=\Gamma_1(S)$ and $\Gamma_2(S)$ is its commutator subgroup. We use the well-known isomorphism $\Gamma_2(S)/\Gamma_3(S)\iso \bwedge^2 H_1(S)$ defined by $[\xi,\xi']\mapsto [\xi]\wedge [\xi']$, and given $a\in H_1(S)$ we denote by $a\wedge H_1(S)$ the subspace spanned by elements of the form $a\wedge y$.
\begin{theorem}[$\KSs$--orbits of nonseparating curves]
\label{thm:intrononsep}
Let $C$ and $D$ be oriented nonseparating curves homologous to $a\in H_1(S)$. The following are equivalent:
\begin{enumerate}[1.]
\item The nonseparating curves $C$ and $D$ are equivalent under $\KSs$.
\item $T_CT_D^{-1}\in \KSs$.
\item For \textbf{some} representatives $\gamma,\,\delta\in\pi_1(S,\ast)$ of the curves $C$ and $D$, the class $[\gamma\delta^{-1}]\in \Gamma_2(S)/\Gamma_3(S)\iso \bwedge^2 H_1(S)$ lies in the subspace $a\wedge H_1(S)$.
\item For \textbf{any} representatives $\gamma,\,\delta$ of $C$ and $D$,  $[\gamma\delta^{-1}]\in \Gamma_2(S)/\Gamma_3(S)\iso \bwedge^2 H_1(S)$ lies in $a\wedge H_1(S)$.
\end{enumerate}
\end{theorem}

Let us apply the theorem to the case when the nonseparating curves  $C$ and $D$  are disjoint and homologous, forming a so-called ``bounding pair''. In this case Johnson proved in \cite[Lemma~4B]{JohnsonAbelian} that $T_CT_D^{-1}\not\in \KSs$, so Condition 2 implies that $C$ and $D$ are not equivalent under $\KSs$; this was previously proved by Farb--Leininger--Margalit~\cite[Proposition~3.2]{FarbLeiningerMargalit}.

As an illustration, we show how Condition 4 would be applied in this case. The union $C\cup D$ necessarily separates $S$ into two components, say $S_{k,2}$ and $S_{g-k-1,3}$, and there exists a standard basis $\{\alpha_1,\beta_1,\ldots,\alpha_g,\beta_g\}$ for $\pi_1(S,\ast)$ so that $[\alpha_1,\beta_1]\cdots[\alpha_k,\beta_k]\alpha_{k+1}$ and  $\alpha_{k+1}$ represent the curves $C$ and $D$. If $\{a_i,b_i\}$ is the induced symplectic basis for $H_1(S)$ we have $[\gamma\delta^{-1}]=a_1\wedge b_1+\cdots+a_k\wedge b_k\in \Gamma_2(S)/\Gamma_3(S)\iso \bwedge^2 H_1(S)$. Since this element certainly does not lie in the subspace $a_{k+1}\wedge H_1(S)$, Condition 4 of Theorem~\ref{thm:intrononsep} is verified, giving another proof that disjoint nonseparating curves $C$ and $D$ never lie in the same $\KSs$--orbit.

Conditions 3 and 4 are most useful in practice, since to check whether $T_CT_D^{-1}\in \KSs$ requires either computing $\tau(T_CT_D^{-1})$ by factoring $T_CT_D^{-1}$ as a product of bounding pair maps, or calculating its action on an entire basis for $\pi_1(S)$.

\para{Orbits of separating curves} Our next theorem describes when two \emph{separating} curves $C$ and $D$ are equivalent under $\KSs$. A separating curve $C$ separates $S-C$ into two components homeomorphic to $S_{k,1}$ and $S_{g-k,2}$ for some $1\leq k\leq g$, and the $\Mod(S)$--orbit of $C$ is determined by the genus $k$. Johnson~\cite[Theorem 1A]{JohnsonConjugacy} proved that the $\ISs$--orbit of $C$ is determined by the rank-$2k$ symplectic subspace of $H_1(S)$ spanned by homology classes supported on the subsurface $S_{k,1}$. 

Given such a symplectic subspace $V<H_1(S)$, we denote by $\omega_V\in \bwedge^2 V$ the restriction of the symplectic form $\omega$ to $V$.
There is a natural surjection $H_1(S)\otimes \bwedge^2 H_1(S)\twoheadrightarrow \Gamma_3(S)/\Gamma_4(S)$ defined by $x\otimes y\wedge z\mapsto [\widetilde{x},[\widetilde{y},\widetilde{z}]]$, where $\widetilde{x},\widetilde{y},\widetilde{z}\in \pi_1(S,\ast)$ are any representatives of $x,y,z\in H_1(S)$.
We denote by $H_1(S)\otimes \omega_V$ the subspace of $\Gamma_3(S)/\Gamma_4(S)$ spanned by the images of elements $\{x\otimes \omega_V\,|\, x\in H_1(S)\}$.

\begin{theorem}[$\KSs$--orbits of separating curves]
\label{thm:introseparating}
Let $C$ and $D$ be separating curves cutting off the same symplectic subspace $V< H_1(S)$. The following are equivalent:
\begin{enumerate}[1.]
\item The separating curves $C$ and $D$ are equivalent under $\KSs$.
\item The separating twists $T_C$ and $T_D$ are conjugate in $\KSs$.
\item For \textbf{some} representatives $\gamma,\,\delta\in\pi_1(S,\ast)$ of the curves $C$ and $D$, the class $[\gamma\delta^{-1}]\in \Gamma_3(S)/\Gamma_4(S)$ lies in the subspace $H_1(S)\otimes \omega_V$.
\item For \textbf{any} representatives $\gamma,\,\delta$ of $C$ and $D$, the class $[\gamma\delta^{-1}]\in \Gamma_3(S)/\Gamma_4(S)$ lies in $H_1(S)\otimes \omega_V$.
\end{enumerate}
\end{theorem}

\para{Defining the Johnson filtration for subsurfaces} 
A standard inductive technique in studying the mapping class group is to reduce to the stabilizers of curves, which amounts to studying the mapping class group $\Mod(S')$ of subsurfaces $S'\subset S$. However, this approach has not been available for the Johnson filtration: the problem is that the restriction of $\Mod_{(k)}(S)$ to a subsurface $S'$ is not intrinsic to $S'$ as an abstract surface, but gives different subgroups of $\Mod(S')$ depending on how  $S'$ is embedded into $S$.

In his thesis \cite{CuttingPasting}, Putman took the first step toward resolving this problem. He showed that the restriction of the Torelli group to a subsurface $S'\subset S$ becomes intrinsic after adding only a small amount of homological data. A \emph{partitioned surface} $\S$ is a surface with nonempty boundary, equipped with a partition of its boundary components. Any subsurface $S'\subset S$ determines a partitioned surface $\S$, where the partition records which of the boundary components of $S'$ become homologous in the larger surface $S$. Putman defined the Torelli group $\IS$ of a partitioned surface as the restriction of $\I(S)$, proved this is well-defined regardless of the embedding $\Sigma\subset S$, and used this to give natural inductive proofs for many key theorems on the Torelli group. An alternate approach, which we take in this paper, is to first define $H(\S)$, a modified version of $H_1(S')$ which serves as the ``first homology group of the partitioned surface $\S$''. It can be thought of as the first homology of the smallest closed surface into which $\S$ embeds; see Section~\ref{sec:Torelli} for details. The mapping class group $\Mod(S')$ then acts on $\HS$, and we define the Torelli group $\IS$ as the subgroup acting trivially on $\HS$.

Based on this evidence, one might expect that to describe how the Johnson filtration $\Mod_{(k)}(S)$ restricts to a subsurface $S'$, it would be necessary to record more and more nilpotent data describing the restriction of the lower central series $\Gamma_k(S)$ to the subsurface $S'$. In Section~\ref{sec:Johnsonfiltration} we prove the surprising result that \emph{no additional data is necessary} to define the Johnson filtration on subsurfaces. Given only the data of a partitioned surface $\S$, we define in Definition~\ref{def:Johnsonfiltration} the \emph{partitioned Johnson filtration} $\MS{k}$. The key property, proved in Theorem~\ref{thm:preserved}, is that $\MS{k}$ is \emph{natural} under inclusions: if $\S$ is a subsurface of a larger surface $S$, then $\MS{k}$ is precisely the subgroup of $\Mod(\S)$ that lies in $\Mod_{(k)}(S)$. This makes it possible to apply inductive arguments to any term of the Johnson filtration. As one example, we prove in Theorem~\ref{thm:BBM} a coherence result for $\KSs$--stabilizers of subsurfaces; this result has already been used in Bestvina--Bux--Margalit~\cite{BBM} to compute the cohomological dimension of $\K(S)$.

\para{Defining and computing the Johnson homomorphism for subsurfaces}
To prove Theorems~\ref{thm:intrononsep} and \ref{thm:introseparating}, it is not enough to understand the Johnson kernel $\K(S)$ for subsurfaces; we also need to understand how the Johnson homomorphism behaves when restricted to subsurfaces. The \emph{Johnson homomorphism} \[\tau\colon \I(S)\to \Hom\big(H_1(S),\bwedge^2 H_1(S)\big),\] defined by Johnson in \cite{JohnsonAbelian}, is constructed from the action on the universal 2-step nilpotent quotient $\pi_1(S)/\Gamma_3(S)$. In particular, the kernel of $\tau$ is  the Johnson kernel $\KSs$ by definition. Johnson proved that the image of $\tau$ is the subspace $\bwedge^3 H_1(S)$, giving a short exact sequence 
\[1\to \KSs\to \ISs\to \bwedge^3 H_1(S)\to 0.\]

The key advance that lets us prove Theorems~\ref{thm:intrononsep} and \ref{thm:introseparating} is an \emph{intrinsic} definition of the Johnson homomorphism for a partitioned surface, without necessarily embedding it into a larger surface.
For any partitioned surface $\S$, we define in Definition~\ref{def:tauS} the \emph{partitioned Johnson homomorphism} \[\tau_\S\colon \IS\to \Hom\big(\HS,\NS\big).\] As in Johnson's original paper, $\tau_\S$ is defined from the action of $\IS$ on a 2-step nilpotent quotient of $\pi_1$, but we replace the lower central series of $\pi_1(S,\ast)$ by a variant depending on the partitioned surface $\S$. In particular, the abelian group $N(\S)$ is a modification of $\Gamma_2(S)/\Gamma_3(S)\cong \bwedge^2 H_1(S)$, just as $H(\S)$ is a modification of $\Gamma_1(S)/\Gamma_2(S)\cong H_1(S)$. We prove in Corollary~\ref{cor:sameasKS} that just as in the classical case, the kernel $\KS=\ker \tau_\S$ is the third term $\MS{3}$ of the partitioned Johnson filtration.

One of our main results on $\tau_\S$ is the exact computation of its image: we prove in Theorem~\ref{thm:WS} that $\im \tau_\S=W_\S$ for a certain explicitly defined subspace $W_\S<\Hom(\HS,\NS)$. If the components of $\partial \S$ are   partitioned into $b$ blocks, the image $W_\S$ can be identified (see \eqref{eq:WSses}) with $\bwedge^3 D(\S)^\perp\oplus (D(\S)^\perp)^{\oplus b-1}$, where $D(\S)^\perp$ is the image of $H_1(\S)$ in $\HS$. This gives a short exact sequence
\[1\to \KS\to \IS\to W_\S\to 0.\]

 In Theorem~\ref{thm:commutes} we prove that this partitioned Johnson homomorphism is natural under inclusions of subsurfaces, so for any mapping class supported on a subsurface, we can compute the Johnson homomorphism locally.
This reduces all of Johnson's classical computations to trivial or nearly-trivial computations.
For example, any separating twist $T_C$ is supported on an annulus $\S$. But if $\S$ is an annulus then $W_\S=0$ by definition, so $\tau_\S(T_C)=0$, and naturality then implies that $\tau_S(T_C)=0$ for any separating curve on any surface. Similarly, any bounding pair $T_CT_D^{-1}$ is supported on a pair of pants $\S$, so the computation of $\tau(T_CT_D^{-1})$ reduces to the computation for a pair of pants, for which $W_\S$ is just $\Z^2$.

The characterization of $\KSs$--orbits in Theorems~\ref{thm:intrononsep} and \ref{thm:introseparating} depends on computing the image under $\tau$ of $\Stab_{\I(S)}(C)$, which is closely related to computing $\im \tau_\S$ for the complementary components of $S-C$.
From the arguments in Section~\ref{sec:KSorbits} it will be clear how to compute this image, and thus the space of $\KSs$--orbits, for other configurations, such as arbitrary collections of  separating curves or nonseparating collections of nonseparating curves. However, there is no guarantee that the resulting classification can still be formulated in terms of $\gamma\delta^{-1}$ in these cases; that this is possible for a single separating curve seems to be a happy coincidence.

\para{First Betti number of $\IS$ and comparison with Putman} By combining the results of this paper with one of the main theorems of Putman~\cite{PutmanKS}, we prove the following theorem.
\begin{theorem}
\label{thm:betti}
Let $\S$ be a surface of genus $g\geq 3$ whose $n\geq 1$ boundary components are partitioned into $b$ blocks, and let $D=2g+n-b$. Then the first Betti number $b_1(\IS)=\dim H_1(\IS;\Q)$ is
\[b_1(\IS)=\binom{D}{3}+D^{b-1}.\]
Moreover, any finite-index subgroup of $\IS$ that contains $\KS$ also has $b_1=\binom{D}{3}+D^{b-1}$.
\end{theorem}
We deduce Theorem~\ref{thm:betti} from \cite[Theorem~1.2]{PutmanKS}, which states that whenever $\S$ has genus at least 3, the rational abelianization $H_1(\IS;\Q)$ is isomorphic to $\im \tau_\S\otimes \Q$, and moreover that this holds for any finite-index subgroup of $\IS$ containing $\KS$.  This means that the calculation of $\im \tau_\S$ in Theorem~\ref{thm:WS} is also a calculation of the first Betti number $b_1(\IS)$ of the Torelli group.

In \cite{PutmanKS}, Putman has independently addressed questions closely related to the focus of this paper, centered around the question of defining the Johnson kernel $\KS$ for a subsurface $\S$ of an ambient closed surface $S$. However, one key difference is that Putman does not prove that $\KS$ is well-defined, which forces him to always work relative to a fixed embedding into a closed surface. Fortunately, Theorem~\ref{thm:preserved} guarantees that our definition of $\KS$ agrees with Putman's definition, so Corollary~\ref{cor:sameasKS} tells us that $\KS$ is indeed well-defined.  In particular, \cite[Theorem~1.1]{PutmanKS} states that  whenever $\S$ has genus at least 2, $\KS$ is generated by separating twists (when $\S=S_{g,1}$, this gives a new proof of the main theorem in Johnson~\cite{JohnsonII}). 

\para{Orbits of curves under the Johnson filtration}
We conclude this introduction with a question that was posed to us by Dan Margalit, inspired by Theorem~\ref{thm:intrononsep}.
\begin{question}
\label{q:margalit}
Let $C$ and $D$ be nonseparating curves on $S$. Is it true that
\[C\text{ and }D\text{ are equivalent under }\Mod_{(k)}(S)\quad\iff\quad T_CT_D^{-1}\in \Mod_{(k)}(S)?\]
\end{question}
For $k=1$ this is trivial, for $k=2$ this was proved by Johnson, and for $k=3$ this is proved in Theorem~\ref{thm:intrononsep} as the equivalence of Conditions 1 and 2. For $k\geq 4$, although the methods of this paper do not suffice to answer Question~\ref{q:margalit}, they do allow us to reduce it to the following question. Let $\tau_k\colon \Mod_{(k)}(S)\to \Hom(H_1(S),\L_k(S))$ denote the $k$th higher Johnson homomorphism, and note that $\Sp(H_1(S))$ acts on the target of $\tau_k$ (see e.g.\ \cite[Section 2]{Morita} for details; these maps will not be used elsewhere in the paper).

\begin{question}
\label{q:tom}
Let $C$ be a nonseparating curve with homology class $c\in H_1(S)$, and let ${t_c\in \Sp(H_1(S))}$ be the symplectic transvection $x\mapsto x+\omega(c,x)c$. Is it the case that
\[\tau_k\big(\Stab_{\Mod_{(k)}(S)}(C)\big)= \ker(t_c-\id)\ \cap\ \tau_k(\Mod_{(k)}(S))?\]
\end{question}
Question~\ref{q:tom} is equivalent to Question~\ref{q:margalit}, as can be shown along the same lines as the proof of Theorem~\ref{thm:intrononsep} in Section~\ref{sec:Kononsep}. Note that the image $\tau_k\big(\Stab_{\Mod_{(k)}(S)}(C)\big)$ is always contained in $\ker(t_c-\id)$, since any element stabilizing $C$  commutes with $T_C$. Therefore the question is whether the right side of the equation is contained in the left side. The difficulty in answering Question~\ref{q:tom} in general  is that for $k\geq 4$, although many partial results have been obtained, we still do not know the image $\tau_k(\Mod_{(k)}(S))$.

\para{Acknowledgements} This paper is based on the author's 2011 Ph.D.{} thesis at the University of Chicago.
I am deeply grateful to my advisor Benson Farb for his constant encouragement, his guidance, and his unwavering support. I thank Vladimir Drinfeld, Aaron Marcus,  and Andy Putman for many helpful conversations, and I thank Benson Farb and Andy Putman for suggesting this problem. I thank Dan Margalit for pointing out that Condition 2 of Theorem~\ref{thm:intrononsep} could be added, and for conversations regarding Question~\ref{q:margalit}. I am extremely grateful to the anonymous referee for their careful reading and  helpful suggestions, which greatly improved the organization and exposition of the paper.

\section{Background}
\subsection{Partitioned surfaces}
\label{section:partitionedsurfaces} Let $S$ be a compact connected
surface with nonempty boundary, with a partition $\P$ of its set of
boundary components $\pi_0(\partial S)$, and a basepoint
$\ast\in \partial S$; we call $\S=(S,\P,\ast)$ a \emph{partitioned
  surface}. This notion was first used by Putman in \cite{CuttingPasting}.
We refer to the elements $P\in \P$ as \emph{blocks} of the partition $\P$; each block is a subset of $\pi_0(\partial S)$.
We distinguish the block $P_0\in \P$
which contains the component containing $\ast$.

\para{Basic terminology} The metaphor underlying all our terminology regarding partitioned surfaces is that $S$ is thought of as being
embedded into a larger surface $S'$, and the partition $\P$ records which components of $\partial S$
can be connected by a path in the complement $S'\setminus S$. (Here and throughout the paper, by $S'\setminus S$ we mean the complement in $S'$ of the interior of the subsurface $S$, so that $S'\setminus S$ is itself a compact surface with boundary.) The data of $\S$ allows us to work intrinsically on $S$, without needing to embed it in a larger surface, or to choose between different embeddings.

We say that two
boundary components are \emph{connected outside $\S$} if they lie in
the same block $P\in\P$. 
A separating curve $\gamma$ on $S$ is called \emph{$\P$--separating}
(or just \emph{separating} if the partition is clear from context) if
each block $P\in \P$ of boundary components lies entirely on one side
or the other of $\gamma$. We say that a boundary component $z$ is
\emph{separating} if $\{z\}\in \P$, and that the partition $\P$ is
\emph{totally separated} if each boundary component is separating.  If 
$\P$ consists of a single block (and $|\pi_0(\partial S)|>1$), we say that $\S$ is
\emph{nonseparating}, since in this case no curve on $S$ which separates any boundary components can be
$\P$--separating.  

\para{Inclusions of partitioned surfaces} If $S$ is a subsurface of a surface $S'$, we say that a path lies
\emph{outside $S$} if it is contained in the complement $S'\setminus{S}$. If
$S'$ is a closed surface, $S$ inherits a partition of its boundary
components from $S'$ by defining two components of $\partial S$ to be
connected outside $\S$ if they can be connected by a path outside
$S$. More generally, if $S$ is a subsurface of a partitioned surface
$\S'=(S',\P',\ast')$, the subsurface $S$ inherits the structure of a partitioned surface
$\S=(S,\P,\ast)$ as follows. The partition $\P$ is defined by saying
that two components $z_1,z_2\in\pi_0(\partial S)$ are connnected
outside $\S$ (lie in the same block $P\in \P$) if either there is a path outside $S$ from $z_1$ to
$z_2$, or there exist components $z_1',z_2'\in\pi_0(\partial S')$ with
paths outside $S$ from $z_i$ to $z_i'$ and such that $z_1'$ and $z_2'$
are connected outside $\S'$ (they lie in the same block $P'\in \P'$). For the basepoint $\ast$ we choose any
point in $\partial S$ that can be connected to $\ast'$ by a path
outside $S$. Although the basepoint is not uniquely defined, the block $P_0\in \P$ containing it is, and for most purposes this is all that is relevant.

\subsection{The Torelli group}\label{sec:Torelli}
In this section, we define the homology of a partitioned surface $\S$,
which we denote by $\HS$; this originally appeared in Putman
\cite{CuttingPasting} using a different but equivalent
definition.

\para{The totally separated surface $\Shat$} Given a partitioned surface $\S$, we construct a totally
separated surface $\Shat=(\Surhat,\Phat,\widehat{\ast})$ with a
canonical embedding $\S\to \Shat$. For each block $P\in \P$ with
$\abs{P}=n$, we take a surface $S_{0,n+1}$ of genus 0 with $n+1$
boundary components, and glue all but one of these to the $n$ boundary
components in $P$. (Notice that when $n=1$ this operation is
effectively trivial.) The resulting surface $\Surhat$ has
$\abs{\pi_0(\partial\Surhat)}=\abs{\P}$; we take the partition $\Phat$ to be the totally
separated partition consisting of singleton blocks. For the basepoint
$\widehat{\ast}\in\partial \Surhat$ we choose any point so that $\ast$ and
$\widehat{\ast}$ lie in the same component of $\Surhat\setminus
{S}$.

The role of the surface $\Shat$ is captured by
the property that those components of $S$ that are connected outside
$\S$ are exactly those that are connected outside $S$ in $\Surhat$. An important consequence is that a curve $\gamma$ in $S$ is $\P$--separating if and only if $\gamma$ is a separating curve in $\Surhat$.
The embedding   $\S\to\Shat$ is universal, in that any embedding $\S\to \S'$ with $\S'$
totally separated factors through $\S\to\Shat$. As a consequence of this universal property, we see that this construction is idempotent: $\widehat{\Shat}=\Shat$.

\para{The homology $\HS$ of a partitioned surface}
The inclusion of $\partial\Surhat$
into $\Surhat$ gives a map from $H_1(\partial\Surhat)$ to
$H_1(\Surhat)$. (All homology groups in this paper are taken with integral coefficients, except in Remark~\ref{rem:rationalhomology} where we explicitly specify otherwise.) We define $\HS$ to be the cokernel of this map: \[\HS\coloneq \coker\big(
H_1(\partial \Surhat)\to H_1(\Surhat)\big).\]
A separating curve in $\Surhat$ is homologous to a collection of boundary components, and thus vanishes in $\HS$. Applying our characterization of $\P$--separating curves above, we conclude that a curve $\gamma$ in $S$ is $\P$--separating if and only if $[\gamma]=0\in \HS$.
The observation above that $\widehat{\Shat}=\Shat$ implies tautologically that $\HS=H(\Shat)$.

\para{The Torelli group $\IS$}
The mapping class group $\Mod(S)$ of $S$ is the group of self-homeomorphisms of $S$ fixing $\partial S$ pointwise, up to isotopy fixing $\partial S$ pointwise.
(We remark that throughout this paper, Dehn twists are twists to the \emph{right}.) 
Given any inclusion $i\colon S\hookrightarrow S'$ of surfaces, a homeomorphism $\phi\in \Mod(S)$ can be extended by the identity on the complement $S'\setminus S$ to obtain $i_*(\phi)\in \Mod(S')$. In particular, for any partitioned surface $\S=(S,\P,\ast)$, the natural inclusion $\S\to \Shat$ induces an embedding $\Mod(S)\to \Mod(\Surhat)$. Since $\Mod(\Surhat)$ naturally acts on $\HS$, composing with this embedding we obtain an action of $\Mod(S)$ on $\HS$.

We
define the \emph{Torelli group} of $\S$ as \[\IS\coloneq
\big\{\,\phi\in\Mod(S)\ \,\big|\ \,\phi\text{ acts trivially on
    }\HS\big\}.\] We obtain an exact sequence
\[1\to \IS\to \Mod(S)\to \Aut(\HS)\] but the latter map is not in
general surjective. (It is possible to show that the image is precisely the symplectic automorphisms preserving the homology classes of all boundary components of $\partial S$, but we will not need this here. For details, see the earlier version of this paper posted at arXiv:1108.4511v1.)

\para{An alternate definition of $\HS$} For future reference, we give another definition of $\HS$. Given a
partitioned surface $\S$, we define $\overline{\S}$ to be the a surface with
one boundary component obtained by gluing a disk to each boundary component of
$\Shat$ except the component containing the basepoint. (Equivalently, $\overline{\S}$ is obtained from $S$ by gluing an $S_{0,\abs{P_0}+1}$ to the boundary components in the block $P_0$ and an
$S_{0,\abs{P}}$ to each other block $P\in \P$ in $\partial S$.)
The Mayer--Vietoris sequence
implies that $\HS\simeq H_1(\overline{\S})$, so the action of $\Mod(S)$ on $\HS$ factors through the action of $\Mod(\overline{\S})$ on $H_1(\overline{\S})$. Since $\overline{\S}$ has only one boundary component, the intersection form on
$H_1(\overline{\S})$ is a $\Mod(\overline{\S})$--invariant symplectic form. In particular, this implies
that $\HS$ is self-dual as a $\Mod(S)$--module.

\subsection{The Torelli category}
\label{sec:Torellicategory} Putman defines a category whose objects
are partitioned surfaces and whose morphisms are inclusions of
subsurfaces respecting the partitions. For our purposes, we will need
the following refinement of this category.

\begin{definition}
Given two partitioned
surfaces $\S_1=(S_1,\P_1,\ast_1)$ and $\S_2=(S_2,\P_2,\ast_2)$ and an inclusion $i\colon S_1\inj S_2$ of their underlying surfaces, we
say that:
\begin{itemize}
\item $i$ \emph{respects the partitions} if $\P_1$--separating and
  $\P_1$--nonseparating curves are taken to $\P_2$--separating and
  $\P_2$--nonseparating curves respectively; and
\item $i$ \emph{preserves basepoints} if $\ast_1$ and $\ast_2$ lie in
  the same component of $S_2\setminus i({S_1})$.
\end{itemize}
\end{definition}
As we described in Section~\ref{section:partitionedsurfaces}, for any inclusion $i\colon S_1\inj S_2$ the subsurface $S_1$ inherits the structure of a
partitioned surface from $\S_2$. An inclusion satisfies these two properties\,---\,that is, it both respects the partitions and preserves basepoints\,---\,exactly when  the inherited structure on $S_1$ is $\S_1$.

\para{The Torelli category} The category $\Surf$ is defined as follows. Its objects are
partitioned surfaces $\S=(S,\P,\ast)$. A morphism $\iota\colon \S_1\to
\S_2$ from $\S_1=(S_1,\P_1,\ast_1)$ to $\S_2=(S_2,\P_2,\ast_2)$ is an inclusion $i\colon S_1\inj S_2$ of the underlying surfaces that respects the
partitions and preserves basepoints, together with an inclusion
$\ihat\colon \Surhat_1\inj \Surhat_2$ extending $i$. (If we liked, we could identify morphisms in $\Surf$ when the underlying inclusions are isotopic; for simplicity we elect not to do this, but everything in this paper would descend nicely to this quotient category.)

The canonical inclusion $S\inj \Surhat$ induces a morphism $\S\to
\Shat$ for any $\S$.
For any morphism $\iota\colon \S_1\to \S_2$, the inclusion $\ihat$ induces a map $H_1(\Shat_1)\to H_1(\Shat_2)$. The fact that $\iota$ respects the partitions implies that 
this descends to a map $\iota_*\colon H(\S_1)\to
H(\S_2)$.

If $\S_2$ is a partitioned surface, any inclusion $i\colon S_1\to S_2$
gives the subsurface the structure of a partitioned surface
$\S_1$. This inclusion always extends to a morphism $\iota\colon
\S_1\to \S_2$, but not canonically; the ambiguity is in the choice of
the map $\ihat\colon\Surhat_1\to \Surhat_2$, or equivalently in the
choice of the inclusion $\iota_*\colon H(\S_1)\to H(\S_2)$.

Given a morphism $\iota\colon \S_1\to \S_2$, extension by the identity
induces a map $\Mod(S_1)\to \Mod(S_2)$, which restricts to a map
$\iota_*\colon \I(\S_1)\to \I(\S_2)$. Putman showed in \cite{CuttingPasting} that the
Torelli group can be regarded as a functor $\I$ from $\Surf$ to the
category of groups and homomorphisms. Our
category $\Surf$ is actually a refinement of the category considered
by Putman; one key benefit of this refinement is that the assignment $\S\to \HS$ becomes functorial. Moreover, this lets us interpret the Johnson homomorphism as a natural transformation, as we will show in Theorem~\ref{thm:naturaltransformation}. 

\para{Non-collapsing inclusions and simple cappings}
When dealing with an inclusion of partitioned surfaces, it is especially convenient if the inclusion does not ``close off'' any block $P\in \P$. Formally, we make the following definition.
\begin{definition}
A morphism $\iota\colon \S_1\to \S_2$ is \emph{non-collapsing} if for each component $U$ of
  $S_2\setminus {S_1}$ we have $\partial U\not\subset\partial
  S_1$.

\end{definition}
In other words, every boundary component in $\partial S_1$ can be connected to $\partial S_2$ by an arc lying outside $S_1$. One convenient property of such inclusions is that if $\iota$ is non-collapsing, the map $\iota\colon \Mod(S_1)\to \Mod(S_2)$ is injective (see Section~\ref{section:Johnsonfiltpreserved}). Of course, not every morphism of partitioned surfaces is non-collapsing; the most basic examples of this are a class of morphisms that we will call ``simple cappings''.

\begin{definition}
A morphism $\iota\colon \S_1\to \S_2$ is a \emph{simple capping} if $S_2\setminus
  {S}_1$ is a single disk.
\end{definition}
Any inclusion can be factored as the composition of a single
non-collapsing inclusion with a sequence of simple cappings, so we can often reduce to considering these special cases separately. 
Note that since a simple capping respects the partitions, the boundary
component which is capped off must be separating.

\section{The lower central series on a subsurface}
\label{sec:charT}
When a subsurface $S$ is embedded in a surface
$S_0$ with one boundary component, restricting the lower central series of $\pi_1(S_0)$ to $\pi_1(S)$ yields a central filtration of $\pi_1(S)$. In this section we show that this filtration of $\pi_1(S)$ depends only on which boundary components of $S$ become homologous in $S_0$; that is, it can be intrinsically defined in terms of the partitioned surface $\Sigma=(S,\P,\ast)$.
One key consequence is that we can define the Johnson filtration for a partitioned surface, which we will show in Section~\ref{sec:Johnsonfiltration} using the results of this section. 

The main technical idea of this section is that if the associated graded Lie algebra of a central filtration on a group happens to be a free Lie algebra, then to describe the filtration it suffices to find a free basis. Moreover, if a purported
basis is known to generate the Lie algebra, we can verify that it is a free basis by mapping to a Lie algebra
already known to be free.

\para{The lower central series}
Given any group $\Gamma$, its lower central series is defined by $\Gamma_1=\Gamma$ and $\Gamma_j=[\Gamma_1,\Gamma_{j-1}]$. 
If we define $\L_j\coloneq \Gamma_j/\Gamma_{j+1}$, the fact that $[\Gamma_i,\Gamma_j]\subset \Gamma_{i+j}$ implies that the commutator bracket on
$\Gamma_1$ descends to a bilinear map $\L_i\otimes \L_j \to
\L_{i+j}$. This makes the associated graded algebra $\L\coloneq \bigoplus
\L_j$ into a graded Lie algebra. (All Lie algebras are over $\Z$ unless
otherwise specified. By a graded Lie algebra, we simply mean a Lie algebra endowed with a grading respected by the bracket; that is, we do not introduce any signs coming from the grading.) It is well-known that if $\Gamma$ is a free group
with basis $\{x_1,\ldots,x_n\}$, then $\L$ is the free Lie algebra on
the same generating set (Witt~\cite{WittD}).

\para{The central series $\Gamma^T_j(\S)$} Given a partitioned surface $\S=(S,\P,\ast)$, let $\pi\coloneq \pi_1(S,\ast)$. We define the normal subgroup $T(\S)$ to be the kernel of the composition $\pi_1(S,\ast)\to H_1(S)\to \HS$.

We define the central series $\Gamma_j^T=\Gamma_j^T(\Sigma)$ by
\[\Gamma_1^T=\pi,\qquad\qquad
\Gamma_2^T=T(\S),\qquad\qquad \Gamma_j^T=\langle
[\Gamma_1^T,\Gamma_{j-1}^T],[\Gamma_2^T,\Gamma_{j-2}^T]\rangle\text{
  for } j\geq 3.\] This is the minimal filtration satisfying
$\Gamma_2^T=T(\S)$ and $[\Gamma_i^T,\Gamma_j^T]\subset
\Gamma_{i+j}^T$.

\para{Explicit generators for $T(\S)$} It will be very useful to have explicit generators for $T(\S$). Let $k=\abs{\P}-1$.  For each block $P_i\in \P$, choose a $\P$--separating curve $\gamma_i$ in $S$ so that the boundary components lying on one side of $\gamma_i$ are exactly those lying in the block $P_i$, and choose $\zeta_i\in \pi$ representing $\gamma_i$. There are of course many such curves $\gamma_i$, and many representatives $\zeta_i$, but the following lemma tells us that any choice of such elements $\zeta_i$ provides generators for $T(\S)$.
\begin{lemma}
\label{lem:TSgen}
The normal subgroup $T(\S)$ is generated by $[\pi,\pi]$ together with the elements $\zeta_1,\ldots,\zeta_k$.
  \end{lemma}
\begin{proof}
  By definition, $\HS$ is the quotient of $H_1(\Surhat)$ by $H_1(\partial \Surhat)$. Each component of $\Surhat$ is a genus 0 homology between the $i$th component of $\partial \Surhat$ and the $|P_i|$ components of $\partial S$ lying in the block $P_i\in \P$. Let $m_i=\abs{P_i}-1$, and let $a_i^j$ be the homology classes in $H_1(S)$ of the boundary  components  in $P_i$ (for $0\leq j\leq m_i$).  The Mayer--Vietoris sequence implies that the image of $H_1(S)$ in $\HS$ is the quotient of $H_1(S)$ by the elements $a_i^0+a_i^1+\cdots+a_i^{m_i}$ for each $i$. But our assumption on $\gamma_i$ guarantees that $[\gamma_i]=a_i^0+a_i^1+\cdots+a_i^{m_i}\in H_1(S)$. It follows that the kernel of the map $H_1(S)\to \HS$ is generated by the homology classes $[\gamma_i]$ for $0\leq i\leq k$. Moreover, the fundamental class of the surface $S$ itself gives the relation $\sum_i,j a_i^j$, which can be rewritten as $[\gamma_0]+[\gamma_1]+\cdots+[\gamma_k]=0$. Thus $\ker(H_1(S)\to \HS)$ is in fact generated by $[\gamma_i]$ for $1\leq i\leq k$. It follows that the kernel of the composition $\pi_1(S)\to H_1(S)\to \HS$ is generated by $[\pi,\pi]$ together with elements $\zeta_i$ representing $\gamma_i$ for $1\leq i\leq k$.
\end{proof}

\begin{remark}
\label{rem:TSoneboundary}
Note that if $S$ is a surface with one boundary component, with $\S=(S,\{P_0\},\ast)$ the associated (trivial) partitioned surface, we have $S=\Surhat$ and so $\HS\iso H_1(S)$. It follows that the kernel $T(\S)$ of the map $\pi\to H_1(S)\iso \HS$ is just the commutator subgroup $[\pi,\pi]$, and so in this case the central series $\Gamma_j^T(\Sigma)$ is simply the lower central series $\Gamma_j(S)$ of the free group $\pi_1(S)$.
\end{remark}

\para{The graded Lie algebra $\L^T(\S)$} We set $\L^T_j(\S)\coloneq \Gamma_j^T(\S)/\Gamma_{j+1}^T(\S)$, and denote by
$\L^T(\S)=\bigoplus \L^T_j(\S)$ the associated graded Lie algebra. The fact that $\Gamma_j^T=\langle
[\Gamma_1^T,\Gamma_{j-1}^T],[\Gamma_2^T,\Gamma_{j-2}^T]\rangle$ for
$j\geq 3$ implies that $\L^T(\S)$ is generated by $\L^T_1(\S)$ and $\L^T_2(\S)$.We begin by constructing a generating set $\mathcal{S}(\S)$ for $\L^T(\S)$; we will eventually prove that $\mathcal{S}(\S)$ is a free basis for $\L^T(\S)$.

\para{The generating set $\cS(\S)$} We first construct a ``standard'' generating set for $\pi=\pi_1(S,\ast)$.
For $0\leq i\leq k$, choose a curve $\gamma_i$ cutting off $P_i$ as above, with the additional assumption that the subsurfaces cut off have genus 0, and that the curves $\gamma_i$ are mututally disjoint. Let $\zeta_i\in \pi_1(S,\ast)$ be a simple loop representing $\gamma_i$, oriented so that the genus 0 subsurface $R_i$ cut off by $\zeta_i$ lies on the left side of $\zeta_i$. Choose simple loops $\alpha_i^0,\alpha_i^1,\ldots,\alpha_i^{m_i}$, disjoint from $\zeta_i$ and from each other, so that $\alpha_i^j$ represents the $j$th boundary component in $P_i$. The elements $\{\alpha_i^0,\ldots,\alpha_i^{m_i}\}$ form a free basis for $\pi_1(R_i)$, and we may (uniquely) reorder these elements so that in $\pi$ we have the relation
  \begin{equation}\label{eq:zetaprime}\zeta_i=\alpha_i^0\alpha_i^1\cdots \alpha_i^{m_i}.\end{equation} 
Let $R_{\text{main}}$ denote the remaining component of $S-\cup \zeta_i$; it has the same genus $g$ as the original surface $S$. Choose simple loops 
  $\xi_1,\ldots,\xi_{2g}$ so that $\{\xi_1,\ldots,\xi_{2g},\zeta_1,\ldots,\zeta_k\}$ form a free basis for $\pi_1(R_{\text{main}})$, and so that in $\pi$ we have the relation
\begin{equation}
\label{eq:xizetaprime}
[\xi_1,\xi_2]\cdots[\xi_{2g-1},\xi_{2g}]\zeta_1\cdots\zeta_k\zeta_0=1.
\end{equation}
Applying Van Kampen's theorem, we conclude that a basis for the free group $\pi=\pi_1(S,\ast)$ is given by the set $\{\xi_1,\ldots,\xi_{2g}\}\cup\{\alpha_i^j\}$, excluding only the element $\alpha_0^0$.

Let $x_i$ and $a_i^j$ be the images of $\xi_i$ and $\alpha_i^j$ in $\L_1^T(\S)$, and let $z_i$ be
  the image of $\zeta_i$ in $\L_2^T(\S)$. 
\begin{proposition}
\label{prop:cS}
  $\L^T(\S)$ is generated by \begin{equation}
  \label{eq:cSdefinition}
  \mathcal{S}(\S)\coloneq \{x_1,\ldots,x_{2g}\}\cup \{a_i^j\}_{i\geq
    0}^{j\geq 1}\cup \{z_i\}_{i\geq 1}.
\end{equation}
\end{proposition}
\begin{proof} Since $\pi$ is generated by $\{\xi_i\}\cup \{\alpha_i^j\}$, the quotient $\L_1^T(\S)$ is spanned by $\{x_1,\ldots,x_{2g}\}$ together with $\{a_i^j\}^{j\geq 0}_{i\geq 0}$. From \eqref{eq:zetaprime} we obtain the
  relation \begin{equation}\label{eq:cancelai}a_i^0+a_i^1\cdots+a_i^{m_i}=0\quad\text{in
  }L_1^T(\S),\end{equation} which lets us eliminate the generator $a_i^0$. Lemma~\ref{lem:TSgen} shows that $T(\S)$ is generated by $[\pi,\pi]$ together with $\zeta_1,\ldots,\zeta_k$, so $\L_2^T(\S)$ is spanned by $[\L_1^T(\S),\L_1^T(\S)]$ together with $\{z_1,\ldots,z_k\}$. 
Since we observed above that $\L^T(\S)$ is generated by $\L^T_1(\S)$ and $\L^T_2(\S)$, this demonstrates that $\L^T(\S)$ is generated by
 $\cS(\S)$.
 \end{proof}

\para{Inclusions of partitioned surfaces and $\L^T(\S)$} Let $\S'=(S',\P',\ast')$ be another partitioned surface, and let $\pi'=\pi_1(S',\ast')$. Given a morphism $\iota\colon \S\to \S'$, the inclusion $i\colon S\hookrightarrow S'$ induces a map $\pi_1(S,\ast)\to \pi_1(S',\ast)$. By concatenating with an arc $A$ in $S'\setminus{S}$ connecting $\ast$ to $\ast'$, we obtain a homomorphism $i_*\colon \pi\to \pi'$.

\begin{lemma}\label{lem:LTcont}
Any morphism $\iota\colon \S\to \S'$ induces a map $\iota_*\colon \L^T(\S)\to \L^T(\S')$ of graded Lie algebras.
\end{lemma}
Note that if we had chosen a different arc $A$ from $\ast$ to $\ast'$, the resulting map $\pi\to \pi'$ would differ from $i_*$ by conjugation in $\pi'$. Since $\Gamma_j^T(\S)$ is a central filtration, this shows that the  map $\iota_*\colon \L^T_j(\S)\to \L^T_j(\S')$ does not depend on the arc $A$.
\begin{proof}
As we noted in Section~\ref{sec:Torellicategory}, any morphism $\iota\colon \S\to \S'$ induces a diagram:
\begin{equation*}
\xymatrix{
 \pi_1(S,\ast)\ar[r]\ar_{i_*}[d] &H_1(S)\ar[r]\ar_{i_*}[d]& H_1(\Surhat)\ar[r]\ar^{\ihat_*}[d]&\HS\ar^{\iota_*}[d]\\
 \pi_1(S',\ast')\ar[r] &H_1(S')\ar[r] &H_1(\Surhat')\ar[r]& H(\S')
 }
\end{equation*}
By induction, it follows that $i_*(\Gamma_j^T(\S))\subset \Gamma_j^T(\S')$ for all $j\geq 1$, so $i_*$ induces a map $\L^T(\S)\to \L^T(\S')$ of graded Lie algebras.
\end{proof}

\para{Injectivity of $\L^T(\S)\to \L^T(\S')$} The map $\iota_*\colon \L^T(\S)\to \L^T(\S')$ is not
always injective; for example, if $\iota$ is a simple capping,
$i_*(\zeta_1)$ is nullhomotopic, so we will have $\iota_*(z_1)=0$. However, this is essentially the
only way that injectivity can fail.
\begin{theorem}\label{thm:LTfree}
  $\L^T(\S)$ is the free Lie algebra on the generating set
  $\mathcal{S}(\S)$ defined in Proposition~\ref{prop:cS}. Furthermore
  any morphism $\S\to \S'$ such that no component of $S'\setminus
  S$ is a disk induces an injection $\L^T(\S)\hookrightarrow
  \L^T(\S')$.
\end{theorem}
\begin{proof}
We will show that $\L^T(\S)$ is free on the claimed basis in the course of proving that $\L^T(\S)\to \L^T(\S')$ is injective. So consider a morphism $\S\to \S'$ such that no component of $S'\setminus S$ is a disk.

We begin by reducing to the case when $S'$ has only one boundary component. Given such a morphism $\S\to \S'$, let $S''$ be obtained from $\Surhat'$ by attaching a surface $S_{1,1}$ to each component of $\Surhat'$ except the one containing the basepoint. This certainly has only one boundary component, so it remains to check that no component of $S''\setminus S$ is a disk. Each component of $S''\setminus\Surhat'$ has genus 1, so any such disk must be contained in $\Surhat'$. Each component of $\Surhat'\setminus S'$ has at least two boundary components by definition, so any disk must be contained in $S'$. This shows that as long as no component of $S'\setminus S$ was a disk,  no component of $S''\setminus S$ is a disk. And of course, if we can prove that the composition $\L^T(\S)\to \L^T(\S')\to \L^T(\S'')$ is injective, then the first map $\L^T(\S)\to \L^T(\S')$ is necessarily injective as well.\\

Assume that $S'$ is a surface with one boundary component, which we may consider as a (trivial) partitioned surface $\S'=(S',\{P_0\},\ast')$. As we noted in Remark~\ref{rem:TSoneboundary}, $\L^T(\S')$ is the free Lie algebra $\L(S')$.
Recall that a subset $Y$ of a Lie algebra is called \emph{independent} if the
  subalgebra generated by $Y$ is free with basis $Y$. Given a morphism $\iota\colon \S\to \S'$ so that no component of $S'\setminus S$ is a disk, we will prove that $\iota_*\colon \L^T(\S)\to \L(S')$ takes the generating set $\cS(\S)$ to an independent subset of $\L(S')$. By the universal property, this implies that $\iota_*$ is an isomorphism of $\L^T(\S)$ onto its image $\iota_*(\L^T(\S))$. This will simultaneously show that $\iota_*$ is injective, and that $\cS(\S)$ is a free basis for $\L^T(\S)$.\\

As in the definition of $\cS(\S)$, for $0\leq i\leq k$ we set $m_i=\abs{P_i}-1$. Let $\{\xi_1,\ldots,\xi_{2g}\}\cup\{\alpha_i^j\}$ (excluding $\alpha^0_0$) be the basis for $\pi$ constructed there. As before, choose disjoint simple closed curves $\gamma_i\subset S$ cobounding a genus 0 surface with the boundary components lying in $P_i$. Let $\delta_i$ be a simple closed curve in the complement $S'\setminus S$ that similarly cobounds a genus 0 surface with the components lying in $P_i$. Together, $\gamma_i$ and $\delta_i$ cobound a surface of genus $m_i$; let $g_i$ be the genus of the subsurface on the other side of $\delta_i$.
Extend the generators  $\{\xi_1,\ldots,\xi_{2g}\}\cup \{\alpha_i^1,\ldots,\alpha_i^{m_i}\}_{0\leq i\leq k}$ to a basis for $\pi_1(S')$ of the form
  \[\{\xi_1,\ldots,\xi_{2g}\}\cup \{\alpha_i^1,\ldots,\alpha_i^{m_i}\}_{0\leq i\leq k}\cup \{\beta_i^1,\ldots,\beta_i^{m_i}\}_{0\leq i\leq k}\cup \{\psi_i^1,\ldots,\psi_i^{2g_i}\}_{0\leq i\leq k}.\] By choosing this basis appropriately, we can ensure that
\begin{equation}\label{eq:etadefinition}\eta_i\coloneq [\psi_i^1,\psi_i^2]\cdots[\psi_i^{2g_i-1},\psi_i^{2g_i}]
\end{equation}
 represents $\delta_i$, and that we have the relation
   \begin{equation}\label{eq:zetaalphabeta}\eta_i=\alpha_i^0\beta_i^1\alpha_i^1\cdots
  \beta_i^{m_i}\alpha_i^{m_i}\overline{\beta_i^{m_i}}\cdots\overline{\beta_i^{1}}\end{equation}
  where $\overline{\beta}$ denotes the inverse $\beta^{-1}$.
  
Let $x_i$, $a_i^j$, $b_i^j$, and $y_i^j$ denote the image in $\L_1(S')$ of $\xi_i$, $\alpha_i^j$, $\beta_i^j$, and $\psi_i^j$ respectively; $\L_1(S')$ is the free Lie algebra on the generating set $\cS(S')=\{x_i\}\cup \{a_i^j\} \cup \{b_i^j\}\cup \{y_i^j\}$.
For any $x_i$, and for any $a_i^j$
  with $i\geq
  1$,  we have $\iota_*(x_i)= x_i$ and
  $\iota_*(a_i^j)=a_i^j$, but for $z_i$ the formula is not so simple.
However, comparing the expression \eqref{eq:zetaprime} for $\zeta_i$ with the expression  \eqref{eq:zetaalphabeta} for $\eta_i$, we see that $\iota_*(z_i)$ can be expressed as
  \[\iota_*(z_i) -[\eta_i]=[a_i^1,b_i^1]+[a_i^2,b_i^1+b_i^2]+\cdots+
  [a_i^{m_i},b_i^1+\cdots+b_i^{m_i}].\]
Combining this with \eqref{eq:etadefinition}, we obtain
\begin{equation}
\label{eq:iotazi}
\iota_*(z_i)=[a_i^1,b_i^1]+[a_i^2,b_i^1+b_i^2]+\cdots+
  [a_i^{m_i},b_i^1+\cdots+b_i^{m_i}]\ \ +\ \ [y_i^1,y_i^2]+\cdots+[y_i^{2g_i-1},y_i^{2g_i}].
  \end{equation}

  Let $Q$ denote the subset $Q\coloneq \iota_*(\cS(\S))=\{x_i\}\cup\{a_i^j\}\cup\{\iota_*(z_i)\}$ of $\L(S')$. We seek to show that $Q$ is
  independent, meaning that $Q$ is a free basis for the Lie subalgebra it generates (namely $\iota(\L^T(\S))$ itself).
    Note that by Shirshov~\cite{Shirshov} and Witt~\cite{WittU}, any subalgebra of the free Lie algebra $\L(S')$ is itself
  free on \emph{some} basis (at least after tensoring with $\Q$).  
  
  Let $\L_\Q=\L(S')\otimes \Q$, and identify each $q\in Q$ with its image
    in $\L_\Q$; since $\L(S')$ is torsion free (Witt~\cite[Theorem
    4]{WittD}), the map $\L(S')\to \L_\Q$ is an injection. The following
    theorem is proved by Shirshov in the course of proving
    \cite[Theorem 2]{Shirshov}: if for each $q\in Q$ we have that
    the leading term of $q$ is not in the subalgebra of $\L_\Q$
    generated by the leading terms of $Q\setminus\{q\}$, then $Q$ is
    independent as a subset of $\L_\Q$. (The leading term of
      $q$ is the highest degree homogeneous component of $q$. For an
      exposition in English of a closely related theorem, see
      Bryant--Kov{\'a}cs--St{\"o}hr~\cite{BKS}.) Since all our elements $q\in Q$ are homogeneous, we must show that $q$ is not in the subalgebra of $\L_\Q$ generated by $Q\setminus \{q\}$.
      
     For $q=x_i$ and $q=a_i^j$, this is easy. Given any subset $X\subset \cS(S')$ of the generating set $\cS(S')$,
      the elimination
    theorem for free Lie algebras implies that as a vector space,
    $\L_\Q$ splits as the direct sum of the algebra
    generated by $X$ with the ideal generated by $\cS(S')\setminus X$
    (see e.g.\ Bourbaki \cite[Chapter II, Section 2.9, Proposition
    10]{Bourbaki}). Since no other element of $Q$ involves the degree 1 generators $x_i$ or $a_i^j$, this implies that in this case $q$ is not even contained in the \emph{ideal} generated by $Q\setminus \{q\}$.
    
 For $q=\iota_*(z_i)$, we argue as follows. It cannot be that both $m_i=0$ and $g_i=0$: if $m_i=0$ we have $\abs{P_i}=1$, in which case the genus $g_i$ must  be at least 1 (otherwise the corresponding component of $S'\setminus S$ would be a disk). Thus at least one of the generators $b_i^1$ or $y_i^1$ appears in the expansion of $\iota_*(z_i)$. Since these generators appear in no other elements of $Q$, the elimination theorem again implies that $q=\iota_*(z_i)$ is not contained in the subalgebra of $\L_\Q$ generated by $Q\setminus\{q\}$. Shirshov's result thus shows that $Q$ is independent in
    $\L_\Q$; since $\L(S')$ is torsion-free, it follows that $Q$ is
    independent in $\L(S')$. This completes the proof that $\iota_*\colon \L^T(\S)\to \L(S')$ is injective, and that $\L^T(\S)$ is a free  Lie algebra with basis $\cS(\S)$.
\end{proof}

\begin{remark}
\label{rem:biggerindep}
For future reference, we note that the independent set $\iota_*(\cS(\S))$ from Theorem~\ref{thm:LTfree} can be extended slightly, by the same proof as above. 
These observations will be used in the proof of Theorem~\ref{thm:preserved}.

First, if $g_1>0$, then $\iota_*(\cS(\S))\cup\{y_1^1\}$ is independent. Since $y_1^1$ has degree 1, it is certainly not contained in $\iota_*(\cS(\S))$. The only generator in $\iota_*(\cS(\S))$ involving $y_1^1$ is $\iota_*(z_1)$, so there is nothing additional to verify except for this generator. But since $\iota_*(z_1)$ remains the only generator involving $y_1^2$, we see that $\iota_*(z_1)$ is still not contained in the subalgebra generated by $\iota_*(\cS(\S))\cup\{y_1^1\}\setminus \{\iota_*(z_1)\}$. By Shirshov's theorem this implies that $\iota_*(\cS(\S))\cup\{y_1^1\}$ is independent.

Second, if $m_1\geq 2$, or if $m_1=1$ and $g_1\geq 1$, then $\iota_*(\cS(\S))\cup\{b_1^1\}$ is independent. Again, the only generator involving $b_1^1$ is $\iota_*(z_1)$. If $m_1\geq 2$, then $\iota_*(z_1)$ is the only generator involving $b_1^2$; if $g_1\geq 1$, then $\iota_*(z_1)$ is the only generator involving $y_1^1$. In either case, $\iota_*(z_1)$ is not contained in the subalgebra generated by $\iota_*(\cS(\S))\cup\{b_1^1\}\setminus \{\iota_*(z_1)\}$, and thus $\iota_*(\cS(\S))\cup\{b_1^1\}$ is independent.
\end{remark}

\para{The filtrations $\Gamma^T_j$ are preserved by inclusions} Given any morphism $\S\to \S'$, we can restrict the filtration $\Gamma^T_j(\S')$ from $\pi_1(S')$ to $\pi_1(S)$. The content of Theorem~\ref{thm:LTfree} is that as long as no component of $S'\setminus S$ is a disk, the induced
filtration is precisely $\Gamma_j^T(\S)$ itself.
\begin{corollary}\label{cor:Tpreimage}
For any
morphism $\iota\colon\S\to\S'$ so that no component of $S'\setminus S$ is a disk, $\Gamma_j^T(\S)=i_*^{-1}(\Gamma_j^T(\S'))$; in other words, $\Gamma_j^T(\S)=\pi_1(S)\cap\Gamma_j^T(\S')$.
\end{corollary}
\begin{proof}   Lemma~\ref{lem:LTcont} states that $\Gamma_j^T(\S)\subset i_*^{-1}(\Gamma_j^T(\S'))$ for all $j$. By Theorem~\ref{thm:LTfree}, $\iota_*$ is injective on each
  $\L^T_j$, which implies that $i_*^{-1}(\Gamma_j^T(\S'))\subset
  \Gamma_j^T(\S)$ for all $j$.
\end{proof}
 In particular, this corollary implies that we could have defined $\Gamma_j^T(\S)$ by embedding $\S$ into an arbitrary surface $S'$ with one boundary component, and restricting the lower central series $\Gamma_j(S')$ to $\pi_1(S)$. However, without Theorem~\ref{thm:LTfree} there would be no reason to think that this definition would be well-defined (independent of the choice of embedding $\S\to S'$).

\para{Totally separated partitioned surfaces}
When $\S$ is a totally separated partitioned surface, the generating set $\cS(\S)$ consists just of $x_1,\ldots,x_{2g}$ in degree 1 and $z_1,\ldots,z_k$ in degree 2. We will need the following proposition in Section~\ref{sec:restrict} when we bound the image of the Johnson homomorphism. Note that for a totally separated surface, $\L^T_1(\S)$ coincides with $\HS$; anticipating the notation of Section~\ref{sec:generaldef}, we write $\NS$ for $\L_2^T(\S)$. 
\begin{proposition}
\label{prop:Tbracket}
If $\S$ is totally separated, the commutator bracket $\L_1^T(\S)\otimes \L_2^T(\S)\to \L_3^T(\S)$ induces the short exact sequence
\[1\to \bwedge^3 \HS\to \HS\otimes
\NS\to \L^T_3(\S)\to 1.\]
\end{proposition}
\noindent
The kernel is simply the Jacobi identity between elements of $\L_1^T(\S)=\HS$. Formally, the embedding $\bwedge^3 \HS\to \HS\otimes \NS$ is defined by sending $x\wedge y\wedge z\in \bwedge^3 \L_1^T(\S)$ to $x\otimes [y,z]+y\otimes [z,x]+z\otimes [x,y]\in \L_1^T(\S)\otimes \L_2^T(\S)$; the Jacobi identity asserts precisely that elements of this form are annihilated by the Lie bracket. 

Proposition~\ref{prop:Tbracket} is not just a corollary of Theorem~\ref{thm:LTfree}, but actually a special case of the theorem: the proposition states that there are no nontrivial relations among the basis elements $\{x_1,\ldots,x_{2g},z_1,\ldots,z_k\}$ in degree 3, while Theorem~\ref{thm:LTfree} states that there are no nontrivial relations among them at all.

\para{$\IS$ acts trivially on $\L^T(\S)$} Theorem~\ref{thm:LTfree} has another important consequence, which we will use in Sections~\ref{sec:Johnsonfiltration} and \ref{sec:generaldef}.
\begin{corollary}\label{cor:TorLTtrivial} For any partitioned surface $\S$, 
the mapping class group $\Mod(S)$ preserves the filtration $\Gamma^T_j(\S)$ on $\pi_1(S,\ast)$. Moreover, the action of $\Mod(S)$ on $\L^T(\S)$ factors through its action on $\HS$; in particular, $\IS$ acts trivially on $\L^T(\S)$.
\end{corollary}
\begin{proof}
As in the proof of Theorem~\ref{thm:LTfree}, given any partitioned surface $\S$ we may construct an inclusion $\S\to \S'$ so that $\S'$ has a single boundary component and  every component of $S'-S$ has genus at least 1. For such an inclusion the map $\Mod(S)\to \Mod(S')$ is injective (see Paris--Rolfsen~\cite[Corollary 4.2(iii)]{ParisRolfsen}). The lower central series of $\pi_1(S')$ is preserved by $\Mod(S')$, and the subgroup $\pi_1(S)$ is preserved by the subgroup $\Mod(S)$. We conclude that the intersection $\pi_1(S)\cap \Gamma_j(S')=\Gamma^T_j(\S)$ is preserved by $\Mod(S)$.

For the second claim, first assume that $\S$ is totally separated, so that $\zeta_i$ represents a boundary component of $S$. Since $\Mod(S)$ fixes the boundary components of $S$, it fixes $\zeta_i$ up to conjugacy, and thus acts trivially on $z_1,\ldots,z_k\in \L^T_2(\S)$. Proposition~\ref{prop:cS} shows that $\L^T_1(\S)$ and $\{z_1,\ldots,z_k\}$ generate $\L^T(\S)$, and so we conclude that the action of $\Mod(S)$ on $\L^T(\S)$ factors through its action on $\L^T_1(\S)=\HS$. In particular, since $\IS$ acts trivially on $\L_1^T(\S)=\HS$ by definition, we see that $\IS$ acts trivially on all of $\L^T(\S)$.  If $\S$ is not totally separated, Theorem~\ref{thm:LTfree} gives us a $\Mod(S)$--equivariant embedding of $\L^T(\S)$ into $\L^T(\Shat)$. Thus the action of $\Mod(S)$ on $\L^T(\S)$ factors through its action on $\HS=H(\Shat)$, as desired.
\end{proof}

\section{The Johnson filtration}
\label{sec:Johnsonfiltration}
Let $\S$ be a partitioned surface, and let $\pi=\pi_1(S,\ast)$ as before. By Corollary~\ref{cor:TorLTtrivial}, the action of $\Mod(S)$ on $\pi$ preserves the central series $\Gamma_i^T(\S)$ defined in Section~\ref{sec:charT}. We will use the action of $\Mod(S)$ on this central series to define the \emph{partitioned Johnson filtration} \[\Mod(S)=\MS{1}>\MS{2}>\MS{3}>\cdots\]
The classical Johnson filtration for a surface with one boundary component consists of those homeomorphisms acting trivially modulo $\Gamma_k(S)$, but for partitioned surfaces we need to impose another condition.

\subsection{The partitioned Johnson filtration}
\label{sec:Johnsonfiltproper}
\para{Action on arcs}
If $A$ is an arc in $S$ from the basepoint $\ast\in\partial S$ to another point $p$ lying in $\partial S$, and $\phi$ is an element of $\Mod(S)$, we define the element $d_A(\phi)\in \pi$ as follows. We denote by $A^{-1}$ the same arc parametrized in reverse, from $p$ to $\ast$. For any $\phi\in\Mod(S)$,  the image $\phi(A)$ is another arc with
the same endpoints as $A$, so we can consider $\phi(A)A^{-1}$ as a loop based at $\ast$. We define $d_A(\phi)\coloneq \phi(A)A^{-1}\in \pi$ to be the resulting element of the fundamental group.

For the following definition, we enumerate the boundary components of $S$, and choose arcs $A_i$ beginning at $\ast$ and ending on the $i$th  component of $\partial S$.

\begin{definition}[The partitioned Johnson filtration]\label{def:Johnsonfiltration} 
For $k\geq 1$, let $\MS{k}$ be the subgroup of $\Mod(S)$ consisting of those $\phi\in \Mod(S)$ satisfying the three following conditions:
\begin{enumerate}[(i)]
\item $\phi$ acts trivially on $\pi_1(S,\ast)$ modulo $\Gamma_k^T(\S)$.
\item for each arc $A_i$, the element $d_{A_i}(\phi)=\phi(A_i)A_i^{-1}$ is contained in $\Gamma_{k-1}^T(\S)$.
\item if $A_i$ and $A_j$ end at two components lying in the same block $P\in \P$,\\ then $d_{A_i}(\phi)\equiv d_{A_j}(\phi)$ modulo $\Gamma^T_k(\S)$.
\end{enumerate}
\end{definition}
\begin{remark}
Conditions (ii) and (iii) do not appear in the definition of the classical Johnson filtration $\Mod_{(k)}(S_{g,1})$, which consists simply of elements $\phi\in \Mod(S_{g,1})$ acting trivially on $\pi_1(S_{g,1})/\Gamma_k(S_{g,1})$, and the reader might wonder if they can be removed. However, if we hope to restrict the Johnson filtration to subsurfaces $S$ with multiple boundary components, such a condition on arcs is unavoidable. We can see this just from considering Dehn twists, as follows.

Consider a subsurface $S\subset S_{g,1}$ with multiple boundary components, let $\gamma$ be a boundary component in $\partial S$ \emph{not} containing the basepoint $\ast$, and consider the Dehn twist $T_\gamma$. We know that no Dehn twist around any curve ever lies in $\Mod_{(4)}(S_{g,1})$. 
 But $T_\gamma$ acts trivially on $\pi_1(S,\ast)$, so we cannot exclude it based solely on its action on $\pi_1(S,\ast)$.
Condition (ii) is what guarantees that Dehn twists behave as we expect: for nonseparating curves we will have $T_\gamma\not\in \MS{2}$, and  for separating curves we will have $T_\gamma\in \MS{3}$ but $T_\gamma\notin \MS{4}$.
\end{remark}

Before moving on, let us establish that Definition~\ref{def:Johnsonfiltration} is well-defined, meaning that it does not depend on the choice of the arcs $A_i$. So let $A$ and $B$ be two arcs with the same endpoints in $\partial S$, and assume that $\phi$ acts trivially on $\pi$ modulo $\Gamma^T_k(\S)$ and that $d_A(\phi)\in \Gamma^T_{k-1}(\S)$. We can compute $d_B(\phi)d_A(\phi)^{-1}$ as
\[(\phi(B)B^{-1})(\phi(A)A^{-1})^{-1}=\phi(BA^{-1})(BA^{-1})^{-1}[BA^{-1},\phi(A)A^{-1}]=\phi(x)x^{-1}[x,d_A(\phi)]\] where $x=BA^{-1}\in \pi$. 
Since $\phi$ acts trivially modulo $\Gamma_k^T(\S)$, the first term $\phi(x)x^{-1}$ lies in $\Gamma_k^T(\S)$. By assumption $d_A(\phi)\in \Gamma_{k-1}^T(\S)$, so $[x,d_A(\phi)]\in \Gamma_k^T(\S)$. This shows that $d_A(\phi)$ and $d_B(\phi)$ agree modulo $\Gamma_{k}^T(\S)$; in particular, if one lies within $\Gamma_{k-1}^T(\S)$, the other does. Similarly, if $d_A(\phi)\equiv d_{A_j}(\phi)$ modulo $\Gamma_k^T(\S)$, then $d_B(\phi)\equiv d_{A_j}(\phi)$ modulo $\Gamma_k^T(\S)$ as well.

\para{Fundamental computation} One motivation for defining $d_A(\phi)$ is the following fundamental computation, which we will use repeatedly throughout the paper. Let $\gamma=A \delta A^{-1}$, where $\delta$ is some loop fixed by $\phi$. (For example, if $S$ is contained in a larger surface $S'$, we might choose $\delta$ contained in $S'\setminus S$.) Then we have
\begin{equation}
\label{eq:fund}
\phi(\gamma)\gamma^{-1}
  =\phi(A)\delta \phi(A)^{-1}A\delta^{-1} A^{-1}=\phi(A)A^{-1}(A\delta A^{-1})A\phi(A)^{-1}(A\delta A^{-1})=[d_A(\phi),\gamma].
  \end{equation}

\begin{remark}
\label{rem:classicalagrees}
If $A_0$ is an arc from $\ast$ to itself that is nullhomotopic, the same is true of $\phi(A_0)$, so $d_{A_0}(\phi)$ is trivial for any $\phi\in \Mod(S)$. Thus condition (iii) implies that for any arc $A$ from $\ast$ to a point in $P_0$ (the block containing the basepoint), $d_A(\phi)\in \Gamma_k^T(\S)$.
Furthermore, when $S$ has only a single boundary component, the only arc is $A_0$, so conditions (ii) and (iii) are vacuous. By Remark~\ref{rem:TSoneboundary}, in this case $\Gamma^T_j(\S)$ is just the lower central series $\Gamma_j(\pi)$. Thus for surfaces with one boundary component, the partitioned Johnson filtration $\MS{k}$ coincides with the classical Johnson filtration $\Mod_{(k)}(S)$.
\end{remark}

\para{Refining the partition}
If $\S=(S,\P,\ast)$ and $\S'=(S,\P',\ast)$ are two partitioned surfaces coming from two different partitions on the same surface, we can compare the resulting filtrations $\Mod_{(k)}(\S)$ and $\Mod_{(k)}(\S')$ on the mapping class group $\Mod(S)$. We have the following comparison result. Say that $\P'$ is \emph{finer than} $\P$ if every block $P'\in \P'$ is contained in a single block $P\in \P$ (but a single block of $\P$ may split into multiple blocks of $\P'$). This encodes the notion that the partition $\P'$ is ``more separated'' than $\P$; for example, the totally separated partition is finer than any other partition. 
\begin{proposition}
\label{prop:refinement}
Given $\S=(S,\P,\ast)$ and $\S'=(S,\P',\ast)$, if $\P'$ is finer than $\P$, then for any $k\geq 1$ we have $\Mod_{(k)}(\S)\subset \Mod_{(k)}(\S')$.
\end{proposition}
\begin{proof}
Given a boundary component $c\in \pi_0(\partial S)$, let $a_c$ be the associated element of $H_1(S)$. We noted in the proof of Lemma~\ref{lem:TSgen} that the kernel of the map $H_1(S)\to H(\S)$ is spanned by the elements $\sum_{c\in P}a_c$ for each block $P\in \P$. Since $\P'$ is finer than $\P$, each block $P\in \P$ is the disjoint union $P=P'_1\cup \cdots\cup P'_\ell$ of blocks  $P'_i\in \P'$. Thus we can regroup this sum as $\sum_{c\in P}a_c=\sum_{i=1}^\ell \sum_{c\in P'_i}a_c$. Each of the latter terms vanishes in $H(\S')$, so we conclude that $\ker(H_1(S)\to H(\S))$ is contained in $\ker(H_1(S)\to H(\S'))$.

It follows that $T(\S)=\ker(\pi_1(S)\to H(\S))$ is contained in $T(\S')=\ker(\pi_1(S)\to H(\S'))$, and thus by induction that $\Gamma_j^T(\S)\subset\Gamma_j^T(\S')$ for all $j\geq 1$. Now if $\phi(x)x^{-1}$ (or $d_A(\phi)$, or $d_A(\phi)d_{A'}(\phi)^{-1}$, respectively) lies in $\Gamma_k^T(\S)$, it also lies in $\Gamma_k^T(\S')$. It follows that $\phi\in \Mod_{(k)}(\S) \implies \phi\in \Mod_{(k)}(\S')$, as desired.\end{proof}

\para{The first terms of the Johnson filtration} By definition, $\MS{1}=\Mod(S)$, and $\MS{2}$ is the Torelli group $\IS$ defined in Section~\ref{sec:Torelli}. (This follows from Theorem~\ref{thm:preserved}, but it is not difficult to verify directly.)
We  denote the next term $\MS{3}$ by $\KS$ and call it the \emph{Johnson kernel of $\S$}. In Section~\ref{sec:Jh} we will define the \emph{partitioned Johnson
homomorphism} $\tau_\S$, and prove in Theorem~\ref{thm:deltas} that $\KS=\ker\tau_\S$. In particular, we will see that when $\S$ is totally separated, $\KS$ is exactly the subgroup of $\Mod(S)$ acting trivially on $\pi_1(S)$ modulo $\Gamma_3^T(\S)$; conditions (ii) and (iii) in Definition~\ref{def:Johnsonfiltration} are not necessary in this case.

\para{Changing the basepoint} The filtration $\MS{k}$ is defined in terms of the partitioned surface $\S=(S,\P,\ast)$, and it is easy to see that this filtration does depend on the partition $\P$ (we will see many examples in Section~\ref{sec:Jh}). However, we have the convenient property that the filtration $\MS{k}$ does \emph{not} depend on the basepoint $\ast$.
\begin{theorem}
\label{thm:basepointirrelevant}
  The Johnson filtration does not depend on the basepoint; that is, if
  $\S$ and $\S'$ differ only in the location of the basepoint, then
  $\Mod_{(k)}(\S)=\Mod_{(k)}(\S')$.
\end{theorem}
\begin{proof}
If $\S=(S,\P,\ast)$, let $\S'=(S,\P,\ast')$. An isomorphism from $\pi=\pi_1(S,\ast)$ to
$\pi'=\pi_1(S,\ast')$ is given by $x\mapsto A^{-1}x A$, where
$A$ is an arc from $\ast$ to $\ast'$. We saw in Lemma~\ref{lem:LTcont} that this isomorphism takes the central series $\Gamma_k^T(\S)$ of $\pi$ to the central series $\Gamma_k^{T}(\S')$ of $\pi'$. In this proof only, let $\overline{x}$ denote $A^{-1}xA$. We compute
\begin{equation}\label{eq:basechange} \begin{array}{rl}
  \phi(\overline{x})\overline{x}^{-1}&=\phi(A^{-1}xA)A^{-1}x^{-1}A\\[4pt]
  &=[\phi(A^{-1})A,A^{-1}\phi(x)A]A^{-1}\phi(x)x^{-1}A\\[4pt]
&=[\overline{d_A(\phi)^{-1}},\overline{\phi(x)}]\overline{\phi(x)x^{-1}}=\overline{[d_A(\phi)^{-1},\phi(x)]\phi(x)x^{-1}}
\end{array}\end{equation}
If $\phi\in \MS{k}$ we know that $\phi(x)x^{-1}\in \Gamma_k^T(\S)$, and that $d_A(\phi)\in \Gamma_{k-1}^T(\S)$, so $[d_A(\phi)^{-1},\phi(x)]\in \Gamma_k^T(\S)$ as well. This shows that $\phi$ acts trivially on $\pi'$ modulo $\Gamma_k^{T}(\S')$, verifying condition (i) of Definition~\ref{def:Johnsonfiltration}. For conditions (ii) and (iii), note that for the arc $A'_j$ from $\ast'$ to the $j$th boundary component we may take $A'_j=A^{-1}A_j$. We compute $d_{A'_j}(\phi)$ as: \[\phi(A'_j){A'_j}^{-1}=\phi(A)^{-1}\phi(A_j)A_j^{-1}A=\overline{A\phi(A)^{-1}\phi(A_j)A_j^{-1}}=\overline{d_A(\phi)^{-1}d_{A_j}(\phi)}\]
Since $d_A(\phi)$ and $d_{A_j}(\phi)$ lie in $\Gamma_{k-1}^T(\S)$, we conclude that $\phi(A'_j){A'_j}^{-1}$ lies in $\Gamma_{k-1}^{T}(\S')$, verifying condition (ii). Finally, if $A_i$ and $A_j$ end in the same block $P\in P$, we have $d_{A_i}(\phi)\cong d_{A_j}(\phi)$ modulo $\Gamma^T_k(\S)$. From the computation above, this implies that $d_{A'_i}(\phi)\cong  d_{A'_j}(\phi)$ modulo $\Gamma^T_k(\S')$. 
\end{proof}

\subsection{The Johnson filtration is preserved by inclusions}
\label{section:Johnsonfiltpreserved}
In this section we prove the fundamental result that the Johnson filtration is sharply preserved by inclusions, meaning that (with a small list of exceptions) any inclusion $\iota\colon \S\to \S'$ satisfies \begin{equation}
\label{eq:sharplypreserved}
\phi\in \MS{k} \iff \iota_*(\phi)\in \Mod_{(k)}(\S').
\end{equation} This implies that the restriction of the Johnson filtration to a subsurface $S$ depends \emph{only} on which boundary components of $S$ become homologous in the larger surface.

Of course, the property \eqref{eq:sharplypreserved} can never hold if the map $\Mod(S)\to \Mod(S')$ is not injective. Fortunately this happens only rarely: by
Paris--Rolfsen~\cite[Corollary 4.2(iii)]{ParisRolfsen}, the map
  $\Mod(S)\to \Mod(S')$ induced by an inclusion $S\to S'$ is an injection if and only if no component of $S'\setminus S$ is a
  disk or annulus with boundary contained in $\partial S$.

\begin{theorem}\label{thm:preserved}
For any inclusion $\iota\colon \S\to \S'$ we have $\iota_*(\Mod_{(k)}(\S))\subset \Mod_{(k)}(\S')$ for all $k\geq 1$. Moreover as long as no component of $S'\setminus S$ is a   disk or annulus with boundary contained in $\partial S$, we have $\Mod_{(k)}(\S)=\iota_*^{-1}(\Mod_{(k)}(\S'))$; in other words,  $\Mod_{(k)}(\S)=\Mod(S)\cap \Mod_{(k)}(\S')$. \end{theorem}

\begin{proof}
Given $\phi\in \Mod(S)$, let $\overline{\phi}$ denote its image in $\Mod(S')$.
We first assume that $\phi\in \Mod_{(k)}(\S)$ and seek to show that $\phibar\in \Mod_{(k)}(\S')$. We begin with conditions (ii) and (iii). Let $C$ be a fixed arc in $S'\setminus S$ from $\ast'$ to $\ast$. For this proof only, given $\xi\in \pi_1(S,\ast)$ or $\xi\in \pi_1(S',\ast)$, let $\overline{\xi}$ denote $C\xi C^{-1}\in \pi_1(S',\ast')$. We verified in Lemma~\ref{lem:LTcont} that the map $\pi_1(S,\ast)\to \pi_1(S',\ast')$ defined by $\xi\mapsto \overline{\xi}$ takes the filtration $\Gamma_j^T(\S)$ faithfully to the filtration $\Gamma_j^T(\S')$.

We can always choose an arc $B_i$ from $\ast'$ to the $i$th component of $\partial S'$ with initial segment $C$, so that $B=CAD$ where $A_j$ is an arc in $S$ and $D$ is contained in $S'\setminus S$. We then have $d_{B}(\phibar)=\phibar(B)B^{-1}=(C\phi(A)D)(CAD)^{-1}=\overline{d_A(\phi)}$.  Thus $d_B(\phibar)\in \Gamma_k^T(\S')$ if and only if $d_A(\phi)\in \Gamma_k^T(\S)$, verifying the claim for condition (ii). Now let $B=CAD$ and $B'=CA'D'$ be two such arcs ending at components lying in the same block $P'\in \P'$; we have $d_B(\phibar)d_{B'}(\phibar)^{-1}=\overline{d_A(\phi)d_{A'}(\phi)^{-1}}$. Since the inclusion $\S\to \S'$ respects the partitions, $A$ and $A'$ necessarily end at components lying in the same block $P\in \P$. Condition (iii) for $\phi$ guarantees that $d_A(\phi)d_{A'}(\phi)^{-1}\in \Gamma^T_k(\S)$, so the above formula shows that $d_B(\phibar)d_{B'}(\phibar)^{-1}\in \Gamma^T_k(\S')$, as desired.

We now handle condition (i). Putting a generic element of $\pi'=\pi_1(S',\ast')$ in general position with respect to $S$ shows that $\pi'$ is generated by loops of the following four forms:
\begin{itemize}
\item first, loops $\gamma$ contained entirely in $S'\setminus S$;
\item second, loops $\gamma=C\xi C^{-1}$ for $\xi\in \pi_1(S,\ast)$;
\item third, loops $\gamma=CA\delta A^{-1} C^{-1}$ for $A$ an arc contained in $S$ and $\delta$ a loop contained in $S'\setminus S$; 
\item fourth, loops $\gamma=CA B A'^{-1}C^{-1}$ for $A, A'$ arcs contained in $S$ and $B$ an arc contained in $S'\setminus S$.
\end{itemize}

Our goal is to show that $\phibar(\gamma)\gamma^{-1}\in \Gamma_k^{T}(\S')$ for any such $\gamma$. In the first case this is trivial, since $\phibar(\gamma)=\gamma$. In the second case we have $\phibar(\gamma)\gamma^{-1}=\overline{\phi(\xi)\xi^{-1}}$; since $\phi(\xi)\xi^{-1}\in\Gamma_k^T(\S)$ by assumption, this implies $\phibar(\gamma)\gamma^{-1}\in \Gamma_k^{T}(\S')$. In the third case, just as in \eqref{eq:fund} we compute:
\[\phibar(\gamma)\gamma^{-1}
  =C\phi(A)\delta \phi(A)^{-1}A\delta A^{-1}C^{-1}=\overline{[\phi(A)A^{-1}, A\delta A^{-1}]}= [\overline{d_A(\phi)},\gamma].\]
By condition (ii) we know $\overline{d_A(\phi)}\in \Gamma^T_{k-1}(\S')$, so $\phibar(\gamma)\gamma^{-1}=[\overline{d_A(\phi)},\gamma]\in\Gamma^T_k(\S')$, as desired. Finally, in the fourth case, we compute:
\begin{align*}\phibar(\gamma)\gamma^{-1}
  &=C\phi(A)B \phi(A')^{-1}A'B^{-1} A^{-1}C^{-1}\\
  &=\overline{d_A(\phi)}\gamma \overline{d_{A'}(\phi)}^{-1}\gamma^{-1}\\&= [\overline{d_A(\phi)},\gamma]\cdot \gamma \overline{d_A(\phi)d_{A'}(\phi)^{-1}}\gamma^{-1}
\end{align*}
As before, by condition (ii) we know $\overline{d_A(\phi)}\in \Gamma_{k-1}^T(\S')$, so the first term lies in $\Gamma^T_k(\S')$. Since $A$ and $A'$ are connected by the arc $B$ lying outside $S$, they must end in the same block. Thus by condition (iii) we have $\overline{d_A(\phi)d_{A'}(\phi)^{-1}}\in \Gamma_k^T(\S')$, so the second term lies in $\Gamma^T_k(\S')$ as well, showing that $\phibar(\gamma)\gamma^{-1}\in \Gamma^T_k(\S')$. This concludes the proof that $\phibar\in \Mod_{(k)}(\S')$.\\

Now assume that no component of $S'\setminus S$ is a disk or annulus. We assume that $\phibar\in \Mod_{(k)}(\S')$ and seek to show that $\phi\in \Mod_{(k)}(\S)$. Condition (i) is easy to check: for any $\xi\in \pi_1(S,\ast)$, we seek to show that $\phi(\xi)\xi^{-1}$ lies in $\Gamma_k^T(\S)$. By Lemma~\ref{lem:LTcont}, this is equivalent to showing that $\overline{\phi(\xi)\xi^{-1}}\in \Gamma_k^T(\S)$. But as we saw above, $\phibar(\overline{\xi})\overline{\xi}^{-1}=\overline{\phi(\xi)\xi^{-1}}$, and the assumption that $\phi\in \Mod_{(k)}(\S')$ says precisely that $\phibar(\overline{\xi})\overline{\xi}^{-1}$ lies in $\Gamma_k^T(\S')$.

Conditions (ii) and (iii) require more work. The difficulty in running the above arguments in reverse is that although \eqref{eq:fund} will show that $[d_A(\phi),\gamma]\in \Gamma_k^T(\S)$, this does \emph{not} imply that $d_A(\phi)\in \Gamma_{k-1}^T(\S)$. For example we could have $d_A(\phi)=\gamma$, in which case $[d_A(\phi),\gamma]=[\gamma,\gamma]=1$. We will need to show that we can choose loops $\gamma\in \pi_1(S')$ that are sufficiently ``independent'' from $\Gamma_k^T(\S)$.

Consider an arc $A$ from $\ast$ to a boundary component of $S$, and let $U$ be the component of $S'\setminus S$ containing the endpoint of $A$. We first handle the straightforward case when $\partial U$ is not contained in $\partial S$. Let $A'$ be another arc in $S$ beginning at $\ast$ and ending in the same block $P\in\P$ as $A$, and let $U'$ be the component of $S'\setminus S$ containing the endpoint of $A'$.

We may choose an arc $B=CAD$ from $\ast'$ to the boundary $\partial S'$ as before. We previously calculated that $d_B(\phibar)=\overline{d_A(\phi)}$, so since $d_B(\phibar)\in \Gamma_k^T(\S')$ by assumption, we conclude that $d_A(\phi)\in \Gamma^T_k(\S)$, verifying condition (ii). Since the inclusion $\S\to \S'$ respects the partitions, we must also have $\partial T'\not\subset \partial S$. Thus we may choose an arc $B'=CA'D'$ from $\ast'$ to $\partial S'$ so that $B$ and $B'$ end in the same block $P'\in \P'$. By assumption $d_B(\phibar)d_{B'}(\phibar)^{-1}=\overline{d_A(\phi)d_{A'}(\phi)^{-1}}$ lies in $\Gamma_k^T(\S')$. Thus $d_A(\phi)d_{A'}(\phi)^{-1}\in \Gamma_k^T(\S)$, verifying condition (iii).

Next we consider the case when $\partial U$ has a single component, which necessarily meets in $\partial S$ in a singleton block of $\P$. Thus Condition (iii) is vacuous for this block.
Since $U$ is not a disk, the genus of $U$ must be positive, so we may choose $\gamma=CA\delta A^{-1} C^{-1}$ where $\delta$ is a free generator for $\pi_1(U)$ descending to a generator of $H_1(U,\partial U)$. In the notation of Remark~\ref{rem:biggerindep}, we may take $\gamma$ to represent the generator $y_1^1$ of $\L^T(\S')$. As in \eqref{eq:fund} we computed above that $\phibar(\gamma)\gamma^{-1}=[\overline{d_A(\phi)},\gamma]$, so by assumption we know that $[\overline{d_A(\phi)},\gamma]\in \Gamma_k^T(\S')$. Using the fact that $d_A(\phi)\in\pi_1(S)$, our goal is to show that this implies $\overline{d_A(\phi)}\in \Gamma_{k-1}^T(\S')$.

Remark~\ref{rem:biggerindep} shows that the subalgebra $\mathcal{A}$ of $\L^T(\S')$ generated by $\iota_*(\L^T(\S))$ and $y_1^1$ is free with basis $\iota_*(\cS(\S))\cup\{y_1^1\}$. By the elimination theorem for free Lie algebras \cite[Ch. II, Sec. 2.9, Corollary]{Bourbaki}, the ideal of $\mathcal{A}$ generated by $\iota_*(\cS(\S))$ is a free Lie algebra on the basis of iterated brackets $\{[y_1^1,[y_1^1,\cdots,[y_1^1,x]]]\big| x\in \iota_*(\cS(\S))\big\}$. Furthermore $\mathcal{A}$ is the direct sum of this ideal with $\langle y_1^1\rangle$. In particular, this implies that $[y_1^1,s]\neq 0$ for any $s\in\mathcal{A}$ satisfying $s\not\in\langle y_1^1\rangle$.

Let $\ell$ be the largest integer such that $\overline{d_A(\phi)}\in \Gamma_\ell^T(\S')$, and let $r$ denote the image of $\overline{d_A(\phi)}$ in $\iota_*(\L^T_\ell(\S))\subset \L^T_\ell(\S')$; by definition, $r\neq 0\in \L^T_\ell(\S')$. Thus the above remark implies that $[y_1^1,r]\neq 0\in \L^T_{\ell+1}(\S')$. However we know that $[\overline{d_A(\phi)},\gamma]\in \Gamma^T_k(\S')$ represents $[y_1^1,r]$. We conclude that $\ell\geq k-1$, and thus that $d_A(\phi)\in \Gamma_{k-1}^T(\S)$, as desired.

Finally, we consider the remaining case: when $\partial U$ has multiple components, all contained in $\partial S$. Since $U$ is neither a disk nor an annulus, it either has positive genus or at least three boundary components. Choose an arc $A'$ from $\ast$ to one of the other components of $\partial U$, and let $\gamma=CABA'^{-1}C^{-1}$ for $B$ an arc in $U$. In the notation of Remark~\ref{rem:biggerindep}, we may assume that $\gamma$ represents the generator $b_1^1$ of $\L^T(\S')$. We computed above that $\phibar(\gamma)\gamma^{-1}=[\overline{d_A(\phi)},\gamma]\cdot \gamma \overline{d_A(\phi)d_{A'}(\phi)^{-1}}\gamma^{-1}$. Choose $\ell$ as before and let $r$ denote the image of $\overline{d_A(\phi)}$ in $\iota_*(\L^T_\ell(\S))$; similarly, let $j$ be the largest integer such that $d_A(\phi)\equiv d_{A'}(\phi)$ modulo $\Gamma^T_j(\S)$, and let $s$ denote the image of $\overline{d_A(\phi)d_{A'}(\phi)^{-1}}$ in $\iota_*(\L^T_j(\S))$. 

Let $\mathcal{A}$ be the subalgebra of $\L^T(\S')$ generated by $\iota_*(\L^T(\S))$ and $b_1^1$, and consider the element $[r,b_1^1]+s\in \mathcal{A}$. Our assumption on $T$ implies that either $g_1\geq 1$ or $m_1\geq 2$, so Remark~\ref{rem:biggerindep} shows that $\mathcal{A}$ is free with basis $\iota_*(\cS(\S))\cup\{b_1^1\}$. Let $I$ denote the ideal of $\mathcal{A}$ generated by $b_1^1$. The elimination theorem states that $\mathcal{A}$ is the direct sum of $\iota_*(\cS(\S))$ and $I$, and that a basis for $I$ as a free Lie algebra is in bijection with a basis for the universal enveloping algebra of $\L^T(\S)$
\cite[Ch. II, Sec. 2.9, Corollary]{Bourbaki}. In particular, the adjoint action of $\iota_*(\L^T(\S))$ on $I$ factors through its action on its universal enveloping algebra, and thus is faithful. Since $r\neq 0\in \iota_*(\L^T_\ell(\S))$, this implies that $[r,b_1^1]\neq 0\in I_{\ell+1}\subset \L^T_{\ell+1}(\S')$. However, we know that $[\overline{d_A(\phi)},\gamma]\cdot \gamma \overline{d_A(\phi)d_{A'}(\phi)^{-1}}\gamma^{-1}\in \Gamma^T_k(\S')$, and so we conclude that $[r,b_1^1]+s\in \L^T_k(\S')$. By the splitting of $\mathcal{A}$ as a direct sum, this implies that $s\in \L^T_k(\S')$ and $[r,b_1^1]\in I_k\subset \L^T_k(\S')$. This shows that $\ell\geq k-1$ and $j\geq k$, which is to say that $d_A(\phi)\in \Gamma_{k-1}^T(\S)$ and that $d_A(\phi)d_{A'}(\phi)^{-1}\in \Gamma_k^T(\S)$. This verifies both conditions (ii) and (iii) in this final case.
\end{proof}

\subsection{An application of the partitioned Johnson filtration}
\label{section:BBMapplication}
In this section we give a specific application of the partitioned Johnson filtration to stabilizers of multicurves.
Theorem~\ref{thm:BBM} below was used by Bestvina--Bux--Margalit in their
calculation of the cohomological dimension of $\K(S_g)$, where it was stated as \cite[Lemma~6.16]{BBM}.

Let $R$ be a surface with 
one boundary component, and let $i\colon S\hookrightarrow R$ be the inclusion of a subsurface so that no component of $R\setminus S$ is a disk. Fix a single boundary component $D$ of $S$, and let
$\overline{S}$ be the surface obtained from $S$ by capping off all
boundary components except for $D$. Let $\pi\colon S\to \overline{S}$ be the natural inclusion of surfaces, and $\pi_*\colon \Mod(S)\twoheadrightarrow \Mod(\overline{S})$ the  surjection obtained by extension by the identity. Given an element $\phi\in \Mod(R)$ which
stabilizes $S$, we may restrict it to an element $\widetilde{\phi}\in \Mod(S)$, and consider its image $\pi_*(\widetilde{\phi})\in\Mod(\overline{S})$.

If any component of $R\setminus S$ is an annulus not containing $\partial R$, the map $i_*\colon \Mod(S)\to \Mod(R)$ will not be injective, so $\widetilde{\phi}$ is not uniquely defined. However \cite[Theorem 4.1(iii)]{ParisRolfsen} states that the kernel of $i_*\colon \Mod(S)\to \Mod(R)$ is generated by certain products of Dehn twists (so-called bounding pairs) around the corresponding boundary components of $S$. All but one of these boundary components becomes nullhomotopic in $\overline{S}$, so the corresponding Dehn twist vanishes in $\Mod(\overline{S})$. We conclude that the element $\pi_*(\widetilde{\phi})$ is
well-defined up to a twist $T_D$ around the boundary component $D$.
\begin{theorem}\label{thm:BBM}For any $\phi\in\K(R)$ stabilizing
  $S$, we have $\overline{\phi}\in\K(\overline{S})$.
\end{theorem}
Note that the twist $T_D$ around the boundary component $D$
  lies in $\K(\overline{S})$ (see e.g.\ Proposition~\ref{prop:st}), so  the desired conclusion is
  well-defined. (Since $T_D\not\in \Mod_{(4)}(\overline{S})$, it would not be well-defined to ask whether $\phi\in \Mod_{(k)}(R)$ implies $\overline{\phi}\in \Mod_{(k)}(\overline{S})$ when $k>3$.)
\begin{proof} Let $\ast'$ be a point in $\partial R$, and let $\S=(S,\P,\ast)$ be the induced partitioned surface  on $S$. Let $\S_1=(S,\P_1,\ast)$ denote the same surface but with the totally separated partition $\P_1$. Finally, choose a point $\ast_D\in D$, and let $\S_2=(S,\P_1,\ast_D)$ denote the same partitioned surface but with basepoint $\ast_D$. Note that the inclusion $\pi\colon S\to \overline{S}$ defines a morphism $\pi\colon \S_2\to \overline{S}$ of partitioned surfaces.

If no component of $R\setminus S$ is an annulus, the theorem follows from Theorem~\ref{thm:preserved} as follows. Since $i_*(\widetilde{\phi})=\phi$ lies in $\K(R)=\Mod_{(3)}(R)$, the second part of Theorem~\ref{thm:preserved} states that $\widetilde{\phi}\in \Mod_{(3)}(\S)$. Since the partition $\P_1$ is finer than $\P$, Proposition~\ref{prop:refinement} states that $\Mod_{(k)}(\S)\subset \Mod_{(k)}(\S_1)$ for all $k$; in particular, $\widetilde{\phi}\in \Mod_{(3)}(\S_1)$. By Theorem~\ref{thm:basepointirrelevant}, $\Mod_{(k)}(\S_1)=\Mod_{(k)}(\S_2)$ for all $k$, so $\widetilde{\phi}\in \Mod_{(3)}(\S_2)$. Finally, $\pi\colon \S_2\to \overline{S}$ is a morphism of partitioned surfaces, so the first part of Theorem~\ref{thm:preserved} implies that $\pi_*(\widetilde{\phi})\in \Mod_{(3)}(\overline{S})$. By Remark~\ref{rem:classicalagrees}, on a surface with one boundary component such as $\overline{S}$, the partitioned Johnson filtration agrees with the classical Johnson filtration, and so $\Mod_{(3)}(\overline{S})$ is just the classical Johnson kernel $\K(\overline{S})$.

The proof in general follows the same structure, but we cannot apply the second part of Theorem~\ref{thm:preserved} if some component of $R\setminus S$ is an annulus. Rather than dealing with $\K(\S)$ itself, we will check by hand that any lift $\widetilde{\phi}$ lies in the larger group $\K(\S_1)=\Mod_{(3)}(\S_1)$.

Let $C$ be a fixed arc in $R\setminus S$ from $\ast'$ to $\ast$, and let $\overline{x}$ denote $CxC^{-1}$. The proof that $\widetilde{\phi}$ acts trivially on $\pi_1(S,\ast)$ modulo $\Gamma_3^T(\S)$ still works without change, as follows. For any $x\in \pi_1(S,\ast)$, we know that $\overline{\widetilde{\phi}(x)x^{-1}}=\phi(\overline{x})\overline{x}^{-1}$ lies in $\Gamma_3(R)$, since $\phi\in \K(R)$. But by Corollary~\ref{cor:Tpreimage}, this implies $\widetilde{\phi}(x)x^{-1}\in \Gamma_3^T(\S)$. We showed in Proposition~\ref{prop:refinement} that $\Gamma_3^T(\S)\subset \Gamma_3^T(\S_1)$, so this verifies condition (i) in showing that $\widetilde{\phi}\in \Mod_{(3)}(\S_1)$. Since $\S_1$ is totally separated, condition (iii) is vacuous, so we need only check that $\widetilde{\phi}$ satisfies condition (ii). Moreover, for any arc not ending at an annular component of $R\setminus S$, the proof of Theorem~\ref{thm:preserved} shows that $d_A(\phi)\in \Gamma_2^T(\S)$, and thus that $d_A(\phi)\in \Gamma_2^T(\S_1)$.

Choose an annular component of $R\setminus S$, and let $D$ be one of its two components. Let $A$ be an arc in $S$ from $\ast$ to $D$, and let $\gamma=CA\delta A^{-1}C^{-1}$ where $\delta$ is a parametrization of the boundary component $D$. In the notation of Theorem~\ref{thm:LTfree}, we may assume that $A\delta A^{-1}$ represents $a_1^1\in \L^T_1(\S)$. As in \eqref{eq:fund} we compute $\phi(\gamma)\gamma^{-1}=C\widetilde{\phi}(A)\delta \widetilde{\phi}(A)^{-1}A\delta A^{-1}C^{-1}=\overline{[d_A(\phi),A\delta A^{-1}]}$. By assumption $\phi(\gamma)\gamma^{-1}\in \Gamma_3(R)$, so $[d_A(\phi),A\delta A^{-1}]\in \Gamma_3^T(\S)$. We will show that this implies that $d_A(\phi)$ is congruent to $A\delta^m A^{-1}$ modulo $\Gamma_2^T(\S)$ for some $m\in \Z$. Since $A\delta A^{-1}$ lies in $\Gamma_2^T(\S_1)$, this will show that $d_A(\phi)\in \Gamma_2^T(\S_1)$, as desired.

By Theorem~\ref{thm:LTfree}, the element $a_1^1$ is part of a free basis for $\L^T(\S)$. Thus just as in the proof of Theorem~\ref{thm:preserved}, applying the elimination theorem to the ideal generated by the other basis elements shows that the kernel of $[-,a_1^1]\colon \L^T(\S)\to \L^T(\S)$ is just $\langle a_1^1\rangle$. (In fact, for a degree 1 element the full strength of the elimination theorem is not necessary, since the bracket yields an isomorphism $\L^T_2(\S)\simeq \bwedge^2 \L^T_1(\S)$.) If we let $y$ denote the image of $d_A(\phi)$ in $\L^T_1(\S)=\pi_1(S,\ast)/\Gamma_2^T(\S)$, the fact that $[d_A(\phi),A\delta A^{-1}]\in \Gamma_3^T(\S)$ implies that $[y,a_1^1]=0\in \L^T_2(\S)$. 
By the elimination theorem, this implies $y=m\cdot a_1^1$ for some $m\in \Z$, as desired.

This shows that any lift $\widetilde{\phi}$ of $\phi$ satisfies $\widetilde{\phi}\in \Mod_{(3)}(\S_1)$. By Theorem~\ref{thm:basepointirrelevant} this is equivalent to $\widetilde{\phi}\in \Mod_{(3)}(\S_2)$, and by Theorem~\ref{thm:preserved} this implies $\pi_*(\widetilde{\phi})\in \Mod_{(3)}(\overline{S})=\K(\overline{S})$, as we claimed.
\end{proof}

\section{The Johnson homomorphism}\label{sec:Jh}
In this section and the next, we compute the quotient $\IS/\KS$. More specifically, we construct the \emph{partitioned Johnson homomorphism} $\tau_\S\colon \IS\to \Hom\big(\HS,\NS\big)$, and show that its kernel is exactly $\KS$. We also prove that the image of $\tau_\S$ is a certain explicitly defined subgroup $W_\S$ of $\Hom(\HS,\NS)$, so that we have a short exact sequence
\[1\to \KS\longrightarrow \IS\overset{\tau_\S}{\longrightarrow} W_\S\to 1.\]

\subsection{The classical Johnson homomorphism}\label{sec:review}
We briefly review Johnson's original construction of the Johnson
homomorphism for a surface $S$ with one boundary component, using
the action of $\Mod(S)$ on the universal two-step nilpotent quotient
of the free group $\pi\coloneq \pi_1(S,\ast)$, where $\ast\in\partial
S$. Let $\Gamma_j(\pi)$ denote the lower central series of the free group $\pi$. We have the
short exact sequence \begin{equation}\label{eq:origseq}1\to \Gamma_2(\pi)\to \pi\to H\to 1.\end{equation} Centralizing
the first term, we get the short exact sequence \[1\to N\to E\to H\to
1\] of Johnson~\cite{JohnsonAbelian}, where $H=\pi/\Gamma_2(\pi)$, $N=\Gamma_2(\pi)/\Gamma_3(\pi)$,
and $E=\pi/\Gamma_3(\pi)$.

Considering the exact sequence \eqref{eq:origseq} as a presentation for $H$, Hopf's
formula says that \[H_2(H;\Z)\iso\frac{\Gamma_2(\pi)\cap
  \Gamma_2(\pi)}{[\Gamma_2(\pi),\pi]}=\frac{\Gamma_2(\pi)}{\Gamma_3(\pi)}=N.\] Since $H$ is free abelian, 
$H_2(H;\Z)$ can be identified with $\bwedge^2 H$. This gives an isomorphism $\bwedge^2 H\iso
N$, which can be described explicitly as follows: $x\wedge y\in \bwedge^2
H$ is sent to $[\tilde{x},\tilde{y}]\in N$, where
$\tilde{x},\tilde{y}\in E$ are lifts of $x$ and $y$. All these
identifications are $\Mod(S)$--equivariant; in particular, $\I(S)$
acts trivially on $N$.

\begin{definition}[\cite{JohnsonAbelian}]
\label{def:classicalJh}
The Johnson homomorphism $\tau\colon \IS\to \Hom(H,N)$ is defined by assigning to $\phi\in \IS$ the homomorphism $\tau(\phi)\colon H\to N$ given by $x\mapsto [f(\tilde{x})\tilde{x}^{-1}]$, where $x\in \pi$ is any lift of $x\in H$. The fact that $N$ is central in $E$ implies that $\tau(\phi)$ is a homomorphism, and the fact that $\IS$ acts trivially on $H$ and on $N$ implies that $\tau$ is a homomorphism.
\end{definition}
Identifying $N$
with $\bwedge^2 H$ and $\Hom(H,N)$ with $H^\ast\otimes\bwedge^2 H$, Johnson proved in \cite{JohnsonAbelian}  that 
the image of this map is exactly $\bigwedge^3 H\leq
H^\ast\otimes\wedgetwo H$; compare with Theorem~\ref{thm:WS} below. The map
$\tau$ is $\Mod(S)$--equivariant with respect to the conjugation actions
on $\I(S)$ and on $\Hom(H,\wedgetwo H)$. The kernel $\ker
\tau\leq \I(S)$ is exactly the subgroup $\Mod_{(3)}(S)$ acting trivially on
$E=\pi/\pi_3$. Johnson proved in \cite{JohnsonII} that for a surface $S=S_{g,1}$
with only one boundary component, $\ker\tau$ is the group $\K(S)$ generated
by separating twists.

\subsection{The partitioned Johnson homomorphism}\label{sec:generaldef}
Our construction of the Johnson homomorphism for a general partitioned
surface follows Johnson's approach closely. For technical reasons, our definition of the partitioned Johnson homomorphism only makes reference to $\Shat$, not $\S$ itself. As a result, in remainder of Section~\ref{sec:Jh} the reader may assume that the partition on $\S$ is totally
separated, so that $\S=\Shat$, without contradiction. (This definition does not require Theorem~\ref{thm:preserved}, since the inclusion $\S\to \Shat$ is canonically defined; however, to say that $\ker \tau_\S$ coincides with $\KS$, as we prove in Corollary~\ref{cor:sameasKS}, does depend on Theorem~\ref{thm:preserved}.)

Let $\pihat\coloneq \pi_1(\Surhat,\ast)$. For a totally separated surface such as $\Shat$, we noted in Section~\ref{sec:charT} that $\L^T_1(\Shat)=\pihat/T(\Shat)$ is isomorphic to $\HS$, yielding an exact sequence
\begin{equation}\label{eq:Tseq}1\to T(\Shat)\to \pihat\to \HS\to 1.\end{equation}
Centralizing the first term, we get the exact sequence
\[1\to T(\Shat)/[T(\Shat),\pihat]\to \pihat/[T(\Shat),\pihat]\to \HS\to 1.\]
We define $\NS\coloneq \L^T_2(\Shat)$ and $\ES\coloneq \pihat/\Gamma_3^T(\Shat)$, so we may rewrite this exact sequence as
\begin{equation*}
  1\to \NS\to \ES\to \HS\to 1.\end{equation*}

\begin{definition}[The partitioned Johnson homomorphism]\label{def:tauS}
The Torelli group $\IS$ acts trivially on $\NS$ by Corollary~\ref{cor:TorLTtrivial}, and on $\HS$ by definition. Thus by the construction described in Definition~\ref{def:classicalJh}, the action of $\IS$ on
$\ES=\pihat/[T(\Shat),\pihat]=\pihat/\Gamma_3^T(\Shat)$ yields a homomorphism
\[\tau_\S\colon \IS\to\Hom\big(\HS,\NS\big)\] which we call the
(\emph{partitioned}) \emph{Johnson homomorphism}. It is given explicitly by
\[\tau_\S(\phi)(x)=[\phi(\widetilde{x})\widetilde{x}^{-1}]\] where $\widetilde{x}\in \pihat$ is a lift of $x\in \HS$.
\end{definition}

Note that $\tau_\S(\phi)=0$ if and only if $\phi$ acts trivially on $E(\S)=\pihat/\Gamma_3^T(\Shat)$. This implies that $\KS=\Mod_{(3)}(\S)$ is contained in $\ker \tau_\S$; we will show in Corollary~\ref{cor:sameasKS} that in fact $\KS=\ker \tau_\S$.

\para{$\Mod(S)$--equivariance of $\tau_\S$} The map $\tau_\S$ is $\Mod(S)$--equivariant, in the following sense. The mapping class
group $\Mod(S)$ acts on $\IS$ by conjugation. The action of $\Mod(S)$ on
$\HS$ and $\NS$ induces an action on $\Hom(\HS,\NS)$; to
be precise, if $f\in \Hom(\HS,\NS)$ and $\phi\in\Mod(S)$, the map $\phi_*f$
is defined by $\phi_* f(x)=\phi\big(f(\phi^{-1}(x))\big)$. The
following lemma is a formal consequence of the definition of $\tau_\S$.
\begin{lemma}
\label{lem:equivariance}
  Let $\phi\in\Mod(S)$ and $\psi\in\IS$. Then
  $\tau_\S(\phi\psi\phi^{-1})=\phi_*\tau_\S(\psi)$.
\end{lemma}
\begin{proof} From the definitions, we have
\begin{align*}
\tau_\S(\phi\psi\phi^{-1})(x)
&=[\phi\psi\phi^{-1}(\tilde{x})\cdot \tilde{x}^{-1}]\\
&=\phi\big[\psi\phi^{-1}(\tilde{x})\cdot \phi^{-1}(\tilde{x})^{-1}\big]\\
&=\phi\big(\tau_\S(\psi)(\phi^{-1}(x))\big)\\
&=\phi_*\tau_\S(\psi)(x)
\end{align*} for any $x\in \HS$.
\end{proof}

\begin{remark}
If $\eta$ is a $\P$--separating curve in $\S$, then $\eta$ is a separating curve in $\Sbar$; it follows that $\eta$ is trivial in $\HS\iso H(\Sbar)$. Thus for any loop $\gamma\in \pi_1(S,\ast)$ representing $\eta$, we have $\gamma\in \Gamma_2^T(\S)$. Moreover, any other such loop $\gamma'$ is conjugate to $\gamma$; since $\Gamma_j^T(\S)$ is a central filtration, they have the same image $[\gamma]=[\gamma']\in \Gamma_2^T(\S)/\Gamma_3^T(\S)=\L^T_2(\S)$. Thus any $\P$--separating curve $\eta$ represents a well-defined class $[\eta]\in \NS$.\end{remark}

\subsection{Action on arcs}\label{sec:arcs}
In this section we use the action of $\IS$ on arcs connecting different boundary components to construct a family of abelian quotients $d_i$ of $\IS$.

\para{The abelian quotients $d_i$} Recall from Section~\ref{sec:Johnsonfiltproper} that for any arc $A$ from the basepoint $\ast$ to a boundary component of $\Surhat$, and any $\phi\in \Mod(\Surhat)$, we denote by $d_A(\phi)$ the element $\phi(A)A^{-1}\in \pihat$. Moreover, if $A'$ is another arc ending at the same boundary component, $d_{A'}(\phi)$ differs from $d_A(\phi)$ by $\phi(x)x^{-1}[x,d_A(\phi)]$, where $x=A'A^{-1}\in \pihat$. Note that  $[x,d_A(\phi)]$ always lies in $T(\Shat)$ while $\phi(x)x^{-1}\in T(\Shat)$ if $\phi\in \IS$. Thus for $\phi\in \IS$, the class of $d_A(\phi)$ in $\pihat/T(\Shat)=\HS$ does not depend on the arc $A$, only on the boundary component.

For a surface $\S=(S,\P,\ast)$ with $\P=\{P_0,P_1,\ldots,P_k\}$, for each $1\leq i\leq k$ we let $z_i$ denote the boundary component of $\Surhat$ corresponding to the block $P_i\in \P$. This boundary component is separating, and its class in $\NS$ is precisely the generator $z_i$ from Proposition~\ref{prop:cS} (justifying the slight abuse of notation).
We define the homomorphism $d_i\colon \IS\to \HS$ by \[d_i(\phi)=[\phi(A_i)A_i^{-1}]\in
H(\S),\] where $A_i$ is an arc from the basepoint $\ast$ to the boundary component $z_i$.

One important difference between $\IS$ and later terms in the Johnson filtration is that these maps $d_i$ factor through the Johnson homomorphism $\tau_\S$. This property is stated below as Theorem~\ref{thm:deltas}, and will be proved in Section~\ref{sec:restrict}.

\para{The homomorphisms $\delta_i$} For each $1\leq i\leq k$ we have a map $\NS\to \Z$, obtained by counting intersections with the arc $A_i$ (with multiplicity). Specifically, if $y\in \NS$ is represented by $\widetilde{y}\in \pihat$, let $(y,A_i)$ denote the algebraic intersection number of $\widetilde{y}$ with the arc $A_i$. This does not depend on the choice of arc $A_i$, since two such arcs from $\ast$ to $z_i$ differ by a cycle plus a multiple of $z_i$. Any cycle has trivial intersection with any commutator or any boundary component $z_j$, as does the boundary component $z_j$. By Lemma~\ref{lem:TSgen}, $\NS$ is generated by $[\pihat,\pihat]$ and the elements $z_j$, so this suffices.
\begin{definition}\label{def:tildedelta} For each $1\leq i\leq k$, the
homomorphism \[\delta_i\colon
\Hom\big(\HS,\NS\big)\to\HS\] is defined by the adjunction
\begin{equation}\big(f(x),A_i\big)=\big(x,\delta_i(f)\big).\end{equation}\end{definition}
 The following proposition will be proved in Section~\ref{sec:restrict}.
\begin{theorem}\label{thm:deltas} The maps $d_i\colon \IS\to \HS$ factor through $\tau_\S$; more specifically, we have $d_i=\delta_i\circ\tau_\S$.
\end{theorem}

\begin{corollary}
\label{cor:sameasKS} The kernel $\ker \tau_\S$ of the partitioned Johnson homomorphism $\tau_\S$ is precisely the subgroup $\KS=\MS{3}$.
\end{corollary}
\begin{proof} By Theorem~\ref{thm:preserved} we have $\KS=\Mod(S)\cap \K(\Shat)$, so it suffices to prove the corollary for $\Shat$; in other words, we may assume that $\S$ is totally separated.  By definition, $\ker\tau_\S$ consists of those $\phi\in \Mod(S)$ acting trivially modulo $\Gamma_3^T(\S)$. When $\S$ is totally separated, condition (iii) of Definition~\ref{def:Johnsonfiltration} is vacuous. Thus $\MS{3}$ consists of those $\phi\in \Mod(S)$ satisfying conditions (i) and (ii); i.e.\ those $\phi\in \ker\tau_\S$ that additionally satisfy $d_{A_i}(\phi)\in \Gamma_2^T(\S)$ for each $1\leq i\leq k$. Since $d_i(\phi)$ is the class of $d_{A_i}(\phi)$ in $\HS=\pi_1(S)/\Gamma_2^T(\S)$, condition (ii) asks that $d_i(\phi)=0$ for all $1\leq i\leq k$.  But Theorem~\ref{thm:deltas} states that $\tau_\S(\phi)=0 \implies d_i(\phi)=0$. We conclude that any $\phi\in \ker\tau_\S$ lies in $\KS=\MS{3}$, as desired.
\end{proof}

\subsection{The image of $\tau_\S$}\label{sec:image}
The definition of $\HS=\L_1^T(\Shat)$, of $\NS=\L_2^T(\Shat)$, and of $\tau_\S$ do not really depend on the partitioned surface $\S$ itself, only on $\Shat$. However the image of $\tau_\S$ certainly depends on $\S$. In this section we describe certain conditions on the image of $\tau_\S$, which together cut out a subspace $W_\S$ of $\Hom(\HS,\NS)$. We will eventually prove in Section Section~\ref{sec:compute} that $W_\S=\im \tau_\S$.

\para{Understanding $\Hom(\HS,\NS)$}
There is a natural quotient $\NS\twoheadrightarrow \bwedge^2 \HS$ defined by sending $z_i\mapsto 0$ for all $1\leq i\leq k$, induced for example by the inclusion $\Shat\hookrightarrow \overline{\S}$ of Section~\ref{sec:Torelli}. In fact, it follows from Theorem~\ref{thm:LTfree} that \begin{equation}\label{eq:NS}N(\S)\simeq\bwedge^2 \HS\oplus \Z^k.\end{equation} where $\bwedge^2 \HS$ is the image of $[\pihat,\pihat]$ and the $\Z^k$ factor is spanned by $z_1,\ldots,z_k$. Note that the intersection $y\mapsto (y,A_i)$ vanishes on $\bwedge^2 \HS$ and satisfies $(z_j,A_i)=\delta_{ij}$.

The projection $\NS\twoheadrightarrow \bwedge^2 \HS$ induces a projection:
\begin{equation}\label{eq:Hw2H}
  \Hom(\HS,\NS)\twoheadrightarrow \Hom(\HS,\bwedge^2 \HS)\simeq \HS\otimes \bwedge^2\HS
\end{equation}
Note that $\bwedge^3
\HS$ embeds into $\HS\otimes \bwedge^2 \HS$ as the ``Jacobi identity'':
\[a\wedge b\wedge c\mapsto a\otimes b\wedge c+b\otimes a\wedge
c+c\otimes a\wedge b\]

\para{The subspaces $D_i\leq \DS\leq \HS$} We denote by $D(\S)\leq \HS$ the isotropic subspace spanned by the homology
classes $a_c$ of all the boundary components $c\in \pi_0(\partial S)$. Similarly, for the single block $P_i\in \P$, we denote by $D_i\leq \DS$ the
subspace spanned by those components $a_c$ for $c\in P_i$. Note that $\DS^\perp$ is exactly the subspace of $\HS$ spanned by $H_1(S)$.
\begin{definition}
\label{def:WS}
The subspace $W_\S\leq \Hom(\HS,\NS)$ consists of
  those elements\linebreak ${f\colon \HS\to \NS}$ satisfying the
  following conditions:
\begin{enumerate}[(I)]
\item the image in $\HS\otimes \bwedge^2 \HS$ of $f$ under the projection \eqref{eq:Hw2H} is contained in the subspace $\bwedge^3 \HS\leq \HS\otimes \bwedge^2 \HS$.
\item for $1\leq i\leq k$ we have $\delta_i(f)\in \DS^\perp$ and furthermore for any $a\in
  D_i$, $f(a)=\delta_i(f)\wedge a$.
\item for any $a\in D_0$, $f(a)=0$.
\end{enumerate}
\end{definition}

The following characterization of the image of $\tau_\S$ is one of the main theorems of the paper.
\begin{theorem}\label{thm:WS}
  $W_\S=\im\tau_\S$.
\end{theorem}
We  prove  that $\im\tau_\S$ is contained in $W_\S$ in the next section (Theorem~\ref{thm:WScontained}). We defer the remaining direction of Theorem~\ref{thm:WS} until Section~\ref{sec:compute}, where we first  compute the value of $\tau_\S$ on various fundamental elements of $\IS$, then use these computations to prove that $\tau_\S$ surjects to $W_\S$. Before moving on, we deduce Theorem~\ref{thm:betti} from Theorem~\ref{thm:WS}.

\begin{proof}[Proof of Theorem~\ref{thm:betti}]
Let $\Gamma$ be an arbitrary finite-index subgroup of $\IS$ satisfying $\KS<\Gamma$. 
Let $\iota\colon \S\hookrightarrow \overline{S}$ be the standard inclusion into the surface $\overline{S}$, which has one boundary component.
Putman~\cite[Theorem~1.2]{PutmanKS} states that as long as $\S$ has genus at least 3, for any such $\Gamma$ we have \[H_1(\IS;\Q)\iso H_1(\Gamma;\Q)\iso \tau_{\overline{S}}(\iota_*(\IS))\otimes \Q.\]  The naturality of $\tau_\S$, proved in Theorem~\ref{thm:commutes} below, means that $\tau_{\overline{S}}(\iota_*(\IS))$ is isomorphic to $\tau_\S(\IS)$, and Theorem~\ref{thm:WS} states that $\tau_\S(\IS)$ is precisely the subspace $W_\S<\Hom(\HS,\NS)$. Thus Putman's theorem implies that the first Betti number $b_1(\IS)=b_1(\Gamma)$ is the dimension $\dim_\Q W_\S\otimes \Q$.

Consider the projection $\Hom(\HS,\NS)\twoheadrightarrow \HS^{\oplus k}$ defined by $f\mapsto \big(\delta_1(f),\ldots,\delta_k(f)\big)$. Condition (II) of Definition~\ref{def:WS} states that when restricted to $W_\S$, this projection
   has image $(\DS^\perp)^{\oplus k}$. It is easy to construct elements of $W_\S$ surjecting to $(\DS^\perp)^{\oplus k}$, and indeed we will do this by hand in proving Theorem~\ref{thm:WS}.
The kernel of this surjection $W_\S\twoheadrightarrow (\DS^\perp)^{\oplus k}$ consists of those $f\in W_\S$ satisfying $\delta_i(f)=0$ for all $i$. This implies that $f$ lies in the subspace $\Hom(\HS,\bwedge^2 \HS)$ of $\Hom(\HS,\NS)$, so by condition (I) we have $f\in \bwedge^3
  \HS<\Hom(\HS,\bwedge^2 \HS)$. Conditions (II) and
  (III) imply that $f$ satisfies $f(a)=0$ whenever $a\in
  \DS$. As an element of $\bwedge^3 \HS$, this means that $f$ is
  contained in \[\big(\DS^\perp\otimes \bwedge^2\HS\big)\cap\bwedge^3 \HS\simeq\bwedge^3 \DS^\perp.\]
Since any $f\in \bwedge^3 \DS^\perp$ clearly satisfies conditions (II) and (III), we conclude that 
$W_\S$ fits into a short exact sequence
  \begin{equation}
  \label{eq:WSses}
  0\to \bwedge^3 \DS^\perp\to W_\S\to (\DS^\perp)^{\oplus k}\to 0.
  \end{equation}
Let $n=\abs{\pi_0(\partial S)}$ be the number of boundary components of $S$, partitioned into $\abs{\P}=k+1$ blocks, and let $D\coloneq 2g+n-k-1$. As abelian groups, we have $H_1(S;\Z)\iso \Z^{2g+n-1}$, $\DS\iso \Z^{n-k-1}$, and $\DS^\perp\iso \Z^D$.
Since $\DS^\perp$ is torsion-free, \eqref{eq:WSses} splits, so $\dim_\Q W_\S\otimes \Q=\binom{D}{3}+D^k$, as desired. \end{proof}

\subsection{Restricting the image of $\tau_\S$}
\label{sec:restrict}
\begin{theorem}\label{thm:WScontained}
The map $\tau_\S$ has image contained in $W_\S$.
\end{theorem}
\begin{proof}[Proof of Theorem~\ref{thm:deltas} and Theorem~\ref{thm:WScontained}.] Consider
$\phi\in\IS$ and let $f=\tau_\S(\phi)\in\Hom(\HS,\NS)$. We first show
that $f$ satisfies condition (I) of Definition~\ref{def:WS}, and at the same time  verify Theorem~\ref{thm:deltas}.

  We will make use of the following identities (where ${}^ab$ denotes
  conjugation):
  \[[aw,b]={}^a[w,b]\cdot[a,b], \qquad\qquad[a,bw]=[a,b]\cdot{}^b[a,w].\] If
 $w\in \Gamma^T_2(\Shat)=T(\Shat)$, both $[w,b]$ and $[a,w]$ are in $\Gamma^T_3(\Shat)$, so we have
  \[[aw,b]\equiv[a,b][w,b]\bmod{\Gamma^T_4(\Shat)},
  \qquad[a,bw]\equiv[a,b][a,w]\bmod{\Gamma^T_4(\Shat)}.\]

Let $\zeta_0\in \pihat$ represent the boundary component $z_0$ of $\Surhat$ that contains the basepoint. The key to our proof is to consider the action of $\phi$ on $\zeta_0$\,---\,even though $\zeta_0$ is contained in the boundary $\partial \Surhat$ and so the action of $\phi$ on $\zeta_0$ is trivial. Choose a basis $\{\alpha_i,\beta_i\}\cup\{\zeta_1,\ldots,\zeta_k\}$ of $\pihat$ as in Proposition~\ref{prop:cS}, so that $\{\alpha_i,\beta_i\}$ descends to the symplectic basis $\{a_i,b_i\}$ of $\HS$, each $\zeta_i$ represents $z_i\in \NS$, and 
  $\zeta_0=[\alpha_1,\beta_1]\cdots[\alpha_g,\beta_g]\zeta_1\cdots\zeta_k$.

  In this proof, let $\eta_\phi(x)\coloneq x^{-1}\phi(x)$. By definition, $\phi\in \IS$ if and only if  $\eta_\phi(x)\in\Gamma^T_2(\Shat)$ for all $x\in \pihat$.   Our fundamental computation \eqref{eq:fund} gives $\phi(\zeta_j)=[d_j(\phi),\zeta_j]\zeta_j$. Recalling that  $\zeta_j\in \Gamma^T_2(\Shat)$ for all $j$, we calculate:
\begin{align*}
  \phi(\zeta_0)&= \phi\left(\prod_{i=1}^g [\alpha_i,\beta_i]\prod_{j=1}^k\zeta_j\right)\\
  &=\prod_{i=1}^g[\phi(\alpha_i),\phi(\beta_i)]\prod_{j=1}^k\phi(\zeta_j)\\
  &=\prod_{i=1}^g[\alpha_i\eta_\phi(\alpha_i),\beta_i\eta_\phi(\beta_i)]
       \prod_{j=1}^k[d_j,\zeta_j]\zeta_j\\
  &\equiv \prod_{i=1}^g[\alpha_i,\beta_i] [\alpha_i,\eta_\phi(\beta_i)]
    [\eta_\phi(\alpha_i),\beta_i] \prod_{j=1}^k[d_j,\zeta_j]\zeta_j \bmod{\Gamma^T_4(\Shat)}\\
  &\equiv
  \left(\prod_{i=1}^g[\alpha_i,\eta_\phi(\beta_i)]
   [\eta_\phi(\alpha_i),\beta_i] \prod_{j=1}^k[d_j,\zeta_j]\right)\cdot\zeta_0\bmod{\Gamma^T_4(\Shat)}
\end{align*} Define \[X=\prod_{i=1}^g[\alpha_i,\eta_\phi(\beta_i)]
   [\eta_\phi(\alpha_i),\beta_i] \prod_{j=1}^k[d_j,\zeta_j]\in \Gamma^T_3(\Shat)\] and consider the class $[X]\in\L^T_3(\Shat)=\Gamma^T_3(\Shat)/\Gamma^T_4(\Shat)$. 
Note that $\eta_\phi(\alpha_i)$ and $\eta_\phi(\beta_i)$ in $\Gamma^T_2(\Shat)$ represent
$\tau(\phi)(a_i)$ and $\tau(\phi)(b_i)$ in $\NS$. 
Thus the following element $Y\in \HS\otimes \NS$ descends to $[X]\in \L^T_3(\Shat)$ under the commutator bracket:
\begin{equation}\label{eq:Y}Y\coloneq\sum_{i=1}^g a_i\otimes \tau_\S(\phi)(b_i)-b_i\otimes \tau_\S(\phi)(a_i)+\sum_{j=1}^k d_j(\phi)\otimes z_j\ \in \HS\otimes \NS\end{equation}
Note that the first summation is exactly the expansion of $-\tau_\S(\phi)\in\Hom(\HS,\NS)$ in $\HS\otimes \NS$. Indeed, since $\{a_i,b_i\}$ form a basis of $\HS$, we can write \[\tau_\S(\phi)=\sum_{i=1}^g a_i^*\otimes \tau_\S(\phi)(a_i)+b_i^*\otimes \tau_\S(\phi)(b_i)\in \HS^*\otimes \NS\] and under the isomorphism $\HS\simeq \HS^*$ we have $a_i^*=b_i$ and $b_i^*=-a_i$ (since $a_i^*=(\cdot,b_i)$ and $b_i^*=(\cdot,-a_i)$).
In particular, it follows from the discussion following \eqref{eq:NS} that the coefficient of $z_j$ in the first summation is $-\delta_j(\tau_\S(\phi))\otimes z_j$. 

We calculated above that $\eta_\phi(\zeta_0)\equiv X \bmod{\Gamma_4^T(\Shat)}$. However, since $\zeta_0$ is contained in the boundary of $\Surhat$, we have $\phi(\zeta_0)=\zeta_0$. Thus $\eta_\phi(\zeta_0)$ is trivial, so $X$ must lie in $\Gamma_4^T(\Shat)$. In other words, $[X]=0\in \L^T_3(\Shat)$; this implies that $Y\in \HS\otimes \NS$ is contained in the kernel of the map $\HS\otimes \NS\to \L_3^T(\Shat)$. We now recall Proposition~\ref{prop:Tbracket}, which states that the bracket $\L_1^T(\Shat)\otimes\L_2^T(\Shat)\to \L_3^T(\Shat)$ induces the short exact sequence \[1\to \bwedge^3 \HS\to \HS\otimes
\NS\to \L^T_3(\Shat)\to 1.\] This has the following implications. 
First, we saw above that the coefficient of $z_j$ in $Y$ is \[\big(d_i(\phi)-\delta_i(\tau_\S(\phi))\big)\otimes z_j.\] But the factor $\HS\otimes \Z^k<\HS\otimes \NS$ intersects $\bwedge^3 \HS$ trivially, so for $Y$ to be contained in $\bwedge^3 \HS$ we must have $d_i(\phi)-\delta_i(\tau_\S(\phi))=0$. This completes the proof of Theorem~\ref{thm:deltas}.

Furthermore, this implies that $Y\in \HS\otimes \NS$ is the difference of $-\tau_\S(\phi)$ and its projection to the $\Z^k$ factor spanned by the $z_j$; in other words, $Y\in \HS\otimes\bwedge^2 \HS$ is the projection of $-\tau_\S(\phi)$ to $\Hom(\HS,\bwedge^2\HS)$ under the map \eqref{eq:Hw2H}.
Now Proposition~\ref{prop:Tbracket} states that this is contained in $\bwedge^3 \HS$, verifying condition (I) of Definition~\ref{def:WS}.\\

Now we can use Theorem~\ref{thm:deltas} to verify conditions (II) and (III). We have just proved that
$\delta_i(f)=d_i(\phi)$. But from the definition of $d_i$ we can see that $d_i(\phi)$ is contained in $\DS^\perp$. Indeed, we have $d_i(\phi)=[\phi(A_i)A_i^{-1}]$. Since $\phi$ is the identity outside $S$,
we have $\phi(A_i)=A_i$ outside $S$, and thus $[\phi(A_i)A_i^{-1}]$ may be
represented by a cycle lying inside $S$. We observed in Section~\ref{sec:image} that the
span of $H_1(S)$ in $\HS$ is exactly $\DS^\perp$, and so we conclude
that $\delta_i(f)=d_i(\phi)\in \DS^\perp$ for all $i$.

To verify the remainder of condition (II) we must show that $\tau_\S(\phi)(a)=d_i(\phi)\wedge a$ for $a\in D_i$. It suffices to check
this when $a$ is the class of a single boundary component in $P_i$. By our fundamental computation \eqref{eq:fund} we have $\phi(\alpha)\alpha^{-1}=[d_i(\phi),\alpha]$, where $\alpha\in \pihat$ represents $a\in D_i$, so $\tau_\S(\phi)(a)=d_i(\phi)\wedge a$ as desired.
A similar argument verifies condition (III). If $a$ is the class of a
single boundary component contained in $P_0$, then it is connected
outside $S$ to the basepoint; in particular, $a$ can be represented by
a loop disjoint from ${S}$. Thus $\phi$ fixes this loop
pointwise, and so $\tau_\S(\phi)(a)=0$. This completes the proof of Theorem~\ref{thm:WScontained}.
\end{proof}
\begin{remark}
For any term in the Johnson filtration, we could define homomorphisms $d_i$ from $\MS{k}$ to $\L_{k-1}^T(\Shat)$ by $d_i(\phi)=[d_{A_i}(\phi)]\in \L^T_{k-1}(\Shat)$.
An argument similar to the above would show that these maps $d_i$ are controlled to some degree by the action of $\MS{k}$ on $\pihat/\Gamma_{k+1}^T(\S)$. However, to show that the $d_i$ actually factor through this action when $k=2$ we needed Proposition~\ref{prop:Tbracket}, which states that the various maps $\L_1^T(\S)\to \L_3^T(\S)$ defined by $x\mapsto [x,z_j]$ have independent images. The corresponding statement for maps $\L_{k-1}^T\to \L_{k+1}^T$ is definitely false for $k>2$: for example, we have relations such as $[z_i,z_j]+[z_j,z_i]=0$ or $[[z_i,z_j],z_k]+[[z_j,z_k],z_i]+[[z_k,z_i],z_j]=0$.
\end{remark}

\subsection{Naturality of $\tau_\S$} \label{sec:naturality}
In this section we show that the partitioned Johnson homomorphism $\tau_\S$ is natural. This means that for any morphism
$\iota\colon \S\to \S'$,
we must define a map \[\iota_*\colon
W_{\S}\to W_{\S'}\] so
that $\tau_{\S'}(\iota_*\phi)=\iota_* \tau_\S(\phi)$ for all $\phi\in \IS$; in other words,  the following diagram commutes:
\begin{equation*}
\xymatrix{
    \I(\S)\ar^{\qquad\tau_{\S}\qquad}[rr]\ar_{\qquad\iota_*}[d]
    &&W_{\S}\ar^{\iota_*}[d]\\
    \I(\S')\ar_{\qquad\tau_{\S'}\qquad}[rr]&&W_{\S'}}\quad\ \ 
\end{equation*}
To define the map $\iota_*\colon W_{\S}\to W_{\S'}$, we consider separately the cases when $\iota$ is non-collapsing and when $\iota$ is a simple capping; since every morphism is a
composition of such inclusions, this suffices.

\begin{definition}[Alternate notation for $W_\S$]
\label{def:alternate}
The map $\iota_*\colon W_\S\to W_{\S'}$ is somewhat unwieldy when expressed in terms of the usual notation for $W_\S<\Hom(\HS,\NS)$. However there is an equivalent way to describe a basis of $W_\S$, and with respect to this basis the definition of the map $\iota_*$ is very simple. Recall from Definition~\ref{def:WS} that $W_\S$ can be thought of as a subspace of $\bwedge^3 \HS\oplus (\HS\otimes \Z^k)$. For the $\bwedge^3 \HS$ component, we have the standard
$a\wedge b \wedge c$ notation. For the $\HS\otimes \Z^k$ factor we write $x\wedge z_i$ as a formal expression for the element $x\otimes z_i$. We remark that in this notation, the homomorphism $\delta_i\colon W_\S\to \HS$ has the form
\[\delta_i\Big(\sum_j a_j\wedge b_j\wedge c_j + \sum_{i=1}^k x_i\wedge z_i\Big)=x_i,\] so we will often write $\delta_i(f)\wedge z_i$ for the components of this second factor $\HS\otimes \Z^k$.
\end{definition}

\begin{definition}
\label{def:WSiota}
For a simple capping  $\iota\colon \S\to \S'$, if $\S'$ is obtained from $\S$ by attaching a disk to the separating component $z_i$, we define $\iota_*$ as follows.
\begin{align*}\label{eq:naturality_refined}
\qquad\qquad \iota_*\colon W_{\S}\ &\to\  W_{\S'}\\
a\wedge b\wedge c\ &\mapsto\  a\wedge b\wedge c\\
x\wedge z_i\ &\mapsto\  0\\
x\wedge z_j\ &\mapsto\ x\wedge z_j \quad\text{for }j\neq i
\end{align*}

For a non-collapsing morphism $\iota\colon \S\to \S'$, decompose $\Shat'\setminus
{\Shat}$ into subsurfaces $U_i$, so that $U_i$
meets $\Shat$ in the component $z_i$ corresponding to the block $P_i$. We can consider $H(U_i)$ as a subspace of $H(\S')$. Let
$\omega_{U_i}\in\bwedge^2 H(U_i)$ represent the intersection form on
$H(U_i)$, and let $z_i^1,\ldots,z_i^l$ be the boundary
components of $U_i$ that are disjoint from $\Shat$. We define $\iota_*$ as follows.
\begin{align*}
\qquad\qquad \iota_*\colon W_{\S}\ &\to\  W_{\S'}\\
a\wedge b\wedge c\ &\mapsto\  a\wedge b\wedge c\\
x\wedge z_i\ &\mapsto\  x\wedge (\omega_{U_i}+z_i^1+\cdots+z_i^l)
\end{align*}
Note that in $N(U_i)$ we have $z_i=\omega_{U_i}+z_i^1+\cdots+z_i^l$, since $z_i$ is the boundary component that would be $z_i^0$, but traversed in the opposite direction.
\end{definition}

\begin{theorem} \label{thm:commutes} For any inclusion $\iota\colon
  \S\to \S'$, the map $\iota_*\colon W_\S\to W_{\S'}$ defined in Definition~\ref{def:WSiota} satisfies $\tau_{\S'}(\iota_*\phi)=\iota_*\tau_\S(\phi)$ for all $\phi\in \IS$. \end{theorem}
\begin{proof}
The argument is similar to the proof of Theorem~\ref{thm:preserved}. In the course of the proof, we will extend $\iota_*\colon W_{\S}\to W_{\S'}$ to a map $\iota'_*\colon\Hom(H(\S),N(\S))\to \Hom(H(\S'),N(\S'))$. We then verify naturality for $\iota'_*$ directly from the definition of $\tau_\S$, and it remains only to check that $\iota'_*$ does in fact restrict to $\iota_*$ on $W_\S$.

 First, consider the case when
$\iota$ is a simple capping, so that $\S'$ is obtained from $\S$ by
attaching a disk to the separating component $z_j$. Since $z_j$ was
separating and thus represented by $\zeta_j\in T(\Shat)$, we have
$H(\S')=H(\S)$, and the natural map $N(\S)\to N(\S')$ is
surjective with kernel generated by $z_j$.
Thus the exact sequences defining $\tau_{\S}$ and $\tau_{\S'}$ are related by the following diagram \eqref{eq:NEdiagram}:

\begin{equation}
\label{eq:NEdiagram}
\begin{split}\xymatrix{
    &1\ar[d]&1\ar[d]&&\\
    &\Z\ar[d]&\Z\ar[d]&&\\
    1\ar[r]&N(\S)\ar[r]\ar[d]&E(\S)\ar[r]\ar[d]
    &H(\S)\ar[r]\ar@{=}[d]&1\\
    1\ar[r]&N(\S')\ar[d]\ar[r]&E(\S')\ar[d]\ar[r]&H(\S')\ar[r]&1\\
    &1&1&&
  }\end{split}\end{equation}
It follows that $\tau_{\S'}(\iota_*\phi)=\iota'_* \tau_{\S}(\phi)$ for $\phi\in \I(\S)$, where $\iota'_*$ is the map \[\iota'_*\colon\Hom(H(\S),N(\S))\to \Hom(H(\S'),N(\S'))\] sending $f\in\Hom(H(\S),N(\S))$ to the
composition
\[\iota'_*f\colon H(\S')\iso H(\S)
\overset{f}{\longrightarrow}N(\S)\twoheadrightarrow N(\S').\] 
The restriction of $\iota'_*$ to $W_{\S}$ has kernel $\HS\otimes z_j$ and thus coincides with $\iota_*$ as defined in Definition~\ref{def:WSiota}, verifying the theorem in this case.\\

Now, consider the case when $\iota$ is non-collapsing. In this case the induced
map $\iota_*\colon N(\S)\to N(\S')$ is an injection by Theorem~\ref{thm:LTfree}. Recall that the $U_i$ are the components of
$\Shat'\setminus {\Shat}$, appropriately labeled so that $U_i$
intersects $\Shat$ in the boundary component $z_i$. Each $U_i$ inherits the structure of a partitioned surface from $\Shat'$, and since the resulting partition is totally separated,
the inclusion morphism $U_i\to \Shat'$ is uniquely determined. Identifying
$H(\S)$ with its image in $H(\S')$, we have an orthogonal splitting
$H(\S')=H(\S)\oplus \bigoplus_i H(U_i)$. We use this to define $\iota'_*$ as follows.
Given a homomorphism $f\in \Hom(H(\S),N(\S))$,
let $\iota'_*f\in \Hom(H(\S'),N(\S'))$ be the homomorphism defined by:
\begin{equation}\label{eq:induced}\iota'_*f(x)=
\begin{cases}\iota_*\big(f(x)\big)&\text{ if }x\in H(\S)\\
\delta_i(f)\wedge x&\text{ if }x\in H(U_i)
\end{cases}
\end{equation} The term $\delta_i(f)\wedge x$ here is taken
in $\bwedge^2 H(\S')\leq N(\S')$.

To verify that $\tau_{\S'}(\iota_*\phi)=\iota'_*\tau_{\S}(\phi)$, we consider two cases separately. If $x\in \HS$, we can represent it by an element $\gamma=C\xi C^{-1}\in \pi_1(\Surhat',\ast)$, where $C$ is an arc from $\ast'$ to $\ast$ in $\Surhat'\setminus \Surhat$, and $\xi$ is an element of $\pi_1(\Surhat,\ast)$ representing $x\in H(\S)$. As in the proof of Theorem~\ref{thm:preserved} we compute $(\iota_*\phi)(\gamma)\gamma^{-1}=C\phi(\xi)\xi^{-1}C^{-1}=i_*(\phi(\xi)\xi^{-1})$. It follows that $\tau_{\S'}(\iota_*\phi)(x)=\iota_*(\tau_\S(x))$, as claimed.
If $x\in H(U_i)$, we can represent it by an element $\gamma=CA\delta A^{-1}C^{-1}$, where $\delta$ is a loop in $\pi_1(U_i)$ and $A$ is an arc in $\Surhat$ from $\ast$ to the boundary component $z_i$. Our fundamental computation \eqref{eq:fund} shows as before that $(\iota_*\phi)(\gamma)\gamma^{-1}=[C d_A(\phi) C^{-1},\gamma]$. Thus $\tau_{\S'}(\iota_*\phi)(x)= d_i(\phi)\wedge x$, and by Theorem~\ref{thm:deltas} this is $\delta_i(\tau_\S(\phi))\wedge x$. Thus for any $x\in H(\S')$ we have $\tau_{\S'}(\iota_*\phi)(x)=\iota'_*(\tau_{\S}(\phi))(x)$.

It remains only  to check that $\iota'_*$ agrees on $W_\S$ with the map $\iota_*$  in Definition~\ref{def:WSiota}. For $x\in \HS$ we have $\iota'_*f(x)=\iota_*(f(x))$, while for $x\in H(U_i)$ we have $\iota'_*f(x)=\delta_i(f)\wedge x$. Considering $\Hom(H(U_i),N(\S'))$ as a subspace of $\Hom(H(\S'),N(\S'))$, if $\{a_i^j,b_i^j\}$ form a symplectic basis for $H(U_i)$, the restriction of $\iota'_*f$ to $H(U_i)$ is $\sum_j a_i^j\otimes b_i^j\wedge \delta_i(f)+b_i^j\otimes \delta_i(f)\wedge a_i^j$. Thus if $f=\sum_i x_i\wedge y_i\wedge z_i+\sum_i \delta_i(f)\otimes z_i$, we have
\[\iota_*f = \sum_i x_i\wedge y_i\wedge z_i + \sum_i \left(\delta_i(f)\otimes z_i+\sum_j a_i^j\otimes b_i^j\wedge \delta_i(f)+b_i^j\otimes \delta_i(f)\wedge a_i^j\right).\]
We noted above that $z_i=\sum_j a_i^j\wedge b_i^j + \sum_k z_i^k$ in $N(\S')$, so the parenthesized term can be written as: \[\delta_i(f)\otimes (\sum_j a_i^j\wedge b_i^j+\sum_k z_i^k)+\sum_j a_i^j\otimes b_i^j\wedge \delta_i(f)+b_i^j\otimes \delta_i(f)\wedge a_i^j,\]
which is precisely $\delta_i(f)\wedge (\sum_j a_i^j\wedge b_i^j+\sum_k z_i^k)$. Thus the restriction of $\iota'_*$ to $W_\S$ is the map induced by $x\wedge y\wedge z\mapsto x\wedge y\wedge z$ and $\delta_i(f)\wedge z_i\mapsto \delta_i(f)\wedge (\sum_j a_i^j\wedge b_i^j+\sum_k z_i^k)$, namely $\iota_*$.
\end{proof}

\para{Moving the basepoint} We saw in Theorem~\ref{thm:basepointirrelevant} that the Johnson filtration of a partitioned surface $\S=(S,\P,\ast)$ does not depend on the location of the basepoint $\ast$, only on the partition $\P$. The Johnson homomorphism $\tau_\S\colon \IS\to W_\S$ does in fact depend on the basepoint, but in a controlled way. Let $\S'=(S,\P,\ast')$ be the same partitioned
surface, except that the basepoint $\ast'$ lies in $P_i\in
\P$ instead of $P_0$, and let $\omega\in \bwedge^2 \HS\iso\bwedge^2 H(\S')$ represent the symplectic form.
\begin{lemma}\label{lem:bpchange}
  If $\S$ and $\S'$ coincide as partitioned surfaces except that
  $\ast'\in P_i$, then
\[\tau_{\S'}(\phi)=\tau_{\S}(\phi)-d_i(\phi)\wedge \omega.\]
\end{lemma}
\begin{proof}
As in the proof of Theorem~\ref{thm:basepointirrelevant}, let $\overline{\xi}$ denote $A_i^{-1}\xi A_i$, where $A_i$ is an arc from $\ast$ to $\ast'$. The assignment $\xi\mapsto \overline{\xi}$ defines an isomorphism $\pi_1(\Surhat,\ast)\iso \pi_1(\Surhat,\ast')$, and thus induces isomorphisms $\HS\iso H(\S')$ and $\NS\iso N(\S')$. We saw in \eqref{eq:basechange} that $\phi(\overline{\xi})\overline{\xi}^{-1}=\overline{[d_i(\phi)^{-1},\phi(\xi)]\phi(\xi)\xi^{-1}}$, which shows that  \[\tau_{\S'}(\phi)(x)=\tau_\S(\phi)(x)-d_i(\phi)\wedge x.\] The homomorphism $x\mapsto -d_i(\phi)\wedge x$ in $\Hom(\HS,\bwedge^2\HS)$ is represented by $-d_i(\phi)\otimes \omega\in \HS\otimes \bwedge^2\HS$, so $\tau_{\S'}(\phi)=\tau_\S(\phi)-d_i(\phi)\wedge\omega\in W_{\S}\iso W_{\S'}$, as desired.
\end{proof}

\para{Viewing $\tau_\S$ as a natural transformation} One way to phrase the naturality of $\tau_\S$ proved in Theorem~\ref{thm:commutes} is to say that $\tau_\S$ is a natural
transformation.  
We have already noted in Section~\ref{sec:Torellicategory} that $\I$ can be
considered as a functor from $\Surf$ to $\Grp$ (the category of
 groups and homomorphisms), defined on objects by $\S\mapsto
\IS$ and on morphisms by $\iota\mapsto \iota_*$. Let $\mathcal{W}$ be the functor from $\Surf$ to $\Grp$ defined on objects by
$\mathcal{W}(\Sigma)= W_\S$ and on morphisms by $\mathcal{W}(\iota)=\iota_*$ as in Definition~\ref{def:WSiota}. Then we can rephrase Theorem~\ref{thm:commutes} as follows:
\begin{theorem}\label{thm:naturaltransformation}
  There is a natural transformation $\tau$ from the Torelli functor $\I$ to the functor $\mathcal{W}$ which assigns to each surface $\S$ the surjective homomorphism
  $\tau_\S\colon \IS\to W_\S$.
\end{theorem}

\section{Fundamental calculations and surjectivity of $\tau_\S$}
\label{sec:compute}
 In this section we calculate $\tau_\S$ on many simple but fundamental elements of $\IS$, including the natural ``point-pushing'' subgroups. Using this, we prove in
Section~\ref{sec:surjectivity} that $W_\S$ is exactly the image of
$\tau_\S$. These
results are also used in Section~\ref{sec:KSorbits}.

\subsection{Separating twists}\label{sec:st}
Let $\S$ be the surface $S_{0,n}$ with the totally separated
partition. The homology group $\HS$ is trivial. It follows that $W_\S=0$, so by Theorem~\ref{thm:WScontained} we have $\tau_\S(\phi)=0$
for any $\phi\in \IS$.

Any Dehn twist $T_\gamma$ is supported on the annulus which is the regular neighborhood of $\gamma$, and the induced partition on this annulus $S_{0,2}$ is  totally separated exactly when $\gamma$ is $\P$--separating.
Applying the naturality of $\tau_\S$, we obtain
the following corollary.
\begin{proposition}\label{prop:st}
  If $\S$ is a partitioned surface and $\gamma$ is a $\P$--separating curve, then $\tau_\S(T_\gamma)=0$.
\end{proposition} 
More generally, we have the following.
\begin{proposition}
 If $\phi\in \IS$ is supported on a totally separated
  genus 0 subsurface, then $\tau_\S(\phi)=0$.
\end{proposition}

\subsection{Bounding pair maps}\label{sec:bp}
Let $\S$ be a surface $S_{0,3}$ of genus 0 with 3 boundary components
$z_0$, $a_1$, and $a_2$, endowed with the partition
$\P=\{\{z_0\},\{a_1,a_2\}\}$ and basepoint $\ast\in z_0$.
 Let $\phi=T_{a_1}T_{a_2}^{-1}$. This is a bounding pair, and $\S$ is
the minimal connected surface on which $\phi$ is supported. The surface
$\Shat$ has genus 1 with 2 boundary components, $z_0$ and
$z_1$. Its fundamental group has rank 3, and we may choose a basis
$\{\alpha,\beta,\zeta\}$ for $\pi_1(\Surhat,\ast)$ so that the first two terms descend to a basis
$\{a,b\}$ for $\HS$, the generator $\alpha$ can be represented by a loop in $S$, and $\zeta$ descends to $z=z_1$ in $\NS$. For an appropriate choice of generators we have
$\phi(\alpha)=\alpha$ and $\phi(\beta)=\beta\zeta^{-1}$, so $\tau_\S(\phi)\in \Hom(\HS,\NS)$ is
defined by $a\mapsto 0$ and $b\mapsto -z$. In the alternate notation of Definition~\ref{def:alternate} for elements of $W_\S$, we
have $\tau_\S(\phi)=a\wedge z\in W_\S$.
The same argument
applies when the basepoint lies in $z_1$. Since every
bounding pair in a partitioned surface $\S$ sits inside at least one
such $S_{0,3}$ (for example, the regular neighborhood of the curves
together with an arc connecting them), we can apply the naturality of
$\tau_\S$ to obtain the following corollary.
\begin{proposition}\label{prop:BP}
  Given a bounding pair $T_\gamma T_{\delta}^{-1}$ defined by
  nonseparating curves $\gamma$ and $\delta$, let $\zeta$ be a
  separating curve that cobounds a pair of pants with
  $\gamma\cup\delta$. Orient these curves so that the pair of pants lies on the left
    side of $\zeta$ and the right side of $\gamma$ and $\delta$ (or
    vice versa). Let $a$ be the homology class of $\gamma$, and let
  $z$ be the class of $\zeta$ in $\NS$.  Then we have
  $\tau_\S(T_\gamma T_{\delta}^{-1})=a\wedge z$ in $W_\S$.
\end{proposition}
 Theorem~\ref{thm:commutes} guarantees that this recipe is well-defined, even though this is not obvious. We can check for example that changing the orientation of all three curves would negate both $a$ and
$z$, and thus would preserve $a\wedge z$. Similarly, there is a curve $\zeta'$ on
the other side of $\gamma\cup\delta$ whose class $z'\in \NS$ differs from $z$ by a term of the form
$a\wedge b$, in which case $a\wedge z'=a\wedge z$. But there are many different curves $\zeta$ bounding pairs of pairs with $\gamma$ and $\delta$; a key strength of the naturality of $\tau_\S$ is that it lets us choose any pair of pants that we like.

\subsection{Lantern core maps}\label{sec:lc}
We say that $\alpha$ and $\beta$ span a \emph{lantern core} if their geometric intersection number is 2 and their algebraic intersection number is 0.
A \emph{lantern core} map is the commutator $[T_\alpha,T_\beta]$ of Dehn twists around two such
curves.  Since $\alpha$ and
$\beta$ have algebraic intersection 0, $T_\alpha$ and $T_\beta$ act
on $\HS$ by commuting transvections. It follows that
$[T_\alpha,T_\beta]\in \IS$. Such maps were first used by Johnson in \cite{JohnsonII};
they are called simply intersecting pair maps in Putman
\cite{CuttingPasting}. 
\begin{figure}[h]
  \label{fig:arcs}
  \begin{center}
    \includegraphics[width=105mm]{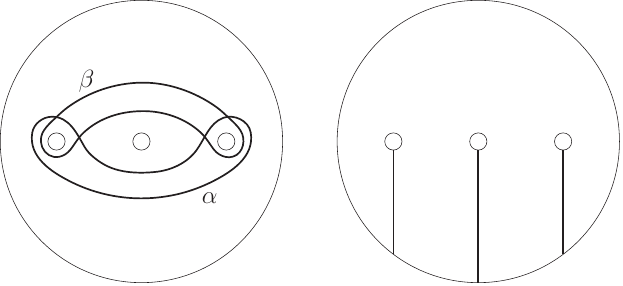}\\
  \end{center}
  \caption{The simple closed curves $\alpha$ and $\beta$, and the
    arcs $B_1,B_2,B_3$.}
\end{figure}

The regular neighborhood of two such curves is
always a lantern, so it suffices to compute $\tau_\S$ for partitioned surfaces $\S=(S_{0,4},\P,\ast)$. To begin, we say that $\alpha$ and $\beta$ span a \emph{nonseparating lantern core} if $\alpha\cup \beta$ is $\P$--nonseparating; that is, the induced partition on the regular neighborhood $S_{0,4}$ of $\alpha\cup\beta$ is the nonseparating partition $\P=\{P_0\}$.
Let $\alpha$ and $\beta$ be as in Figure~\reffigurearcs, and let $a_1,a_2,a_3$ be the homology classes in $\HS$ of the three boundary components in the center.

\begin{proposition}
\label{prop:lanterncore}
If $\alpha$ and $\beta$ span a nonseparating lantern core, then $\tau_\S([T_\alpha,T_\beta])=-a_1\wedge a_2\wedge a_3$.
\end{proposition}
\begin{proof}
If $\S=(S_{0,4},\{P_0\},\ast)$, the surface $\Surhat$ has genus 3 and 1 boundary component. A basis for
$\pi_1(\Surhat,\ast)$ is given by curves $\alpha_1,\alpha_2,\alpha_3$ traveling clockwise around the three central boundary components, together with curves $\beta_1,\beta_2,\beta_3$ whose
intersection with $\S$ are the arcs $B_1,B_2,B_3$ respectively, oriented
bottom-to-top. These descend to a basis
$\{a_1,a_2,a_3,b_1,b_2,b_3\}$ for $\HS$.

Let $\phi=[T_\alpha,T_\beta]$. Since $a_1$, $a_2$, and $a_3$ are
contained in $D_0$, from condition (III) of Definition~\ref{def:WS} we know that \[\tau_\S(\phi)(a_1)=\tau_\S(\phi)(a_2)=\tau_\S(\phi)(a_3)=0.\] The action of
$\phi$ on the arcs $B_i$ is displayed in Figure~\reffigurearcsmoved.
\begin{figure}
\label{fig:arcsmoved}
\begin{center}
    \includegraphics[width=162mm]{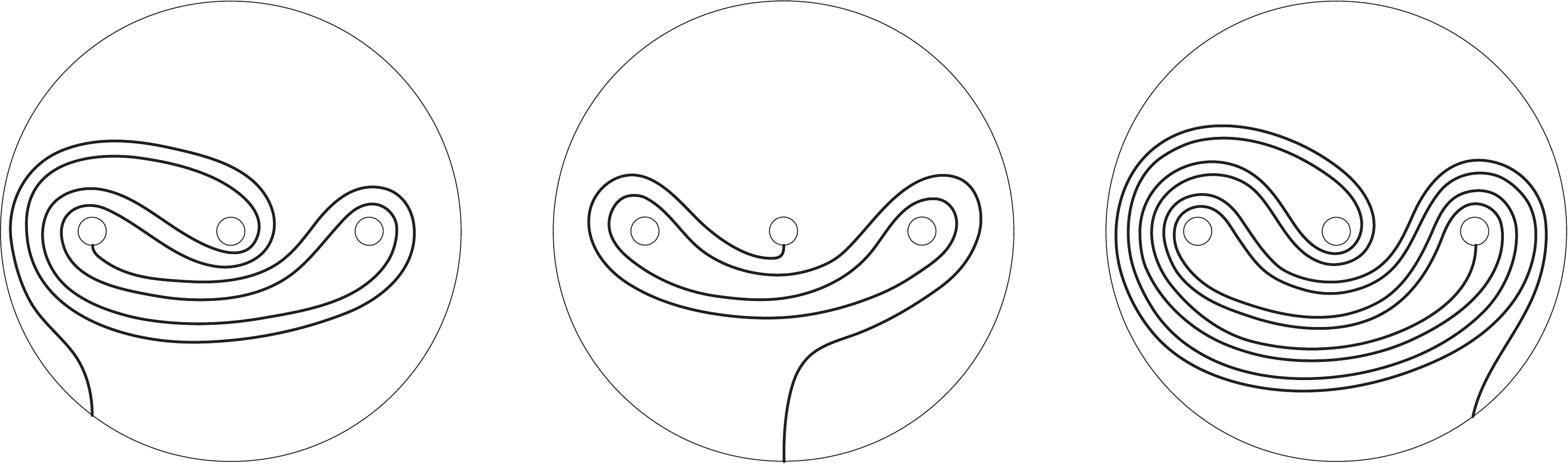}\end{center}
  \caption{The arcs $\phi(B_1), \phi(B_2)$ and $\phi(B_3)$.}
\end{figure}
 Thus we have:
\begin{equation}
\label{eq:lanternarcs}
\begin{array}{rl}
  B_1&\mapsto [\alpha_1\alpha_2\alpha_1^{-1},\alpha_3]B_1\\
  B_2&\mapsto [\alpha_3^{-1},\alpha_1^{-1}]B_2\\
  B_3&\mapsto
  [\alpha_3^{-1}\alpha_1^{-1}\alpha_3,\alpha_1\alpha_2^{-1}\alpha_1^{-1}]B_3
\end{array}\end{equation}
It follows that:
\begin{align*}
  \tau_\S(\phi)(b_1)&=a_2\wedge a_3\\
\tau_\S(\phi)(b_2)&=(-a_3)\wedge (-a_1)=a_3\wedge a_1\\
\tau_\S(\phi)(b_3)&=(-a_1)\wedge (-a_2)=a_1\wedge a_2
\end{align*} As an element of $W_\S=\bwedge^3 \HS$, this is the element
$\tau_\S(\phi)=-a_1\wedge a_2\wedge a_3$. The naturality of $\tau_\S$ implies that the same formula holds for a nonseparating lantern core in any surface $\S$, as claimed.
\end{proof}

From the same computations we can deduce that any other lantern core map is contained in the Johnson kernel $\KS$. Indeed for any other surface $\S=(S_{0,4},\P,\ast)$ for which $\P$ is not the nonseparating partition, the rank of $\DS^\perp$ is at most 2. Thus $\bwedge^3 \DS^\perp=0$, so the exact sequence \eqref{eq:WSses} implies that $\tau_\S(\phi)\in W_\S$ is determined by the $d_i(\phi)\in \HS$. But our computations of $\phi(B_i)$ in \eqref{eq:lanternarcs} show that $\phi(B_i)B_i^{-1}\in [\pi_1(\Surhat),\pi_1(\Surhat)]$ for each $i$, and so $d_i(\phi)=0\in \HS$ for all $i$.
\begin{corollary}
\label{prop:separatinglanterncore}
If $\alpha$ and $\beta$ span a  lantern core which is not nonseparating, then $\tau_\S([T_\alpha,T_\beta])=0$.
\end{corollary}

\subsection{Disk-pushing subgroups}\label{sec:pp}
In this section we determine the restriction of $\tau_\S$ to certain
``disk-pushing'' subgroups of $\I(\S)$. Let $\S\to \S'$ be a simple
capping obtained by attaching a disk to a separating boundary component $z=z_i$. In particular, we assume that $\S$ has at least two boundary components, and that $\ast$ is not contained in $z$. Such a capping induces a surjection
$\Mod(S)\to \Mod(S')$, whose kernel is isomorphic to
$\pi_1(UTS',v)$, where $v$ is a unit vector at the center of the disk
glued to $z$ (Johnson~\cite{JohnsonI}).

\begin{remark}
The conventions for composition in $\Mod(S)$ and in $\pi_1(UTS',v)$
unfortunately disagree; in the former we take composition of functions and in the
latter we take concatenation of paths. As a result, the isomorphism of
$\ker(\Mod(S)\to \Mod(S'))$ with $\pi_1(UTS',v)$ is defined as follows.
Given $\phi\in \ker(\Mod(S)\to \Mod(S'))$, extend $\phi$ by the identity to
$S'$; by definition, there is an isotopy $h^\phi_t$ of $S'$ from
$h^\phi_0=\phi$ to $h^\phi_1=\id$. Since $\phi$ fixes $v$, the image of $v$ under this isotopy
is a loop of tangent vectors from $h^\phi_0(v)=\phi(v)=v$ to $v$. We associate to $\phi$ the corresponding element
$\gamma_\phi\coloneq h^\phi_t(v)\in \pi_1(UTS',v)$. Note that given
another map $\psi$ fixing $v$, $h^\phi_t\circ \psi$ satisfies
$h^\phi_t\circ \psi(v)=h^\phi_t(v)$. Thus an isotopy from $\phi\circ
\psi$ to $\id$ can be obtained by concatenating the isotopy
$h^\phi_t\circ \psi$ from $\phi\circ \psi$ to $\psi$ with the isotopy
$h^\psi_t$ from $\psi$ to $\id$; the path this determines is the
concatenation $\gamma_\phi\cdot \gamma_\psi$. This verifies that
$\gamma_{\phi\circ \psi}=\gamma_\phi\cdot \gamma_\psi$. We remark that
this isomorphism is the opposite of the identification naively
suggested by the ``disk-pushing'' label.
\end{remark}

\begin{proposition}
\label{prop:seppushTor}
For a separating component, the entire disk-pushing subgroup $\pi_1(UTS',v)$ is contained in $\IS$.
\end{proposition}
\begin{proof}
Note that we have $H(\S')\simeq H(\S)$, since $z$ is separating and thus vanishes in homology. For any curve $\gamma$ in $S$ and any $\phi\in\pi_1(UTS',v)$, the curve $\phi(\gamma)$ is homotopic to $\gamma$ in $S'$, which implies that $[\phi(\gamma)]$ is homologous to $[\gamma]$ in $H(\S)$. This shows that $\phi\in\IS$ as desired.
\end{proof} 

\begin{proposition}
\label{prop:seppush}
  The restriction $\tau_\S\colon \pi_1(UTS',v)\to W_\S$ of the Johnson homomorphism to the disk-pushing subgroup
  determined by $z_i$ is the composition
  \[\ker\big(\IS\to \I(\S')\big)\iso \pi_1(UTS',d)
  \twoheadrightarrow \pi_1(S',d) \to H(\S')\iso
  \HS\hookrightarrow W_\S,\] where the last map is the embedding $\HS\hookrightarrow W_\S$ defined by
  $x\mapsto -x\wedge z_i$.
\end{proposition}

 Let $d\in S'$ be the projection of $v\in UTS'$, and for $\gamma\in\pi_1(UTS',v)$, let $\overline{\gamma}$ denote its projection to $\pi_1(S',d)$.
\begin{proof}
For any $\phi\in \ker(\IS\to \I(\S'))$, the naturality of $\tau_\S$
implies that $\tau_\S(\phi)\in \ker(W_\S\to W_{\S'})$. For a simple
capping $\S\to \S'$, by Definition~\ref{def:WSiota} the map $W_\S\to W_{\S'}$ is induced by the
quotient $\NS\twoheadrightarrow \NS/\langle z_i\rangle=N(\S')$. We thus
know \emph{a priori} that $\tau_\S(\phi)$ must be contained in the
subspace $\{x\wedge z_i|x\in \HS\}$. By Theorem~\ref{thm:deltas} we thus have \[\tau_\S(\phi)=d_i(\phi)\wedge z_i.\] Thus it suffices to prove that if $\phi\in\ker(\IS\to\I(\S'))$ corresponds to $\gamma_\phi\in\pi_1(UTS',v)$, we have \begin{equation}
\label{eq:PPdelta}
d_i(\phi)=-[\overline{\gamma_\phi}].
\end{equation}

Recall that $d_i(\phi)$ is the homology class of $\phi(A)A^{-1}$
  where $A=A_i$ is an arc from the basepoint to $z_i$. 
Since $\HS\iso H(\S')$,
  we may consider $A$ as an arc in $S'$ from the basepoint to $d$. If $\overline{\gamma_\phi}$ is simple and disjoint from $A$
  as a based loop, then $\phi(A)$ will be homotopic in $S'$ to the
  concatenation $A\cdot \overline{\gamma_\phi}^{-1}$. Note that if
  $\overline{\gamma_\phi}$ is simple, we can always choose $A$
  disjoint from it. For such elements, we have
  $\phi(A)A^{-1}=A\cdot\overline{\gamma_\phi}^{-1}\cdot A^{-1}$. This is
  homologous in $S'$ to
  $[\overline{\gamma_\phi}^{-1}]=-[\overline{\gamma_\phi}]$, and since $H(\S')\simeq \HS$ this implies that $[\phi(A)A^{-1}]=-[\overline{\gamma_\phi}]\in \HS$. Thus for such
  elements
  \eqref{eq:PPdelta} holds and $\tau_\S(\phi)=-[\overline{\gamma_\phi}]\wedge z$.

  The kernel $\ker(\pi_1(UTS',v)\to \pi_1(S',d))$ is generated by a
  twist around the boundary component $z_i$. By Proposition~\ref{prop:st} $\tau_\S$ vanishes for any separating twist, so \eqref{eq:PPdelta} holds for this element as well. The group
  $\pi_1(S',d)$ is normally generated by (in fact, generated by)
  elements represented by simple loops. It follows that
  $\pi_1(UTS',v)$ is normally generated by $\gamma$ for which
  $\overline{\gamma}$ is either simple or trivial. Since \eqref{eq:PPdelta} holds for generators of either form, it holds
  for all elements of $\ker(\IS\to \I(\S'))\iso \pi_1(UTS',v)$. This completes the proof of the proposition.
\end{proof}
\begin{remark}
\label{rem:rationalhomology}
  Note that for many surfaces $\S$, the subgroup $\pi_1(UTS',v)$ is generated by bounding pairs, so we could obtain another proof of the proposition 
  by applying Proposition~\ref{prop:BP}. This was the approach used by Johnson to prove the proposition in the 
  classical case when $S$ has a single boundary component. In fact the proposition is almost automatic in this case: the
  restriction of $\tau_S$ to the point-pushing subgroup
  $\pi_1(UTS',v)$ is a map to the torsion-free abelian group $\bwedge^3 H_1(S;\Z)$, so it
  factors through the surjection $\pi_1(UTS',v)\twoheadrightarrow
  H_1(S;\Z)$, and moreover extends to a $\Mod(S)$--equivariant map $H_1(S;\Q)\to \bwedge^3 H_1(S;\Q)$. Since $H_1(S;\Q)$ is an irreducible $\Sp_{2g}\Z$--module and $\bwedge^3 H_1(S;\Q)$ contains a unique submodule isomorphic to $H_1(S;\Q)$ (embedded by $x\mapsto x\wedge \omega$), Schur's lemma implies that the proposition holds up to a
  multiplicative constant. This constant can be detected by computing $\tau_S$ for a single element.
\end{remark}

\subsection{Surjectivity of $\tau_\S$}
\label{sec:surjectivity} We proved in Theorem~\ref{thm:WScontained} that $\im \tau_\S$ is contained in $W_\S$. To complete the proof of Theorem~\ref{thm:WS}, it remains to show that
$\IS$ surjects to $W_\S$ under $\tau_\S$.

\begin{proof}[Proof of Theorem~\ref{thm:WS}]
The surjection $W_\S\twoheadrightarrow (\DS^\perp)^k$ defined by $f\mapsto (\delta_1(f),\ldots,\delta_k(f))$ has kernel equal to $\bwedge^3 \DS^\perp$, as we saw in the exact sequence \eqref{eq:WSses}: \[0\to \bwedge^3 \DS^\perp\to W_\S\to (\DS^\perp)^k\to 0\] 
Fix a symplectic basis $\{a_i,b_i\}\cup\{a^i_j,b^i_j\}$ for $\HS$ so that $\{a^i_j\}$ provides a basis for $D_j$, each $a^i_j$ is represented by a boundary component, and $\{a_i,b_i\}\cup\{a^i_j\}$ provides a basis for $\DS^\perp$.

\para{The image $\im\tau_\S$ surjects to $(\DS^\perp)^k$} Fix $J\geq 1$ and let $\eta_J$ be the $\P$--separating curve cutting off a subsurface of genus 0 bounded by $\eta_J$ together with all the boundary components lying in $P_J$. We will show  for any $x\in \{a_i,b_i\}\cup\{a^i_j\}$ that we can find $\phi\in \IS$ with $\tau_\S(\phi)=x\wedge [\eta_J]$. Note that $\eta_J=[a^1_J,b^1_J]+\cdots+[a^{m_J}_J,b^{m_J}_J]+z_J$ in $\NS$, so $f=x\wedge [\eta_J]$ satisfies $\delta_J(f)=x$ and $\delta_j(f)=0$ for $j\neq J$. Since these $x$ form a basis for $\DS^\perp$, this will verify that  $\im\tau_\S$ surjects to $(\DS^\perp)^k$.

First, consider the case when $x=a^i_J$ for some $i$. Let $\gamma$ be a curve homotopic to $a^i_J$, and let $\delta$ be the band sum of $\gamma$ with $\eta_J$. (The band sum of two simple closed curves is their connected sum along some simple arc connecting the two curves.) The homologous curves $\gamma$ and $\delta$ determine a bounding pair $T_\gamma T_\delta^{-1}$ that cobounds a pair of pants with $\eta_J$, so Proposition~\ref{prop:BP} states that $\tau_\S(T_\gamma T_\delta^{-1})=a^i_J\wedge [\eta_j]$.

For the remaining cases, let $S_J$ be the component of $S-\eta_J$ containing the basepoint, and consider the disk-pushing subgroup of $\I(\S_J)$ determined by $\eta_J$. By Proposition~\ref{prop:seppush}, we can find $\phi\in \I(\S_J)$ with $\tau_{\S_J}(\phi)=x\wedge [\eta_J]$ for any $x\in H(\S_J)$ in the image of $\pi_1(S_J)$. But the image of $H_1(S_J)$ in $\HS$ is spanned by $\{a_i,b_i\}\cup \{a^i_j,b^i_j\,|\, j\neq J\}$. Applying Theorem~\ref{thm:commutes} to the non-collapsing inclusion of $\S_J$ into $\S$, we obtain a disk-pushing map $\phi\in\IS$ with $\tau_\S(\phi)=x\wedge [\eta_J]$ for the remaining basis elements $x$, as desired.

\para{The image $\im\tau_\S$ contains $\bwedge^3\DS^\perp$} Consider the natural basis of $\bwedge^3 \DS^\perp$ induced by the basis $\{a_i,b_i\}\cup\{a^i_j\}$ of $\DS^\perp$. First, consider the basis elements of the form $x\wedge a_I\wedge b_I$ for some $I$ and some other basis element $x$. Realize $a_I$ and $b_I$ by simple closed curves intersecting once, and let $\eta_I$ be the regular neighborhood of their union. Let $S_I$ be the component of $S-\eta_I$ containing the basepoint, and consider the disk-pushing subgroup of $\I(\S_I)$ determined by $\eta_I$. Any $x$ in our basis distinct from $a_I$ and $b_I$ can be realized by a loop in $S_I$, so by Proposition~\ref{prop:seppush} we can find a disk-pushing map $\phi\in \I(\S_I)$ with $\tau_{\S_I}(\phi)=x\wedge [\eta_I]$. Since $[\eta_I]=a_I\wedge b_I$ in $\NS$, applying Theorem~\ref{thm:commutes} to the non-collapsing inclusion of $\S_I$ into $\S$ implies that $\tau_\S(\phi)=x\wedge a_I\wedge b_I$.

Any other basis element is of the form $x\wedge y\wedge z$ where $(x,y)=(y,z)=(z,x)=0$. We may thus realize these homology classes by disjoint curves $\alpha,\beta,\gamma$. Let $\delta$ be the band sum of all three curves (to be precise we should specify orientations, but this will only change the answer by a sign). We obtain a lantern $S_{0,4}$ bounded by $\alpha$, $\beta$, $\gamma$, and $\delta$. If $\phi$ is a lantern core map supported on this lantern, then by Proposition~\ref{prop:lanterncore}, $\tau_\S(\phi)=\pm x\wedge y\wedge z$. This verifies that a basis for $\bwedge^3 \DS^\perp$ is contained in $\im \tau_\S$. Together with Theorem~\ref{thm:WScontained}, this completes the proof of Theorem~\ref{thm:WS} and shows that the image of $\tau_\S$ is exactly $W_\S$.\end{proof}

\section{Orbits of curves under $\KS$}
\label{sec:KSorbits}
In this section, we combine the characterization of $\im \tau_\S$ from Section~\ref{sec:Jh} with the computations from Section~\ref{sec:compute} to describe the orbits of simple closed curves under the Johnson kernel $\K(S)$.

\subsection{Orbits of nonseparating twists}\label{sec:Kononsep}
Let $S=S_{g,1}$ and $\S=(S,\{P_0\},\ast)$, so that $\HS=H(S)=H_1(S)$. In this section we consider nonseparating curves in $S$; we always consider curves to be oriented. It was known to Johnson~\cite{JohnsonConjugacy} that two nonseparating curves are in the same
$\I(S)$--orbit if and only if they are homologous. We thus
consider two homologous nonseparating curves $C$ and $D$ with
homology class $a\in H(S)$, and ask when they are in the same
$\K(S)$--orbit.
 We first recall the theorem.

\begin{theorem:nonsep}[$\KSs$--orbits of nonseparating curves]
Let $C$ and $D$ be nonseparating curves homologous to $a\in H_1(S)$. The following are equivalent:
\begin{enumerate}[1.]
\item The nonseparating curves $C$ and $D$ are equivalent under $\KSs$.
\item $T_CT_D^{-1}\in\KSs$.
\item For some representatives $\gamma,\,\delta\in\pi_1(S,\ast)$ of the curves $C$ and $D$, the class $[\gamma\delta^{-1}]\in \Gamma_2(S)/\Gamma_3(S)\iso \bwedge^2 H_1(S)$ lies in the subspace $a\wedge H_1(S)$.
\item For any representatives $\gamma,\,\delta$ of $C$ and $D$,  $[\gamma\delta^{-1}]\in \Gamma_2(S)/\Gamma_3(S)\iso \bwedge^2 H_1(S)$ lies in $a\wedge H_1(S)$.
\end{enumerate}
\end{theorem:nonsep}
The hardest step will be to prove that (3 $\implies$ 1), and this is where the  results of this paper are required. To simplify the proof, we introduce two  equivalent restatements of Condition 2. Fix $b\in H(S)$ satisfying $\omega(a,b)=1$. This induces a direct splitting $H(S)=\langle b\rangle \oplus a^\perp$, which induces the decomposition (the second summand is embedded by $b\otimes x\wedge y\mapsto b\wedge x\wedge y$):\begin{equation}
\label{eq:WSsplit}
W_S=\bwedge^3 H(S)=\bwedge^3 a^\perp\ \oplus\ \langle b\rangle \otimes\bwedge^2 a^\perp
\end{equation} 
Then the four conditions above are also equivalent to the following two conditions.
\begin{enumerate}[\emph{1.}]
\item[\emph{5.}] \emph{For any $\phi\in \I(S)$ with $\phi(C)=D$, the element $\tau(\phi)\in W_S$ lies in the subspace $\bwedge^3 a^\perp\ \oplus\ \langle b\rangle \otimes a\wedge a^\perp$.}
\item[\emph{6.}] \emph{For any $\phi\in \I(S)$ with $\phi(C)=D$ we have $\tau(\phi)(a)\in a\wedge a^\perp$.}
\end{enumerate}
By Johnson~\cite{JohnsonConjugacy}, there exists some $\phi\in \I(S)$ with $\phi(C)=D$, so these conditions are not vacuous.
\begin{proof}[Proof of Theorem~\ref{thm:intrononsep}]
The Dehn twist $T_C$ acts on $H(S)$ by the transvection $t_a\in \Sp(H(S))$ defined by $t_a(x)=x+\omega(a,x)a$. Since $T_D$ also acts on $H(S)$ by $t_a$, we always have $T_CT_D^{-1}\in \I(S)$. It was known to Johnson~\cite{JohnsonConjugacy} that for any homologous nonseparating curves $C$ and $D$ there exists $\phi\in \I(S)$ satisfying $\phi(C)=D$.
Given such a $\phi\in \I(S)$, we can write $T_C T_D^{-1}=T_C T_{\phi(C)}^{-1}= T_C \phi T_C^{-1} \phi^{-1}$. Since $\tau$ is $\Mod(S)$--equivariant (Lemma~\ref{lem:equivariance}), we have $\tau(T_C\phi T_C^{-1})=t_a\cdot \tau(\phi)$ and thus 
\begin{equation}
\label{eq:tauCD}
\tau(T_C T_D^{-1})=\tau(T_C \phi T_C^{-1})-\tau(\phi)=(t_a-\id)\cdot \tau(\phi).
\end{equation}

\para{(1 $\implies$ 2)} If $\phi(C)=D$ for some $\phi\in \K(S)$, the  equation \eqref{eq:tauCD} becomes
\[\tau(T_CT_D^{-1})=(t_a-\id)\cdot \tau(\phi)=(t_a-\id)\cdot 0=0,\]
and thus $T_CT_D^{-1}\in \KSs$.

\para{(2 $\iff$ 5)} 
Given $\phi\in \IS$ so that $\phi(C)=D$, \eqref{eq:tauCD} implies that \[\tau(T_CT_D^{-1})=0\qquad \iff\qquad \tau(\phi)\in \ker\big(t_a-\id\colon W_S\to W_S\big)\] so we must identify $\ker(t_a-\id)$. The transvection $t_a$ acts as the identity on $a^\perp$, so $(t_a-\id)$ acts by 0 on the first summand $\bwedge^3 a^\perp$ of \eqref{eq:WSsplit}. On the second summand, since $t_a(b)=b+\omega(a,b)a=b+a$ we have $(t_a-\id)(b\otimes x\wedge y)=a\wedge x\wedge y$. Since $H(S)$ is torsion-free, the kernel of the map $\bwedge^2 a^\perp\to \bwedge^3 a^\perp$ given by $x\wedge y\mapsto a\wedge x\wedge y$ is precisely $a\wedge a^\perp$. We conclude that $\ker(t_a-\id)=\bwedge^3 a^\perp\, \oplus\, \langle b\rangle\otimes a\wedge a^\perp$, as desired.

\para{(5 $\iff$ 6)}
Under the isomorphism $\Hom(H(S),\bwedge^2 H(S))\iso H(S)\otimes \bwedge^2 H(S)$, the projection of $f\in \Hom(H(S),\bwedge^2 H(S))$ onto the summand $\langle b\rangle \otimes \bwedge^2 H(S)$ of the latter group is precisely $b\otimes f(a)$. Moreover, this projection vanishes on the summand $\bwedge^3 a^\perp$ of \eqref{eq:WSsplit}, and restricts to an isomorphism from the summand $\langle b\rangle \otimes \bwedge^2 a^\perp$ of \eqref{eq:WSsplit} to $\langle b\rangle \otimes \bwedge^2 a^\perp < \langle b\rangle \otimes \bwedge^2 H(S)$.  It follows that any $f\in W_S$ satisfies $f(a)\in \bwedge^2 a^\perp$, and moreover that $f$ lies in $\bwedge^3 a^\perp\,\oplus\, \langle b\rangle\otimes a\wedge a^\perp$ if and only if $f(a)\in a\wedge a^\perp$. 
  
\para{(4 $\implies$ 6 $\implies$ 3)} 
Choose a representative $\gamma\in\pi_1(S,\ast)$ of $C$.
Given $\phi\in \I(S)$ satisfying $\phi(C)=D$, let $\delta=\phi(\gamma)$ be the resulting representative of $D$ . Since $\gamma$ represents $a\in H(S)$ we may compute $\tau(\phi)(a)$ by
\[\tau(\phi)(a)=[\phi(\gamma)\gamma^{-1}]=[\delta\gamma^{-1}]=-[\gamma\delta^{-1}]\in N(S)\iso \bwedge^2 H(S).\] 
Condition 6 states that $\tau(\phi)(a)\in a\wedge a^\perp$, so we have $[\gamma\delta^{-1}]\in a\wedge a^\perp \subset a\wedge H(S)$, verifying Condition 3. Conversely, Condition 4 states that any representatives $\gamma,\delta$ satisfy $[\gamma\delta^{-1}]\in a \wedge H(S)$, so $\tau(\phi)(a)\in a\wedge H(S)$. We noted in the previous paragraph that $\tau(\phi)(a)$ always lies in $\bwedge^2 a^\perp$, so we conclude that $\tau(\phi)(a)\in a\wedge a^\perp$, verifying Condition 6.

\para{(3 $\implies$ 4)} Fix representatives $\gamma$ and $\delta$ of the curves $C$ and $D$. Any other  representative $\delta'$ of the same curve $D$ satisfies
  $\delta'=\xi\delta\xi^{-1}$ for some $\xi\in \pi_1(S)$. The difference between $\gamma\delta^{-1}$ and
  $\gamma\delta'^{-1}$ is thus equal modulo $\Gamma_3(S)$ to
  $[\delta',\xi]$, which corresponds to $a\wedge [\xi]$ under the
  identification $\Gamma_2(S)/\Gamma_3(S)\iso \bwedge^2 H(S)$. The same argument
  applies to $\gamma$, and so we conclude that the class $[\gamma\delta^{-1}]$ is well-defined modulo $a\wedge H_1(S)$. In particular, if $[\gamma\delta^{-1}]\in a\wedge H_1(S)$, then any other representatives $\gamma'$ and $\delta'$ satisfy $[\gamma'\delta'^{-1}]\in a\wedge H_1(S)$ as well.

\para{(5 $\implies$ 1)} Given two homologous nonseparating curves $C$ and $D$, choose $\phi\in \I(S)$ satisfying $\phi(C)=D$. Condition 4 tells us that $\tau(\phi)$ lies in $\bwedge^3 a^\perp\oplus \langle b\rangle\otimes a\wedge a^\perp$. Our goal will be to show that there exists $\psi\in \Stab_{\I(S)} C$ satisfying $\tau(\psi)=\tau(\phi)$, since then $\phi\psi^{-1}$ satisfies $\phi\psi^{-1}(C)=D$ and $\tau(\phi\psi^{-1})=0$, and this demonstrates that $C$ and $D$ are equivalent under $\K(S)$.

We thus need to show that $\tau_S(\Stab_{\I(S)}(C))$ is all of $\bwedge^3 a^\perp\oplus \langle b\rangle \otimes a\wedge a^\perp$. We do this by considering  $S-C$ as a partitioned subsurface of $S$. Formally, let $\tS=(S_{g-1,3},\tPP,\ast)$ where $\tPP$ is of the form $\{\{z_0\},\{a_1,a_2\}\}$ and $\ast\in z_0$. There is a natural inclusion $\iota\colon\tS\to S$ as the complement of a regular neighborhood of $C$. Any mapping class stabilizing $C$ lifts (non-uniquely) to $\tS$. Conversely, extension by the identity gives a surjection $\iota_*\colon\Mod(\tS)\to \Stab_{\Mod(S)}(C)$. By Paris--Rolfsen~\cite[Theorem 4.1(iii)]{ParisRolfsen}, the kernel of this surjection is cyclic, generated by $T_{a_1}T_{a_2}^{-1}$. Note that $T_{a_1}T_{a_2}^{-1}\in \I(\S')$, so $\I(\tS)$ surjects to $\Stab_{\IS(S)}(C)$. Let $z_1$ be the boundary component of $\widehat{\tS}$ corresponding to $\tP_1=\{a_1,a_2\}$. By Theorem~\ref{thm:WS}, noting that $D(\tS)=\langle a\rangle$, the short exact sequence \eqref{eq:WSses} becomes:
\begin{equation}
\label{eq:WSsub}
0\to \bwedge^3 a^\perp\to W_{\tS}\to a^\perp\to 0
\end{equation}

By the naturality proved in Theorem~\ref{thm:commutes}, $\tau_S(\Stab_{\IS}(C))$ is the image of the map $\iota_*\colon W_{\tS}\to W_S$ that sends $\delta_1(f)\wedge z_1\mapsto 0$ and is the identity on other factors.
But the kernel of this map restricted to $W_{\tS}$ is quite small, since condition (II) of Definition~\ref{def:WS} implies $\iota_*(f)(a)=\delta_1(f)\wedge a$. In particular, this means $\iota_*(f)$ cannot be 0 unless $\delta_1(f)=na$ for some $n\in \Z$. And indeed  $\tau_{\tS}(T_{a_1}T_{a_2}^{-1})=a\wedge z_1$ by Proposition~\ref{prop:BP}, and this element certainly lies in the kernel of $\iota_*$, so the kernel of $\iota_*$ is the cyclic subgroup spanned by $a\wedge z_1$. Under this projection the first factor $\bwedge^3 a^\perp<W_{\tS}$ of \eqref{eq:WSsub} is mapped isomorphically to $\bwedge^3 a^\perp<W_S$ of \eqref{eq:WSsplit}. The other factor $a^\perp$ of \eqref{eq:WSsub} is represented by $\delta_1(f)$, and by Theorem~\ref{thm:commutes} this factor is mapped by $\delta_1(f) \mapsto b\otimes \delta_1(f)\wedge a$, and this is a surjection onto $\langle b\rangle \otimes a\wedge a^\perp$ inside the factor $\langle b\rangle \otimes \wedge^2 a^\perp$ of \eqref{eq:WSsplit}.
\end{proof}

\subsection{Orbits of separating twists}\label{sec:Kosep}
As before we take  $S=S_{g,1}$ and $\S=(S,\{P_0\},\ast)$, and now consider separating curves in $S$. A separating curve $C$ has a canonical orientation, by taking the boundary $\partial S$ to be on the right side of $C$. The curve $C$ separates $S$ into two components $S_1$ and $S_2$, where $S_2$ contains $\partial S$, and we define $V_C\leq H(S)$ to be the subspace of $H(S)$ spanned by $H_1(S_1)$; this induces an orthogonal splitting $H(S)=V_C\oplus V_C^\perp$. Johnson~\cite[Theorem 1A]{JohnsonConjugacy} proved that two separating curves $C$ and $D$ are in the same
$\I(S)$--orbit if and only if the subspaces $V_{C}$ and $V_{D}$ coincide. We thus
consider two separating curves $C$ and $D$ with
$V_C=V_D=V$, and ask when they are in the same
$\K(S)$--orbit.

Recall from Proposition~\ref{prop:Tbracket} that the map $\beta\colon H(S)\otimes \bwedge^2 H(S)\to \Gamma_3(S)/\Gamma_4(S)\iso \L_3(S)$ induced by the commutator bracket has kernel $\bwedge^3 H(S)$. Let $H_1(S)\otimes \omega_V$ denote the subspace of $\L_3(S)$ which is the image under $\beta$ of the subspace $\{x\otimes \omega_V\,|\,x\in H(S)\}$ (we will see below that $\beta$ is injective on this subspace). Here $\omega_V\in \bwedge^2 H(S)$ represents the restriction of the symplectic form $\omega$ to the symplectic subspace $V$.

\begin{theorem:sep}[$\KSs$--orbits of separating curves]
Let $C$ and $D$ be separating curves cutting off the same symplectic subspace $V< H_1(S)$. The following are equivalent:
\begin{enumerate}[1.]
\item The separating curves $C$ and $D$ are equivalent under $\KSs$.
\item The separating twists $T_C$ and $T_D$ are conjugate in $\KSs$.
\item For some representatives $\gamma,\,\delta\in\pi_1(S,\ast)$ of the curves $C$ and $D$, the class $[\gamma\delta^{-1}]\in \Gamma_3(S)/\Gamma_4(S)$ lies in the subspace $H_1(S)\otimes \omega_V$.
\item For any representatives $\gamma,\,\delta$ of $C$ and $D$, the class $[\gamma\delta^{-1}]\in \Gamma_3(S)/\Gamma_4(S)$ lies in $H_1(S)\otimes \omega_V$.
\end{enumerate}
\end{theorem:sep}
As before, it will be useful to establish an additional equivalent condition. The splitting $H(S)=V\oplus V^\perp$ induces a decomposition
\begin{equation}
\label{eq:WSsplitsep}
W_S=\bwedge^3 H(S)=\bwedge^3 V^{\perp}\ \oplus\ (V\otimes \bwedge^2 V^\perp)\ \oplus\ (V^\perp \otimes \bwedge^2 V)\ \oplus\ \bwedge^3 V.
\end{equation}
Let $V^\perp\wedge \omega_V$ denote the subspace of $W_S$ spanned by $w\wedge \omega_V$ for $w\in V^\perp$. In the decomposition \eqref{eq:WSsplitsep}, $V^\perp\wedge \omega_V$ is contained in the third factor.
\begin{enumerate}[\emph{1.}]
\item[\emph{5.}] \emph{For any $\phi\in \I(S)$ with $\phi(C)=D$, the element $\tau(\phi)\in W_S\iso\bwedge^3 H(S)$ lies in the subspace $\bwedge^3 V\oplus \bwedge^3 V^\perp\oplus(V^\perp\wedge \omega_V)$.}
\end{enumerate}

In contrast with the nonseparating case, any representative $\gamma$ of $C$ is trivial in homology, so there is no reason \emph{a priori} that $\phi(\gamma)\gamma^{-1}$ should be related to $\tau(\phi)$. The following lemma, which is fundamental to the proof of Theorem~\ref{thm:introseparating}, shows that in fact $\phi(\gamma)\gamma^{-1}$ captures a significant portion of $\tau(\phi)$.
Since $C$ is separating, any representative $\gamma$ lies in $\Gamma_2(S)$, and Corollary~\ref{cor:TorLTtrivial} implies that $[\phi(\gamma)]=[\gamma]\in \L_2(S)$ for any $\phi\in \I(S)$. Thus $\phi(\gamma)\gamma^{-1}\in \Gamma_3(S)$ and we may consider its class $[\phi(\gamma)\gamma^{-1}]\in \L_3(S)$.
\begin{lemma}\label{lem:composition}
The homomorphism $\I(S)\to \L_3(S)$ defined by $\phi\mapsto [\phi(\gamma)\gamma^{-1}]\in \L_3(S)$
 is equal to the composition
\begin{equation}\label{eq:composition}\theta\colon\I(S)\overset{\tau}{\to}
\Hom(H(S),\bwedge^2 H(S))\overset{f|_V}{\to} \Hom(V,\bwedge^2 H(S))
\overset{\iso}{\to}V\otimes \bwedge^2 H(S)\overset{\beta}{\to}
\L_3(S)\end{equation}
Here $f\mapsto f|_V$ simply restricts the homomorphism $f$ to the subspace $V<H(S)$, and as before $\beta\colon \L_1(S)\otimes \L_2(S)\to \L_3(S)$ is the Lie bracket.
\end{lemma}
\begin{proof}
Choose a basis $\{\alpha_1,\beta_1,\ldots,\alpha_g,\beta_g\}$ for $\pi_1(S_1)$ with homology classes $\{a_i,b_i\}$ so that $\{a_1,b_1,\ldots,a_k,b_k\}$ form a symplectic basis for $V$ and 
  $\gamma=[\alpha_1,\beta_1]\cdots[\alpha_k,\beta_k]$. (For future reference, this implies that $[\gamma]=\omega_V\in \L_2(S)\iso \bwedge^2 H(S)$.)
  We  compute $\theta(\phi)=[\phi(\gamma)\gamma^{-1}]$ just as in the proof of Theorem~\ref{thm:WScontained}. Set $\eta_\phi(x)=x^{-1}\phi(x)$ and note that
  $\eta_\phi(x)\in\Gamma_2(S)$ for all $x\in \pi_1(S)$. We compute:
\begin{align*}
  \phi(\gamma)= \phi\left(\prod_{i=1}^k [\alpha_i,\beta_i]\right)
  &=\prod_{i=1}^k[\phi(\alpha_i),\phi(\beta_i)]\\
  &=\prod_{i=1}^k[\alpha_i\eta_\phi(\alpha_i),\beta_i\eta_\phi(\beta_i)]\\
  &\equiv \prod_{i=1}^k[\alpha_i,\beta_i] [\alpha_i,\eta_\phi(\beta_i)]
    [\eta_\phi(\alpha_i),\beta_i] \bmod{\Gamma_4(S)}\\
  &\equiv
  \left(\prod_{i=1}^k[\alpha_i,\eta_\phi(\beta_i)]
   [\eta_\phi(\alpha_i),\beta_i]\right)\cdot\gamma\bmod{\Gamma_4(S)}
\end{align*} Thus $\theta(\phi)=[\phi(\gamma)\gamma^{-1}]\in \Gamma_3(S)/\Gamma_4(S)$ is
represented by
$\prod_{i=1}^k[\alpha_i,\eta_\phi(\beta_i)][\eta_\phi(\alpha_i),\beta_i]$.
Since $\eta_\phi(\alpha_i)$ and $\eta_\phi(\beta_i)\in \Gamma_2(S)$ represent
$\tau(\phi)(a_i)$ and $\tau(\phi)(b_i)\in \bwedge^2 H(S)$, we see that\begin{equation}\label{eq:sumik}\theta(\phi)=\beta\left(\sum_{i=1}^k
a_i\otimes \tau(\phi)(b_i)-b_i\otimes \tau(\phi)(a_i)\right)\in \L_3(S).
\end{equation} The element $\sum_{i=1}^k
a_i\otimes \tau_S(\phi)(b_i)-b_i\otimes \tau_S(\phi)(a_i)$ represents the homomorphism $H(S)\mapsto \bwedge^2 H(S)$ defined by \[a_i\mapsto \begin{cases}\tau(\phi)(a_i)&1\leq i\leq k\\0&k<i\leq g\end{cases}\qquad b_i\mapsto \begin{cases}\tau(\phi)(b_i)&1\leq i\leq k\\0&k<i\leq g\end{cases}\] and this is precisely the restriction of $\tau(\phi)$ to $V$. This shows that that $[\phi(\gamma)\gamma^{-1}]\in \L_3(S)$ agrees with the composition $\theta(\phi)$ of \eqref{eq:composition}, and completes the proof of the lemma.
\end{proof}

\begin{proof}[Proof of Theorem~\ref{thm:introseparating}]
\para{(1 $\iff$ 2)} Since $C$ and $D$ are separating, we have $T_C\in \KSs$ and $T_D\in \KSs$ by Proposition~\ref{prop:st}. For any $\phi\in \Mod(S)$ satisfying $\phi(C)=D$ we have $\phi T_C \phi^{-1}=T_{\phi(C)}=T_D$. Thus there exists $\phi\in \KSs$ satisfying $\phi(C)=D$ if and only if $T_C$ and $T_D$ are conjugate in $\KSs$.

\para{(1 $\implies$ 3)}
If Condition 1 holds, there exists $\phi\in \K(S)$ satisfying $\phi(C)=D$. Choose a representative $\gamma$ of $C$, and let $\delta=\phi(\gamma)$ be the resulting representative of $D$. Lemma~\ref{lem:composition} shows that $[\gamma\delta^{-1}]\in \L_3(S)$ can be computed as $\theta(\phi)$. Since $\theta$ factors through the Johnson homomorphism $\tau$ and $\tau(\phi)=0$, we conclude that  $[\gamma\delta^{-1}]=0\in \L_3(S)$. This shows that $[\gamma\delta^{-1}]\in H_1(S)\otimes \omega_V$ for this choice of representatives $\gamma$ and $\delta$, verifying Condition 3.

\para{(3 $\implies$ 4)}
Fix representatives $\gamma$ and $\delta$ of the curves $C$ and $D$. We saw in the proof of Lemma~\ref{lem:composition} that $[\gamma]=[\delta]=\omega_V\in \L_2(S)\iso \bwedge^2 H(S)$. Any other representative $\delta'$ satisfies
$\delta'=\xi\delta\xi^{-1}$  for some $\xi\in \pi_1(S)$. 
The difference between $\gamma\delta^{-1}$ and
  $\gamma\delta'^{-1}$ is thus equal modulo $\Gamma_4(S)$ to
  $[\delta',\xi]$, which represents $-\beta([\xi]\otimes \omega_V)$ in $\L_3(S)$. The same argument
  applies to $\gamma$, and so we conclude that the class $[\gamma\delta^{-1}]$ is well-defined modulo $H_1(S)\otimes \omega_V$. Thus if $[\gamma\delta^{-1}]\in H_1(S)\otimes \omega_V$, then any other representatives $\gamma'$ and $\delta'$ satisfy $[\gamma'\delta'^{-1}]\in H_1(S)\otimes \omega_V$ as well.

\para{(4 $\implies$ 5)} Let $\mu$ be the composition \[\mu\colon \Hom(H(S),\bwedge^2 H(S))\overset{\iso}{\to}H(S)\otimes \bwedge^2 H(S)\overset{\pi_V\otimes \id}{\twoheadrightarrow} V\otimes \bwedge^2 H(S)\overset{\beta}{\to}\L_3(S),\] where $\pi_V$ is the orthogonal projection $H(S)\twoheadrightarrow V$. The dualization $H(S)^*\to V^*$ of $\pi_V$ simply restricts a functional on $H(S)$ to the subspace $V$:
\[f\mapsto f|_V\colon H(S)^* \iso H(S)\overset{\pi_V}{\twoheadrightarrow}V\iso V^*\] It follows that $\theta=\mu\circ \tau$.

For any $\phi\in \I(S)$ with $\phi(C)=D$, choose a representative $\gamma$ of $C$.  Lemma~\ref{lem:composition} gives that $\theta(\phi)=[\phi(\gamma)\gamma^{-1}]$, but Condition 4 states that $[\phi(\gamma)\gamma^{-1}]\in H_1(S)\otimes \omega_V$, so we must have $\theta(\phi)=\mu\circ \tau(\phi)\in H_1(S)\otimes \omega_V$.
To prove that Condition 4 implies Condition 5, it therefore suffices to show that when restricted to $W_S=\bwedge^3 H(S)<\Hom(H(S),\bwedge^2 H(S))$, the preimage $\mu^{-1}(H_1(S)\otimes \omega_V)$ under $\mu$ of the subspace $H_1(S)\otimes \omega_V<\L_3(S)$ is precisely $\bwedge^3 V\oplus \bwedge^3 V^\perp\oplus V^\perp\wedge \omega_V$.

We consider each factor of the decomposition \eqref{eq:WSsplitsep} in term. 
 The embedding of $\bwedge^3 H(S)$ into $H(S)\otimes \bwedge^2 H(S)$ sends $x\wedge y\wedge z$ to $x\otimes y\wedge z+y\otimes z\wedge x+z\otimes x\wedge y$. In particular, if $w,w',w''$ lie in $V^\perp$ and thus vanish under $\pi_V$, the element $w\wedge w'\wedge w''$ is annihilated by $\pi_V\otimes \id$. This shows that  $\bwedge^3 V^\perp \leq \ker \pi_V\otimes \id \leq \ker \mu$. Conversely, any $v\in V$ has $\pi_V(v)=v$, so if $v,v',v''\in V$ we have $\pi_V\otimes \id(v\wedge v'\wedge v'')=v\wedge v'\wedge v''$. However, by Proposition~\ref{prop:Tbracket} we have $\ker \beta=\bwedge^3 H(S)$. Thus $\beta(v\wedge v'\wedge v'')=0$, showing that $\bwedge^3 V\leq \ker \mu$ as well.
 
If $v\in V$ while $w,w'\in V^\perp$, two terms of $v\wedge w\wedge w'$ are annihilated by $\pi_V\otimes \id$, leaving $\pi_V\otimes \id(v\wedge w\wedge w')=v\otimes w\wedge w'$. Thus $\mu(v\wedge w\wedge w')=\beta(v\otimes w\wedge w')$. Similarly, if $w\in V^\perp$ and $v,v'\in V$, we have $\pi_V\otimes \id(w\wedge v\wedge v')=v\otimes v'\wedge w+v'\otimes w\wedge v$. We therefore have $\mu(w\wedge v\wedge v')=[v,[v',w]]+[v',[w,v]]$. By the Jacobi identity, this is equal to $-[w,[v,v']]=-\beta(w\otimes v\wedge v')$. These two factors are thus embedded by applying $\beta$ to $V\otimes \bwedge^2 V^\perp$ and $V^\perp\otimes \bwedge^2 V$. Since $\ker \beta=\bwedge^3 H(S)$ is disjoint from these subspaces, we conclude that $\mu^{-1}(H_1(S)\otimes \omega_V)$ intersects these factor only in $V^\perp\wedge \omega_V$.

\para{(5 $\implies$ 1)} As before, our goal is to prove that $\tau_S(\Stab_{\I(S)}(C))$ is equal to the subspace $\bwedge^3 V\oplus \bwedge^3 V^\perp\oplus V^\perp \wedge \omega_V$ of Condition 5. Given this, we can choose any $\phi\in \I(S)$ with $\phi(C)=D$. Condition 5 states that $\tau_S(\phi)\in \bwedge^3 V\oplus \bwedge^3 V^\perp\oplus V^\perp \wedge \omega_V$, so we can find $\psi\in \Stab_{\I(S)}(C)$ such that $\tau_S(\psi)=\tau_S(\phi)$. Then $\phi\psi^{-1}$ lies in $\K(S)$ and satisfies $\phi\psi^{-1}(C)=D$, demonstrating that $C$ and $D$ are equivalent under $\K(S)$.

Let $\S_1$ and $\S_2$ be the components of $S-C$, where $\S_2$ contains the basepoint. The partitioned surface $\S_1$ has underlying surface $S_{k,1}$, while the partitioned surface $\S_2$ has underlying surface $S_{g-k,2}$ and the totally separated partition. The natural inclusions  $\iota_1\colon \S_1\to S$ and $\iota_2\colon \S_2\to S$ induce natural identifications $(\iota_1)_*\colon H(\S_1)\iso V$ and $(\iota_2)_* H(\S_2)\iso V^\perp$. The stabilizer $\Stab_{\Mod(S)}(C)$ decomposes as \[\Stab_{\Mod(S)}(C)\simeq \Mod(S_1)\times_{\langle T_C\rangle} \Mod(S_2).\] The action of $\Mod(S_1)$ on $H(S)$ is the identity on $V^\perp$, and that of $\Mod(S_2)$ is the identity on $V$. Since $T_C\in \I(S)$, we obtain a decomposition \[\Stab_{\I(S)}(C)\simeq \I(\S_1)\times_{\langle T_C\rangle} \I(\S_2).\] Since $C$ is separating, we have $\tau_S(T_C)=0$, so $\tau_S(\Stab_{\I(S)}(C))$ is generated by $\tau_S(\I(\S_1))$ and $\tau_S(\I(\S_2))$.  Theorem~\ref{thm:WS} implies that that $\tau_{\S_1}(\I(\S_1))=W_{\S_1}\simeq \bwedge^3 H(\S_1)$, and Definition~\ref{def:WSiota} states that $(\iota_1)_*\colon W_{\S_1}\to W_S$ is just the inclusion of $\bwedge^3 V$ into $\bwedge^3 H(S)$. Since $\S_2$ has the totally separating partition, Theorem~\ref{thm:WS} states that $W_{\S_2}\iso \bwedge^3 H(\S_2)\oplus H(\S_2)$, where the second factor is spanned by $\delta_1(f)\otimes z_1$. Since $\iota_2(z_1)$ is homotopic to $C$, and $[C]=\omega_V\in N(S)$, Definition~\ref{def:WSiota} tells us that $(\iota_2)_*\colon W_{\S_2}\to W_S$ embeds the first factor by the inclusion $\bwedge^3 V^\perp \to \bwedge^3 H(S)$, and the second factor by $w\otimes z_1\mapsto w\wedge \omega_V$. Applying Theorem~\ref{thm:commutes}, we conclude that $\tau_S(\Stab_{\I(S)}(C))$ is the subspace spanned by $\tau_S(\I(\S_1))=(\iota_1)_*(W_{\S_1})=\bwedge^3 V$ and $\tau_S(\I(\S_2))=(\iota_2)_*(W_{\S_2})=\bwedge^3 V^\perp\oplus V^\perp\wedge \omega_V$, as desired.
\end{proof}

\begin{remark} More care must be taken when considering orbits of separating curves and multicurves when the partition on $\S$ is not totally separated. Consider the lantern $S=S_{0,4}$ depicted in Figure~\reffigurearcs, but with partition given by $\P=\{\{A_0,A_2\},\{A_1,A_3\}\}$. Here $A_1,A_2,A_3$ are the three boundary components in the center and $A_0$ is the outer boundary component. The two curves $C=\alpha$ and $D=\beta$ depicted in Figure~\reffigurearcs\ are both $\P$--separating. The complementary components $U_1$ and $U_2$ of $S-C$ span the subspaces $\langle a_1={-a_3}\rangle$ and $\langle a_2={-a_0}\rangle$ respectively of $H_S$, and the complementary components $V_1$ and $V_2$ of $S-D$ determine the same subspaces.

Nevertheless $C$ and $D$ do not lie in the same $\IS$--orbit, as can be seen by considering the larger surface $\Shat$. The complementary components of $\Surhat-D$ determine the subspaces $\langle a_1,b_3-b_1\rangle$ and $\langle a_2,b_2\rangle$ of $H(\S)$, while the complementary components of $\Surhat-C$ determine the subspaces $\langle a_1,b_3-b_1-a_2\rangle$ and $\langle a_2,b_2+a_1\rangle$. This shows there is no element of $\I(\Shat)$, and thus certainly no element of $\IS$, taking $C$ to $D$.
\end{remark}

\small

\noindent
Department of Mathematics\\ Stanford University\\
450 Serra Mall\\
Stanford, CA 94305\\
E-mail: church@math.stanford.edu


\begin{thebibliography}{ABCD}

\bibitem[BBM]{BBM}
M.~Bestvina, K.-U.~Bux, and D.~Margalit, \emph{The dimension of the Torelli group}, Jour. Amer. Math. Soc. \textbf{23} (2010) 1, 61--105. Available at arXiv:0709.0287.

\bibitem[B]{Bourbaki}
N.~Bourbaki, \emph{Lie groups and {L}ie algebras. {C}hapters 1--3}, Elements of
  Mathematics (Berlin), Springer-Verlag, Berlin, 1989, Translated from the
  French, Reprint of the 1975 edition.

\bibitem[BKS]{BKS}
R.~M. Bryant, L.~G. Kov{\'a}cs, and R.~St{\"o}hr, \emph{Subalgebras of free
  restricted {L}ie algebras}, Bull. Austral. Math. Soc. \textbf{72} (2005),
  no.~1, 147--156.

\bibitem[FLM]{FarbLeiningerMargalit}
B.~Farb, C.~Leininger, and D.~Margalit, \emph{The lower central series and pseudo-Anosov dilatations}, Amer. J. Math. \textbf{130} (2008), no. 3, 799--827. Available at arXiv:math/0603675.

\bibitem[HO]{Hahn-OMeara}
A.~Hahn and O.~T. O'Meara, \emph{The classical groups and {$K$}-theory},
  Grundlehren der Mathematischen Wissenschaften, vol. 291, Springer-Verlag,
  Berlin, 1989.

\bibitem[J1]{JohnsonAbelian}
D.~Johnson, \emph{An abelian quotient of the mapping class group
  {$\mathcal{I}_{g}$}}, Math. Ann. \textbf{249} (1980), no.~3, 225--242.

\bibitem[J2]{JohnsonConjugacy}
D.~Johnson, \emph{Conjugacy relations in subgroups of the mapping class group and
  a group-theoretic description of the {R}ochlin invariant}, Math. Ann.
  \textbf{249} (1980), no.~3, 243--263.

\bibitem[J3]{JohnsonI}
D.~Johnson, \emph{The structure of the {T}orelli group. {I}. {A} finite set of
  generators for {$\mathcal{I}$}}, Ann. of Math. \textbf{118} (1983), no.~3,
  423--442.

\bibitem[J4]{JohnsonII}
D.~Johnson, \emph{The structure of the {T}orelli group. {II}. {A} characterization
  of the group generated by twists on bounding curves}, Topology \textbf{24}
  (1985), no.~2, 113--126.

\bibitem[MP]{MeeksPatrusky}
W.~Meeks, III and J.~Patrusky, \emph{Representing homology classes by embedded
  circles on a compact surface}, Illinois J. Math. \textbf{22} (1978), no.~2,
  262--269.

\bibitem[M]{Morita}
S.~Morita, \emph{The extension of {J}ohnson's homomorphism from the {T}orelli group to the mapping class group}, Invent. Math. \textbf{111} (1993), no.~1, 197--224.

\bibitem[PR]{ParisRolfsen}
L.~Paris and D.~Rolfsen, \emph{Geometric subgroups of mapping class groups}, J.
  Reine Angew. Math. \textbf{521} (2000), 47--83. Available at arXiv:math/9906122.

\bibitem[P1]{CuttingPasting}
A.~Putman, \emph{Cutting and pasting in the {T}orelli group}, Geom. Topol.
  \textbf{11} (2007), 829--865. Available at arXiv:math/0608373.
  
\bibitem[P2]{PutmanKS}
A.~Putman, \emph{The Johnson homomorphism and its kernel}, preprint, arXiv:0904.0467v3, revised March 2013.

\bibitem[Sh]{Shirshov}
A.~I. Shirshov, \emph{Podalgebry svobodnykh lievykh algebr ({R}ussian:
  {S}ubalgebras of free {L}ie algebras)}, Mat. Sbornik N.S. \textbf{33(75)}
  (1953), 441--452.

\bibitem[W1]{WittD}
E.~Witt, \emph{Treue {D}arstellungen {L}iescher {R}inge ({G}erman: {F}aithful
  representations of {L}ie rings)}, J. Reine Angew. Math. \textbf{177} (1937),
  152--210.

\bibitem[W2]{WittU}
E.~Witt, \emph{Die {U}nterringe der freien {L}ieschen {R}inge ({G}erman:
  {S}ubrings of free {L}ie rings)}, Math. Z. \textbf{64} (1956), 195--216.

\end{thebibliography}
\end{document}